 \magnification=\magstep1 
\newcount\sectno
\newcount\subsectno
\newcount\parno
\newcount\equationno
\newif\ifsubsections
\subsectionsfalse

\def\sectnum{\the\sectno} 
\def\subsectnum{\sectnum\ifsubsections .\the\subsectno\fi} 
\def\parnum{\subsectnum .{\the\parno}}
\def\eqnum{\subsectnum .\the\equationno}

\def\abstract#1\endabstract
{
{\abstractfont
    \baselineskip=9pt
    \leftskip=4pc  \rightskip=4pc
    \bigskip
    \noindent
     ABSTRACT.\ #1
\medskip} 
}
 
\def\thanks[#1]#2\endthanks{\footnote{$^#1$}{\footnotefont\kern-6pt #2}}

\newcount\minutes
\newcount\scratch

\def\timestamp{%
\scratch=\time
\divide\scratch by 60
\edef\hours{\the\scratch}
\multiply\scratch by 60
\minutes=\time
\advance\minutes by -\scratch
\the \month/\the\day$\,$---$\,$\hours:\null
\ifnum\minutes< 10 0\fi
\the\minutes}

\def\today{\ifcase\month\or
January\or February\or March\or April\or May\or June\or
July\or August\or September\or October\or November\or December\fi
\space\number\day,\number\year}

\outer\def\newsection #1.\par{\vskip1.5pc plus.75pc \penalty-250
     \subsectno=0
     \parno=0
     \equationno=0
     \advance\sectno by1
     \leftline{\smalltitlefont \sectnum.\hskip 1pc  #1}
                \nobreak \vskip.75pc\noindent}

\outer\def\newsectiontwoline #1/#2/.\par{\vskip1.5pc plus.75pc \penalty-250
     \subsectno=0
     \parno=0
     \equationno=0
     \advance\sectno by1
     \leftline{\smalltitlefont \sectnum.\hskip 1pc  #1}
     \leftline{\smalltitlefont \hskip 22pt #2}
                \nobreak \vskip.75pc\noindent}

\outer\def\newsubsection #1.\par{\vskip1pc plus.5pc\penalty-250
     \parno=0
     \equationno=0
     \advance\subsectno by1
     \leftline{{\bf \subsectnum}\hskip 1pc  #1.}
                \nobreak \vskip.5pc\noindent}

\def\newpar #1.{\advance \parno by1
     \par
 \medbreak \noindent 
      {\bf \parnum. #1.} \hskip 6pt}

\long\def \newclaim #1. #2\par {\advance \parno by1
    \medbreak \noindent 
     {\bf \parnum \hskip 6 pt #1.\hskip 6pt} {\sl #2} \par \medbreak}

\def\eq $$#1$${\global \advance \equationno by1 $$#1\eqno(\eqnum)$$}

\def\rmarginsert[#1]{\hglue 0pt\vadjust
{\null\vskip -\baselineskip\rightline{\abstractfont\rlap{\hfil\  #1}}}}

\def\lmarginsert[#1]{\hglue 0pt\vadjust
{\null\vskip -\baselineskip\leftline{\abstractfont\llap{#1\ \hfill}}}}

\newif\ifproofmode
\proofmodefalse

\def\refpar[#1]#2.{\advance \parno by1
     \par
 \medbreak \noindent 
      {\bf \parnum \hskip 6 pt #2.\hskip 6pt}%
\expandafter\edef\csname ref#1\endcsname
{\parnum}\ifproofmode\rmarginsert[\string\ref#1]\fi}

\long\def \refclaim[#1]#2. #3\par {\advance \parno by1
    \medbreak \noindent 
{\bf \parnum \hskip 6 pt #2.\hskip 6pt}%
\expandafter\edef\csname ref#1\endcsname
{\parnum}\ifproofmode\rmarginsert[\string\ref#1]\fi
{\sl #3} \par \medbreak}

\def\refer[#1]{%
\expandafter\xdef\csname ref#1\endcsname
{\parnum}\ifproofmode\rmarginsert[\string\ref#1]\fi}

\def\refereq[#1]$$#2$$ {%
\eq$$#2$$%
\expandafter\xdef\csname ref#1\endcsname{(\eqnum)}%
\ifproofmode\rmarginsert[\string\ref#1]\fi
}

\def\refeq{\refereq}

\def \Definition #1\\ {\vskip 1pc \noindent 
      {\bf #1. Definition. \hskip 6pt}\vskip 1pc}

\def\proof{{PROOF.} \enspace}

\def\qedmark{\hbox{\vrule height 4pt width 3pt}}
\def\qedskip{\vrule height 4pt width 0pt depth 1pc}
\def\qed{\penalty 1000\quad\penalty 1000{\qedmark\qedskip}}

\def \i {\par\noindent \hskip 10pt}
\def \ii {\par\noindent \hskip 20pt}

\def \a {\alpha}
\def \b {\beta}

\def \d {\delta}
\def \D {\triangle}
\def \e {\epsilon}

\def \g {\gamma}

\def \K {\nabla}
\def \l {\lambda}

\def \n {\,\vert\,}
\def \N {\,\Vert\,}
\def \o {\theta}
\mathchardef\p="011E    

\def \s {\sigma}

\def \w {\omega}
\def \W {\Omega}

\def\II{I\!I}     

\def \reals {{\mbi R}}
\def \Cx {{\mbi C}}

\def\Gtwo{{\mathop{{{\mbi G\/}}\kern-.5pt_{{}_2}}}}
\def\Ffour{{\mathop{{{\mbi F\/}}\kern-2.5pt_{{}_4}}}}
\def\Esix{{\mathop{{{\mbi E\/}}\kern-.5pt_{{}_6}}}}
\def\Eseven{{\mathop{{{\mbi E\/}}\kern-.5pt_{{}_7}}}}
\def\Eeight{{\mathop{{{\mbi E\/}}\kern-.5pt_{{}_8}}}}


\def\ca{{\liefont A}}

\def\cf{{\liefont F}}

\def\cg{{\liefont G}}

\def\ck{{\liefont K}}
\def\cK{{\liefont K}}

\def\cm{{\liefont M}}
\def\cM{{\liefont M}}

\def\co{{\liefont O}}

\def\cp{{\liefont P}}
\def\cP{{\liefont P}}

\def\cu{{\liefont U}}
\def\cU{{\liefont U}}

\def\sdp{
\mathop{\hbox{$\raise 1pt\hbox{$\scriptscriptstyle |$}\kern-2.5pt\times$}}
        }

\def\dag{{\raise 1 pt\hbox{{$\scriptscriptstyle \dagger$}}}}

\def\*{\raise 1.5pt \hbox{*}}

\def\sech{\mathop{\rm sech}\nolimits}

\def\exp{\mathop{\rm exp}\nolimits}

\def\diag{\mathop{\rm diag}\nolimits}

\def\ad{\mathop{\rm ad}\nolimits}

\def\diag {\mathop{\rm diag}\nolimits}

\def \li{\langle}
\def \ri{\rangle}

                                                          
\def\ifundefined#1{\expandafter\ifx\csname#1\endcsname\relax}

\def\cross{\times}
\def\longerrightarrow{-\kern-5pt\longrightarrow}

\def\star{\lower 1pt\hbox{*}}
\def \nulset {
\raise 1pt\hbox{
\hskip -3pt$\not$\kern -0.2pt \raise .7pt\hbox{${\scriptstyle\bigcirc}$}}}
\def \norm|#1|{\Vert#1\Vert}

\def \interior(#1){#1\kern -6pt \raise 7.5pt 
      \hbox{$\scriptstyle \circ$}{}\hskip 2pt}

\def\twist_#1{\kern -.15em\cross\kern -.30 em{}_{{}_{#1}}\kern .07 em}

\font\cmr=cmr10 at 10pt
\font\cmrviii=cmr8
\font\cmrvi=cmr6

\font\cmrXIV=cmr12 at 14 pt
\font\cmrXX=cmr17 at 20 pt
\font\cmrXXIV=cmr17 at 24 pt
\font\cmbxXII=cmbx12
\font\cmbxsl=cmbxsl10
\font\cmbxslviii=cmbxsl10 at 8pt
\font\cmbxslv=cmbxsl10 at 5pt

      \font\tenrm=cmr10 at 10.3 pt
      \font\sevenrm=cmr7 at 7.21 pt
      \font\fiverm=cmr5 at 5.15 pt
      \font\teni=cmmi10 at 10.3 pt
      \font\seveni=cmmi7 at 7.21 pt
      \font\fivei=cmmi5 at 5.15 pt     
      \font\tensy=cmsy10 at 10.3 pt
      \font\sevensy=cmsy7 at 7.21 pt   
      \font\fivesy=cmsy5 at 5.15 pt    
      \font\tenex=cmex10 at 10.3 pt
      \font\tenbf=cmbx10 at 10.3 pt
      \font\sevenbf=cmbx7 at 7.21 pt    
      \font\fivebf=cmbx5 at 5.15 pt  

\def\UseComputerModern   
{
\textfont0=\tenrm \scriptfont0=\sevenrm \scriptscriptfont0=\fiverm
\def\rm{\fam0\tenrm}
\textfont1=\teni \scriptfont1=\seveni \scriptscriptfont1=\fivei
\def\mit{\fam1} \def\oldstyle{\fam1\teni}
\textfont2=\tensy \scriptfont2=\sevensy \scriptscriptfont2=\fivesy
\def\cal{\fam2}
\textfont3=\tenex \scriptfont3=\tenex \scriptscriptfont3=\tenex
\def\it{\fam\itfam\tenit} 
\textfont\itfam=\tenit
\def\sl{\fam\slfam\tensl} 
\textfont\slfam=\tensl
\def\bf{\fam\bffam\tenbf} 
\textfont\bffam=\tenbf \scriptfont\bffam=\sevenbf
\scriptscriptfont\bffam=\fivebf
\def\tt{\fam\ttfam\tentt} 
\textfont\ttfam=\tentt
\def\abstractfont{\cmrviii}
\def\footnotefont{\cmrviii}
\def\tinyfont{\cmrvi}
\def\smalltitlefont{\cmbxXII}  
\def\titlefont{\cmrXIV}
\def\bigtitlefont{\cmrXX}
\def\verybigtitlefont{\cmrXXIV} 
\textfont9=\cmbxsl \scriptfont9=\cmbxslviii \scriptscriptfont9=\cmbxslv
\def\mbi{\fam9}
\cmr
}  

\def\liefont{\cal}

\def \bs {\bigskip}
\def \ms {\medskip}
\def \ss {\smallskip}

\def \ni {\noindent}

\def\enditem{\item{}\par\vskip-\baselineskip\noindent}
\def\ei{\enditem}

   \baselineskip=14 true pt 
   \hsize 37 true pc \hoffset= 22 true pt
   \voffset= 0 true pt
   \vsize  54 true pc

\def\ni{\noindent}
\def\bs{\bigskip}
\def\ms{\medskip}

\def\ss{\smallskip}

   \baselineskip=14 true pt 
   \hsize 35 true pc \hoffset= 25 true pt
   \vsize  52 true pc
\UseComputerModern
\subsectionsfalse
\font\smalltitlefont=cmbx10 at 11 pt
\proofmodetrue
%


\def\Bibliography
{

\font\TRten=cmr10 at 10 true pt
\font\TIten=cmti10 at 10 true pt
\font\TBten=cmbx10 at 10 true pt

\def\ourindent{\hfil\vskip-\baselineskip}

 \frenchspacing
 \parindent=0pt

 \def \keyfnt{\TRten}
 \def \authornamefnt{\TRten}
 \def \booktitlefnt{\TIten}
 \def \articletitlefnt{\TRten}
 \def \journalnamefnt{\TIten}
 \def \volumefnt{\TBten}
 \def \publishernamefnt{\TRten}
 \def \pagesfnt{\TRten}
 \def \yearfnt{\TRten}
 \def \commentfnt{\TRten}

 \def \bookitem //##1//##2//##3//##4//##5//##6//##7//##8//
      { \goodbreak{\par\hskip-40pt{\keyfnt [##1]}\ourindent{\authornamefnt ##2,}}
              {\booktitlefnt ##3.\/}\thinspace
              {\publishernamefnt ##4,}
              {\yearfnt ##6.}
              {\commentfnt ##8}
       }

\def \b{\bookitem}

\def \articleitem //##1//##2//##3//##4//##5//##6//##7//##8//
      { \goodbreak{\par\hskip-40pt{\keyfnt [##1]}\ourindent{\authornamefnt ##2,}}
              {\articletitlefnt ##3},
              {\journalnamefnt ##4\/}
              {\volumefnt ##5}
              {\hbox{\yearfnt(\hskip -1pt ##6)}},
              {\pagesfnt ##7.}
              {\commentfnt ##8}
       }
\def \a{\articleitem}

\def \preprintitem //##1//##2//##3//##4//##5//##6//##7//##8//
      { \goodbreak{\par\hskip-40pt{\keyfnt [##1]}\ourindent{\authornamefnt ##2,}}
              {\articletitlefnt ##3},
              {\commentfnt ##8}
       }
\def \p{\preprintitem}

   \vskip 1in
   \centerline{References}
   \vskip .5in
}


\def\bu{\bullet}

\def\Reduce#1  
{
\font\stenrm=cmr10 scaled #1
\font\sninerm=cmr9 scaled #1
\font\seightrm=cmr8 scaled #1
\font\ssevenrm=cmr7 scaled #1
\font\ssixrm=cmr6 scaled #1
\font\sfiverm=cmr5 scaled #1

\font\steni=cmmi10 scaled #1
\font\sninei=cmmi9 scaled #1
\font\seighti=cmmi8 scaled #1
\font\sseveni=cmmi7 scaled #1
\font\ssixi=cmmi6 scaled #1
\font\sfivei=cmmi5 scaled #1

\font\stenit=cmti10 scaled #1
\font\snineit=cmti9 scaled #1
\font\seightit=cmti8 scaled #1
\font\ssevenit=cmti7 scaled #1

\font\stensy=cmsy10 scaled #1
\font\sninesy=cmsy9 scaled #1
\font\seightsy=cmsy8 scaled #1
\font\ssevensy=cmsy7 scaled #1
\font\ssixsy=cmsy6 scaled #1
\font\sfivesy=cmsy5 scaled #1

\font\stenbf=cmbx10 scaled #1 
\font\sninebf=cmbx9 scaled #1
\font\seightbf=cmbx8 scaled #1
\font\ssevenbf=cmbx7 scaled #1
\font\ssixbf=cmbx6 scaled #1
\font\sfivebf=cmbx5 scaled #1

\font\stentt=cmtt10  scaled #1
\font\sninett=cmtt9 scaled #1
\font\seighttt=cmtt8 scaled #1

\font\stenex=cmex10 scaled #1

\font\stensl=cmsl10 scaled #1
\font\sninesl=cmsl9 scaled #1
\font\seightsl=cmsl8 scaled #1

\textfont0=\stenrm \scriptfont0=\ssevenrm \scriptscriptfont0=\sfiverm
\def\rm{\fam0\stenrm}
\textfont1=\steni \scriptfont1=\sseveni \scriptscriptfont1=\sfivei
\def\mit{\fam1} \def\oldstyle{\fam1\steni}
\textfont2=\stensy \scriptfont2=\ssevensy \scriptscriptfont2=\sfivesy
\def\cal{\fam2}
\textfont3=\stenex \scriptfont3=\stenex \scriptscriptfont3=\stenex
\def\it{\fam\itfam\tenit} 
\textfont\itfam=\stenit
\def\sl{\fam\slfam\stensl} 
\textfont\slfam=\stensl
\def\bf{\fam\bffam\stenbf} 
\textfont\bffam=\stenbf \scriptfont\bffam=\ssevenbf
\scriptscriptfont\bffam=\sfivebf
\def\tt{\fam\ttfam\stentt} 
\textfont\ttfam=\stentt
\rm
}  


\input BoxedEPS
\SetRokickiEPSFSpecial
\HideDisplacementBoxes

\nopagenumbers \headline{\tenrm \hfill \folio}
\def\na{\natural}
\def\ti{\tilde}
\def\bu{$\bullet$\hskip 4pt}

\proofmodefalse


\def\ti{\tilde}

\def\p{\partial}

\font\twelvebf=cmbx12

\centerline {\twelvebf The Submanifold Geometries associated to  Grassmannian
Systems\/}

\bs\bs \centerline{Martina Br\"uck, Xi Du, Joonsang Park, and
Chuu-Lian Terng\footnote{$^*$} {Research supported in part by 
NSF Grant DMS 9626130}\/} 
\bs \bs

\centerline{\bf Abstract}
\bs

{\Reduce{900}
There is  a hierarchy of commuting soliton
equations associated to each symmetric space $U/K$. When $U/K$
has rank $n$, the first $n$ flows in the hierarchy give rise to a natural first order
non-linear system of partial differential equations in n variables, the so called 
$U$/K-system. Let $G_{m,n}$ denote the Grassmannian of $n$-dimensional linear
subspaces in $R^{m+n}$, and $G_{m,n}^1$ the Grassmannian of space-like
$m$-dimensional linear subspaces in the Lorentzian space $R^{m+n,1}$.   In this
paper, we use techniques from soliton theory to study submanifolds in space forms
whose Gauss-Codazzi equations are gauge equivalent to the $G_{m,n}$-system or
the $G_{m,n}^1$-system. These include submanifolds with constant sectional
curvatures, isothermic surfaces, and submanifolds admitting principal curvature
coordinates. The dressing actions of simple elements on the space of solutions of the
$G_{m,n}$ and  $G_{m,n}^1$
systems correspond to B\"acklund, Darboux and Ribaucour transformations for 
submanifolds. 
}
\bs\bs

\centerline {\bf Table of Contents\/}
\ms

\def\dotting{\leaders\hbox to  1em{\hfil. \hfil}\hfil}
\line{1. Introduction \dotting \ \ 2}
\line{2. The $U/K$-system \dotting \ \ 8}
\line{3. $G_{m,n}$-systems \dotting 13}
\line{4.  $G_{m,n}^1$-systems \dotting 20}
\line{5. The moving frame method for submanifolds \dotting 23}
\line{6. Submanifolds associated to  $G_{m,n}$-systems \dotting 27}
\line{7. Submanifolds associated to  $G_{m,n}^1$-systems \dotting 37}
\line{8. $G_{m,1}^1$-system and isothermic surfaces \dotting 42}
\line{9.  Loop group action for $G_{m,n}$-systems\dotting 50}
\line{10. Ribaucour transformations for $G_{m,n}$-systems \dotting 59}
\line{11. Loop group action for  $G_{m,n}^1$-systems  \dotting 69}
\line{12. Ribaucour transformations for $G_{m,n}^1$-systems\dotting 73}
\line{13. Darboux transformations for $G_{m,1}^1$-systems\dotting 75}
\line{14. B\"acklund transformations \dotting 80}
\line{15. Permutability formula for Ribaucour transformations\dotting 88}
\line{16. The $U/K$-hierarchy and finite type solutions \dotting 91}

\bs

\newsection  Introduction.\par

Some of the high points in classical differential geometry are the study of surfaces in
$R^3$ with special geometric properties, to find good coordinates so that the
corresponding Gauss-Codazzi equations have specially nice forms, and to construct
explicit examples and deformations of these surfaces.   Surfaces with
negative constant Gaussian curvature, surfaces with constant mean curvature, and
isothermic surfaces are some of the well-known examples.  The Gauss and Codazzi
equations of these surfaces are now known to be ``soliton'' equations.  In recent
years, modern geometers have found that these equations admit ``Lax pairs'', i.e.,
they can be written as the condition of a family of connections to be flat.  The
existence of a Lax pair is one of the characteristic properties of soliton equations,
and it often gives rise to an action of an infinite dimensional group on the space of
solutions (the dressing action).  The geometric transformations found for these
surfaces by classical geometers such as B\"acklund, Darboux, and
Ribaucour transformations, often arise as the dressing action of some simple rational
elements.  For more detail see [Bo2], [Bo3], [TU3].  In this approach, we start with a
class of surfaces in $R^3$.  If there are methods to construct an infinite parameter
family of solutions from a given one, then it hints that we may be able to find a good
coordinate system and a Lax pair. Geometers have used this method to construct
soliton equations involving n variables (cf. [TT], [Te1], [FP2]).  But there is no
uniform algorithm to achieve this or determine whether a geometric equation for a
certain class of submanifolds is a soliton equation.   

It is also known that we can associate to each symmetric space $U/K$ a hierarchy of
soliton equations (cf. [TU1]).  For example, the $SU(2)$-hierarchy is the hierarchy
for the non-linear Schr\"odinger equation and the $SU(2)/SO(2)$-hierarchy is the
hierarchy for the modified KdV equation.  If the rank of the symmetric space $U/K$
is $n$, Terng [Te2] put the $n$ first flows together to construct a natural non-linear
first order system, the $U/K$-system, and initiated the project of identifying the
submanifold geometry associated to these systems.  This
means to find submanifolds in certain symmetric space
$M$ whose Gauss-Codazzi equation is given by the
$U/K$-system  and  to find the geometric
transformations corresponding to the dressing actions of certain simple elements. 
This direct approach may provide ways to find Lax pairs for some known class
of submanifolds, and also may give new interesting class of submanifolds.  The main
goal of this paper is to carry out this project for the real Grassmannian manifolds of
space-like $m$-dimensional linear subspaces in
$R^{m+n}$ and in $R^{m+n,1}$.    

Below we give a short review of some known facts and
outline our results.  

\ms\ni \bu {\bf The $U/K$-system}
\ss

Let $U$ be a semi-simple Lie group, $\s$ an involution on $G$, and $K$ the fixed point
set of $\s$.  Then $U/K$ is a symmetric space.  The Lie algebra $\ck$ is the
$+1$ eigenspace of the differential $\s_\ast$ of $\s$ at the identity.  Let
$\cp$ denote the $-1$ eigenspace of $\s_\ast$.  Then 
$\cu= \ck \oplus \cp$ and 
$$[\ck,\ck]\subset \ck, \quad [\ck, \cp]\subset \cp, \quad [\cp, \cp]\subset \cp.$$
  Let $\ca$ be a maximal abelian subalgebra in $\cp$, $a_1,
\dots , a_n$ a basis for $\ca$, and $\ca^\perp$ the orthogonal complement of
$\ca$ in $\cU$ with respect to the Killing form $< , >$. 
We recall that $n=\dim(\ca)$ is called the {\it rank\/} of
the symmetric space.  {\it The n-dimensional
system associated to $U/K$\/},  defined by Terng in [Te2], is the following first
order non-linear partial differential equation  for $v: R^n\to \cP\cap\ca^\perp $:
\refeq[at]$$[a_i, v_{x_j}]-[a_j, v_{x_i}]=\big[\,[a_i, v], [a_j, v]\,\big], \ \ \ \
1\le
i\ne j \le n,$$
where $v_{x_i}= {\partial v\over\partial x_i}$.   We will call such system
the $U/K$-{\it system\/}.  

\ms\ni \bu {\bf The Lax connection}
\ss

Recall that a $\cg$-valued connection ${\p\over \p
x_i} + A_i$ is flat if its curvature is zero, i.e.,
$$\left[{\p\over \p x_i} + A_i, \ \ {\p\over \p x_j} + A_j\right] =0$$ 
for all $i,j$.  
Let $V$ be a linear subspace of a Lie algebra $\cg$. 
A partial differential equation (PDE) for $v:R^n\to V$ admits a {\it Lax connection\/}
if there exists a family of $\cg$-valued connection 
$${\p\over \p x_i} + A_i(v, dv, d^2v, \cdots, d^k v, \l)$$ such that the
PDE for
$v$ is given by the flatness of these connections for all $\l$ in some open domain in
$C$.  Equation
\refat{} admits a Lax connection, ${\p\over \p x_i} + a_i\l +  [a_i, v]$. In other words,
$v$ is a solution of \refat{} if and only if 
\refeq[au]$$\left[{\p\over \p x_i} + a_i\l + [a_i,v], \ \ {\p\over \p x_j} + a_j \l +
[a_j,v] \right]=0$$ for all $i,j$ and $\l\in C$.  When $n=2$, a Lax connection gives rise
to a pair of commuting operators. This was first observed by Lax for the KdV
equation and is called a {\it Lax pair\/} in the soliton literature. 

Note that the connection ${\p\over \p x_i} + A_i$ is flat if and only if the connection
1-form $\w=\sum_i A_i dx_i$ is flat, i.e., $d\w= -\w\wedge \w$.  In particular, $v$ is a
solution of the $U/K$-system \refat{} if and only if 
\refeq[au]$$\o_\l = \sum_{i=1}^n (a_i\l + [a_i, v]) dx_i$$ is flat for all $\l\in C$. 
We will also call $\o_\l$ the Lax connection of the $U/K$-system.

\ms\ni \bu {\bf Dressing action\/}\ss

The existence of a Lax connection for an equation often gives rise to an action of
certain subgroup of germs of holomorphic maps at a suitable point  on the space of
local solutions of the equation.  This is called ``dressing action" in the soliton
literature.  Below we give a rough sketch of the construction of the dressing action
(cf. [ZS], [Ch], [TU2]).  

 If $v$ is a solution of
\refat{}, then its Lax connection
$\o_\l$ is flat for all $\l\in C$, so there exists a $U_C$-valued map $E(x,\l)$ such that 
\refeq[ir]$$E^{-1}dE= \o_\l, \quad E(0,\l)=I.$$
Since $\o_\l$ is holomorphic for $\l\in C$, so is $E(x,\l)$.
Now let $g(\l)$ be a holomorphic map defined from a neighborhood of $\l=\infty$ in
$S^2=C\cup \{\infty\}$ to $U_C$ that satisfies certain $U/K$-reality condition
(defined later) and  
$g(\infty)=I$. It follows from the classical Birkhoff factorization theorem (cf. [PS])
that 
 there exist uniquely
$\tilde E(x,\l)$ and $\tilde g(x,\l)$ so that 
\refeq[jg]$$g(\l) E(x,\l) = \tilde E(x,\l) \tilde g(x,\l),$$
 $\tilde E(x,\l)$ is holomorphic for $\l\in C$, $\tilde g(x,\l)$ is holomorphic
near $\l=\infty$ and $\tilde g(x,\infty)=I$. Calculate the residue at $\l=\infty$ to
conclude that 
$$\tilde E^{-1} d\tilde E= \sum_{i=1}^n (a_i\l + [a_i, \tilde v])dx_i$$ for some
$\tilde v$.  So $\tilde v$ is a new solution of \refat{}.  The solution
$\tilde v$ can also be obtained from
$\tilde g$ as follows: Expand 
$$\tilde g(x,\l)\ = \ I + m_1(x) \l^{-1} + m_2(x)\l^{-2} + \cdots$$
at $\l=\infty$.  Then 
$$\tilde v= v - p_0(m_1),$$
where $p_0$ is the projection onto $\ca^\perp\cap \cp$.  
The map 
$$g\sharp v= \tilde v$$
defines the dressing action of the group of germs of holomorphic maps on the space of
local solutions of the $U/K$-system.  

If $g(\l)$ is a meromorphic map on $S^2$ with $g(\infty)=I$, then the factorization 
\refjg{} can be done explicitly by  calculating the residues at poles of $g(\l)$.  
Note that system \refat{} has
a trivial solution $v=0$ and 
$$E(x,\l) = \exp\left(\sum_{i=1}^n a_i x_i\l\right)$$
is the solution of the corresponding linear system \refir{}.  Therefore $g\sharp 0$ can
be computed explicitly.  These explicit solutions correspond to the ``{\it pure
solitons\/}'' in the theory of soliton equations. 

If $g(\l)$ is a holomorphic map defined in a neighborhood of $\infty$ in $S^2$ such
that  $g(\l)a_1g(\l)^{-1}$ is a polynomial in $\l^{-1}$, then the solution $g\sharp 0$
can  be obtained by solving a system of ordinary differential equations on a finite
dimensional linear space.  These solutions are the so called ``{\it finite type
solutions\/}''.  Finite type solutions have been used successfully to construct constant
mean curvature tori in $R^3(c)$ by Pinkall and Sterling in [PiS], in 3-dimensional
space forms by Bobenko in [Bo1], and harmonic maps from a torus to a symmetric
space by Burstall, Ferus, Pedit and Pinkall in [BFPP].  

\ms\ni \bu {\bf Cauchy problem\/}\ss

If $a_1$ is regular, then the linear map $\ad(a_1):\cp\cap \ca^\perp\to
\ck$ is injective.  It follows from Cartan-K\"ahler Theorem
 that if $v_0:(-\d, \d)\to \cp\cap\ca^\perp$ is real analytic, then system
\refat{} has a unique local analytic solution such that $v(x_1, 0, \cdots,
0)=v_0(x_1)$.  If the initial data $v_0$ is not real analytic but is rapidly decaying, then
we can use the inverse scattering method developed by Beals and Coifman [BC] to
solve the initial value problem (cf. [TU1], [Te2]).  

\ms\ni \bu {\bf Gauge equivalent systems\/}\ss

Let $\o_\l$ be the flat connection \refau{} associated to the solution $v$ of the
$U/K$-system
\refat{}, and  $g:R^n\to U_C$ a smooth map. Then the gauge transformation 
$$g\ast \o_\l= g\o_\l g^{-1} - dg g^{-1}$$
is again a flat connection 1-form for all $\l\in C$.  However, the differential equation
given by the condition that $g\ast \o_\l$ is flat for all $\l$ has a different
form.  We say this new equation is gauge equivalent to system \refat{}. 
For example: 
\item {(i)} Since $\o_\l$ satisfies the $U/K$-reality condition, $\o_0$ is a
$\ck$-valued flat connection 1-form.  Hence there exists $g$ such that
$g^{-1}dg=\o_0$. A direct computation shows that the gauge transformation of
$\o_\l$ by $g$ is
$$g\ast \o_\l= g\o_\l g^{-1} - dg g^{-1} = \sum_{i=1}^n ga_i g^{-1}\l dx_i.$$
Write $A_i= ga_ig^{-1}$.  The equation given by the flatness of $g\ast \o_\l$ is the
{\it curved flat system\/} studied by Ferus and Pedit in [FP1].  
\item {(ii)} Suppose $K= K_1\times K_2$.  Let $v$ be a solution of the
$U/K$-system, and $g^{-1}dg = \o_0$.  Since $\o_0\in \ck$,  $g(x)\in
K=K_1\times K_2$.  So we can write
$g=(g_1,g_2)\in K_1\times K_2$.  A direct computation shows that the
coefficients of $\l^0$ in $g_1\ast \o_\l$ and $g_2\ast \o_\l$ are in $\ck_2$ and
$\ck_1$ respectively.  The equations given by the flatness of $g_1\ast\o_\l$ and
$g_2\ast\o_\l$ are called the $U/K$-system I and II respectively.

\ms\ni \bu {\bf Gauss-Codazzi equations for Submanifolds in space forms\/}
\ss 

Let $O(n,1)$ denote the group of all $g\in GL(n+1)$ that preserves the bilinear form
$$x_1^2+ \cdots + x_n^2 - x_{n+1}^2.$$
Henceforth in this paper, we use the following notations: 
$$G_{m,n}= O(m+n)/O(m)\times O(n), \quad G_{m,n}^1= O(m+n,1)/O(m)\times
O(n,1).$$
 
Let $N^n(c)$ denote the $n$-dimensional space form with curvature $c$, i.e., the
complete, simply connected Riemannian manifold with constant sectional curvature
$c$. So $N^n(c)$ is $R^n, S^n$ and $H^n$ for $c=0, 1, -1$ respectively.  The
Levi-Civita connection
$1$-form of
$N^n(c)$ can be read from the flatness of a
$so(n)$, $so(n+1)$ and $o(n,1)$- valued  connection 1-form. The Gauss-Codazzi
equation of a submanifold in $N^n(c)$ is given by the flatness of the restriction of this
connection 1-form to the submanifold.   The Fundamental Theorem of Submanifolds
states that each solution of the Gauss-Codazzi equation correspond to a submanifold in
$N^n(c)$, unique up to ambient isometry. So if $v$ is  a solution of the $G_{m,n}$- or
$G^1_{m,n}$-system I or II, then the corresponding  Lax connection $\o_\l$ at
$\l=1$  gives rise to a submanifold of a certain space form.  Using the method of moving
frames,  special properties of the flat connection
$\o_1$ can be translated easily to geometric properties of the corresponding
submanifolds.  

\ms\ni \bu {\bf Submanifolds corresponding to the $G_{m,n}$- and
$G_{m,n}^1$-system I\/}
\ss

In [Te2], Terng proved that solutions of the $G_{n,n}$-, $G_{n,n+1}$- and
$G_{n,n}^1$-system I correspond to local isometric immersions of the space form
$N^n(c)$ in
$N^{2n}(c)$ with flat normal bundle for $c=0, 1$, and $-1$ respectively. 
We generalize this result to the $G_{m,n}$- and $G_{m,n}^1$-system I for any
$m\geq n$. They give rise to local isometric immersions of $N^n(c)$ into
$N^{n+m}(c)$. 

\ms\ni \bu {\bf Submanifolds corresponding to $G_{m,n}$- and $G_{m,n}^1$-system
II\/}
\ss

In order to explain the submanifold geometry corresponding to the
$G_{m,n}$- and $G_{m,n}^1$-system II, we first need to review some classical
surface theory.  Let $M$ be a surface in $R^3$ with curvature $K= -1$, and $e_3$ its
unit normal field. Then there exists a line of curvature coordinate system  $(x,y)$
such that the two fundamental forms of $M$ are 
$$I_1=\cos^2 u \ dx^2 + \sin^2 u\ dy^2,\quad \II_1=\sin u\cos
u\ (dx^2-dy^2).$$
Since the unit sphere is totally umbilic, the fundamental forms for $S^2$
in $(x,y)$ coordinates via the parametrization $e_3(x,y)$ are
$$I_2 = \sin^2 u \ dx^2 +\cos^2 u\ dy^2, \quad  \II_2= -(\sin^2 u \ dx^2
+\cos^2 u \ dy^2).$$  The Gauss-Codazzi equations
for $M$ (curvature
$-1$) and
$e_3$ (curvature $1$) are the same sine-Gordon equation
$$u_{xx}-u_{yy} = \sin u\cos u, \eqno({\rm SGE})$$ and the tangent plane of $M$ at
$(x,y)$ is the same as the tangent plane of the sphere at $e_3(x,y)$.  A direct
computation shows that such
$u$ gives rise to a solution of the
$G_{3,2}$-system II. This is a special case. 
In fact, we show that each solution of the
$G_{3,2}$-system II corresponds to a pair of surfaces $(X_1, X_2)$ in $R^3$ with
common line of curvature coordinates $x,y$ such that the tangent plane at $X_1(x,y)$
is equal to the tangent plane at $X_2(x,y)$ and their fundamental forms are
$$\eqalign{&I_1= \cos^2 u \ dx^2 + \sin^2 u \ dy^2, \quad \II_1= g_1\cos u\ dx^2 +
g_2
\sin u\ dy^2,\cr & I_2= \sin^2 u\ dx^2 + \cos^2 u\ dy^2, \quad \II_2= -g_1\sin u
\ dx^2 + g_2
\cos u\ dy^2\cr}$$ for some functions $u, g_1, g_2$.  
Moreover, the Gaussian curvature 
$$K_1(x,y)= -K_2(x,y).$$ 
For general $m>n$, we prove that each 
solution of the $G_{m,n}$-system II ($G_{m,n-1}^1$-system II respectively) gives rise
to an n-tuple of n-dimensional submanifolds  $(X_1, \cdots, X_n)$ in
$R^m$ with flat normal bundles and  common line of curvature coordinates $(x_1,
\cdots, x_n)$ such that fundamental forms of $X_j$ are 
$$I_j=\sum_{i=1}^n a_{ji}^2(x) dx_i^2, \quad \II_j= \sum_{i=1, k=1}^{n,m-n}
a_{ji}g_{ki} \ dx_i^2 e_{n+k}$$ for some
$(a_{ij}(x))\in O(n)$ ($\in O(n-1,1)$ respectively) and $g_{ki}(x)$.    

\ss

If we use another standard form of $O(n+1,1)$, the group of $g\in GL(n+2)$ that
leaves the bi-linear from 
$$x_1^2+ \cdots + x_n^2 + 2 x_{n+1} x_{n+2}$$
invariant, 
then the corresponding 2-tuple $(X_1, X_2)$ in $R^n$ for the $G_{n,1}^1$-system II
has the property that
$X_1$ is an isothermic surface and $X_2$ is a Christoffel dual of $X_1$. This was
proved by Burstall, Hertrich-Jeromin, Pedit and Pinkall in [BHPP] for $n=2$ and by
Burstall in [Bu] for general $n$. 

\ms\ni \bu {\bf B\"acklund transformations and dressing action\/}
\ss

 Let $M, M^*$ be two surfaces in $R^3$. A diffeomorphism $\ell:M\to M^*$ is called a
{\it B\"acklund transformation\/} with constant $\o$ if for all $p\in M$, 
\item {(a)} $\overline{pp^*}$ is tangent to both $M$ and $M^*$ at $p$ and
$p^*=\ell(p)$,
\item {(b)} d$(p,p^*)=\sin \o$,
\item {(c)} the angle between $TM_p$ and $TM^*_{p^*}$ is $\o$.

\ni B\"acklund proved ([Ba]) that if $\ell$ is a B\"acklund transformation, then both
$M$ and $M^*$ have curvature $-1$. Moreover, if $M$ is a surface in $R^3$ with
$K=-1$,
$0<\o<\pi$ a constant, and
$v_0\in TM_{p_0}$ a unit vector that is not a principal direction, then there exist a
unique surface $M^*$ and a B\"acklund transformation
$\ell:M\to M^*$ such that $\ell(p_0)= p_0+ \sin\o \  v_0$. 
Analytically, this gives a method of constructing new solution of SGE from a given one.  
 B\"acklund transformations have been generalized to
isometric immersions of $N^n(c)$ in $N^{2n-1}(c+1)$ by Terng and  Tenenblat for
$c=-1$ in [TT] and by Tenenblat for $c=0, 1$ in [Ten].  

 Terng and Uhlenbeck proved in [TU2]  that the dressing
action of a meromorphic map with one pole on the space of solutions of SGE
gives rise exactly to the classical B\"acklund transformations. We generalize this result
to  $G_{n,n}$ and $G_{n,n}^1$-systems.  

\ms\ni \bu {\bf Ribaucour transformations and dressing action\/}
\ss

Let  $M, \tilde M$ be two surfaces in $R^3$.  A diffeomorphism $\ell:M\to \tilde M$
is called a {\it Ribaucour transformation\/} if  for all $p\in M$
\item {(a)} $TM_p=T\tilde M_{\ell(p)}$,
\item {(b)} the normal line at $p$ to $M$ meets the normal line at $\ell(p)$ to $\tilde
M$ at equidistance $r(p)$,
\item {(c)} the line through $p$ in the principal direction $e$ of $M$ meets the line
through $\ell(p)$ in the direction $\ell_\ast(e)$ at a point at equidistance $s(p)$.  

\ni The notion of Ribaucour transformations has a natural generalization to
submanifolds in space forms with flat normal bundle ([DT]).  We show that the
dressing action of rational maps with two simple poles on the solutions of the
$G_{m,n}$- and $G_{m,n}^1$-system I and II correspond to Ribaucour
transformations for submanifolds.  

\ms\ni \bu {\bf Organization of the paper\/}
\ss

 We review some general facts about the $U/K$-system
in section 2, write down the $G_{m,n}$-systems explicitly in section 3, and the
$G_{m,n}^1$-systems in section 4.  We review the method of moving frames in section
5. We describe submanifolds associated to various $G_{m,n}$-systems and
$G_{m,n}^1$-systems in section 6 and 7 respectively.  In section 8, we  study relations
between constant mean curvature in 3-dimension space forms, isothermic surfaces and
$G_{m,1}^1$-systems.  The dressing action of a rational map with two simple poles on
solutions of the
$G_{m,n}$- and $G_{m,n}^1$-systems are written down explicitly in section 9 and 11
respectively. The corresponding geometric transformations are given in section 10 and
12.  Burstall [Bu] gave a generalization of isothermic surfaces in $R^3$ and their
Darboux transformations to isothermic surfaces  in $R^n$.  In section 13,  we show
that  the dressing action of a rational map with two poles on the space of solutions
of the
$G_{n,1}^1$-system II gives rise to these Darboux transformations.   In section 14,
we give a relation between
the dressing action of loop with one simple pole and B\"acklund transformations. 
A permutability formula for Ribaucour transformations  is explained in section 15, and
relation between dressing action and finite type solutions of the $U/K$-system are
given in the last section. 

\bs

\newsection The $U/K$-System.\par

A connection of a trivial principal $U$-bundle over $R^n$ is:
$${\p\over \p x_i} + A_i, \quad 1\leq i\leq n,$$ for some smooth
maps $A_i:R^n\to \cu$. The curvature of this connection is
$$\W_{ij}= \left[{\p\over \p x_i}+ A_i, {\p\over \p x_j}+
A_j\right]= (A_j)_{x_i}-(A_i)_{x_j} + [A_i, A_j].$$
A connection is flat if its curvature is zero.  

The following Proposition, which is well-known, gives several equivalent conditions for
a connection to be flat.  The proof
follows from a direct computation.

\refclaim[bc] Proposition.  Let $A_1, \cdots, A_n:R^n\to \cU$ be
smooth maps. The following statements are equivalent:
\item {(1)} the connection ${\p\over \p x_i} + A_i(x)$ is flat, i.e.,
$[{\partial\over
\partial x_i} +A_i, {\partial\over
\partial x_j} +A_j]=0$,
\item {(2)} $(A_j)_{x_i}- (A_i)_{x_j} + [A_i, A_j]=0$,
\item {(3)} the connection 1-form $\o=\sum\limits_{i=1}^n A_idx_i$
is flat, i.e., $d\theta + \theta\wedge \theta =0$.
\item {(4)}
\refeq[bn]$$E_{x_i}=EA_i, \quad 1\leq i\leq n,$$ is solvable for
$E:R^n\to U$.\ei

\ni
$E$  is called a {\it trivialization\/} of the flat
connection $\o=\sum_i A_i dx_i$ if it is a solution of \refbn{} or equivalently if
$$E^{-1}dE=\sum_{i=1}^n A_i dx_i.$$ \ss

It follows from a direct computation and Proposition \refbc{} that

\refclaim[as] Proposition ([Te2]). The following statements are
equivalent for a map $v:R^n\to \cp\cap \ca^\perp$:
\item {(i)}  $v$ is a solution of the $U/K$-system \refat{},
\item {(ii)}
\refeq[ba]$$\bigg[{\partial\over\partial x_i}+\lambda a_i +[a_i,
v],\,\, {\partial\over\partial x_j}+\lambda a_j +[a_j, v]\bigg]=0,
\ \ \ \forall\,\, \lambda \in C,$$
\item {(iii)} $\o_\l$ is a flat $\cu_C=\cu\otimes C$-connection 1-form on $R^n$
for all $\l\in C$, where \refeq[au]$$\o_\l=\sum_{i=1}^n (a_i\l +
[a_i, v]) dx_i,$$
\item {(iv)} there exists $E$ so that $E^{-1}dE= \o_\l$.

The one parameter family of flat connections \refba{} or \refau{}
is called the {\it Lax connection\/} of the $U/K$-system \refat{}.

\ms

An element $b\in \cp$ is called {\it regular\/} if the orbit at
$b$ for the Ad$(K)$-action is a principal orbit. Let $\ca$ be a
maximal abelian subalgebra in $\cp$, $$\ck_\ca=  \{\xi\in \ck\n
[\xi, a]=0 \quad \forall\,\, a\in \ca\},$$ and $\ck_\ca^\perp$ the
orthogonal complement of $\ck_\ca$ in $\ck$.  It follows from
standard theory of symmetric space (cf. [H]) that if $b\in \ca$ is
regular, then $\ad(b)$ maps $\ck_\ca^\perp$ and $\cp\cap\ca^\perp$
isomorphically to $\cp\cap\ca^\perp$ and $\ck_\ca^\perp$
respectively.

\refclaim[bd] Proposition.  Let $a_1, \cdots, a_n$ be a basis of a
maximal abelian subalgebra $\ca$ in $\cp$,  and $u_i:R^n\to
\ck_\ca^\perp$ smooth maps for $1\leq i\leq n$.  If
\refeq[bv]$$\o_\l=\sum\limits_{i=1}^n (a_i\l + u_i)dx_i$$ is a
flat connection 1-form on $R^n$ for all $\l\in C$, then there
exists a unique map $v:R^n\to \cp\cap \ca^\perp$ such that
$u_i=[a_i,v]$.

\proof
 Choose a basis $b_1, \cdots, b_n$ of $\ca$ such that each $b_i$ is
regular.  Write $b_j=\sum_ic_{ij}a_i$.  Make a change of
coordinates $x_i=\sum_j c_{ij}y_j$.  Then
$$\o_\l=\sum_i(a_i\l+u_i)dx_i= \sum_j (b_j\l + \tilde u_j)dy_j,$$
where $\tilde u_j= \sum_i c_{ij} u_i\in \ck\cap\ca^\perp$. Note that
$\o_\l$ is flat for all $\l$ if and only if
 \refeq[be]$$\cases{[b_i, \tilde
u_j]= [b_j, \tilde u_i], &\cr (\tilde u_j)_{x_i} - (\tilde
u_i)_{x_j} + [\tilde u_i, \tilde u_j]=0.&\cr}$$ Because $b_1,
\cdots, b_n$ are regular, $\ad(b_j)$ maps $\cp\cap\ca^\perp$
isomorphically to $\ck\cap\ca^\perp$.  Hence there exists a unique
$v_j\in \cp\cap\ca^\perp$ such that $\tilde u_j=\ad(b_j)(v_j)$ for
$1\leq j\leq n$.   Then the first equation of \refbe{} implies
that $$\ad(b_i)\ad(b_j)(v_j)=\ad(b_j)\ad(b_i)(v_i).$$ Since
$[b_i,b_j]=0$, $\ad(b_i)\ad(b_j)=\ad(b_j)\ad(b_i)$.   But $\ad(b_i)$
is injective on $\cp\cap\ca^\perp$ implies that
$v_i=v_j$, which will be denoted by $v$.  We compute directly to
get $$\eqalign{u_j&=\sum_ic^{ij} \tilde u_i = \sum_i c^{ij} [b_i,
v]\cr & =\sum_ic^{ij}\left[\sum_k c_{ki}a_k, v\right]= [a_j,
v],\cr}$$ where $(c^{ij})$ is the inverse of $(c_{ij})$. \qed

The following
Proposition is immediate:

\refclaim[bb] Proposition. If $\o$ is a flat $\cg$-valued connection 1-form and
$g:R^n\to G$ a smooth map, then the gauge transformation of $\o$ by $g$,
$$g\ast\o = g\o g^{-1} - dg g^{-1},$$ is also flat. Moreover, if $E$ is a trivialization of
$\o$, then
$Eg^{-1}$ is a trivialization of $g\ast \o$.

Consider the system of PDE for $$(A_1, \cdots, A_n, B_1, \cdots,
B_n):R^n\to \prod_{i=1}^n \cp \times \prod_{i=1}^n \ck$$ given by
the condition that $$\sum_{i=1}^n (A_i\l+ B_i) dx_i$$ is a flat
connection on $R^n$ for all $\l\in C$. Or equivalently,
$$(A_j\l+B_j)_{x_i} -(A_i\l + B_i)_{x_j} + [A_i\l + B_i, A_j\l +
B_j]=0$$ for all $\l\in C$.  By comparing coefficients of $\l^2,
\l$ and the constant term, we get
\refeq[br]$$\cases{[A_i,A_j]=0,&\cr (A_i)_{x_j}-(A_j)_{x_i} =
[A_i, B_j]+[B_i, A_j],&\cr (B_i)_{x_j}-(B_j)_{x_i} = [B_i,
B_j].&\cr}$$

\refclaim[cw] Proposition.  Let $\W_\l=\sum_i(A_i\l + B_i)dx_i$,
$A_i\in \cp$ and $B_i\in \ck$. If $[A_i,A_j]=0$ for all $1\leq i,
j\leq n$, then $\W_\l$ is flat for all $\l\in C$ if and only if
$\W_{\l_0}$ is flat for some non-zero real or pure imaginary $\l_0$.

\proof If $\W_{\l_0}$ is flat, then
\refeq[cq]$$(\l_0A_j+B_j)_{x_i}-(\l_0A_i+B_i)_{x_j}+[\l_0A_i+B_i,
\l_0A_j+B_j]=0.$$ So both  the $\ck$ and $\cp$  components of the left
hand side of the equation must be zero. Since $U/K$ is a symmetric
space, $$[\ck,\ck]\subset \ck, \quad [\ck,\cp]\subset \cp, \quad
[\cp, \cp]\subset \ck.$$ Equate the $\ck$ and $\cp$ components of
\refcq{} to get $$\cases{(B_j)_{x_i} -(B_i)_{x_j} + [B_i,
B_j]=0,&\cr (A_j)_{x_i} -(A_i)_{x_j} + [A_i, B_j] + [B_i,A_j]=0.
&\cr}$$ But this is the equation for $\W_\l$ to be flat for all
$\l\in C$. \qed

Restrict system \refbr{} to the case when all $B_i=0$ to get a
system for maps $(A_1,\cdots, A_n):R^n\to\cp\times \cdots \times
\cp$: \refeq[bg]$$\cases{[A_i, A_j]=0,&\cr
(A_i)_{x_j}=(A_j)_{x_i}, & for all $i\not=j$.\cr}$$ This is the
{\it Curved Flat system\/} associated to $U/K$ defined by Ferus and
Pedit in [FP]. Its Lax connection is $$\left[ {\p\over \p x_i}+A_i\l,
\,\, {\p\over \p x_j} + A_j\l \right]=0, \quad \forall\,\, \l\in
C.$$ The second equation of \refbg{} implies that $\sum_i A_i
dx_i$ is exact.  So we get

\refclaim[dg] Proposition. Let $A_i(x)\in \cp$. Then $\sum_i
\l A_i(x) dx_i$ is flat for all $\l\in C$ if and only if $[A_i, A_j]=0$ for all $1\leq i,
j\leq n$ and there exists a map
$X:R^n\to \cp$ such that $$dX=\sum_{i=1}^n A_i dx_i.$$

Let  $\o_\l$ be defined as in \refbv{}.  If $\o_\l$ is flat for
all $\l\in C$, then $\o_0=\sum_i u_i dx_i$ is flat.  Let $g$ be a
trivialization of $\o_0$.   A direct computation shows that the
gauge transformation of $\o_\l$ by $g$ is $$g\ast \o_\l=
\sum_{i=1}^n ga_ig^{-1}\l dx_i.$$ In other words, we have gauged
away the $\ck$-part of the Lax connection $\o_\l$ and the corresponding PDE is the
curved flat system.  If
$\ck=\ck_1\oplus
\ck_2$, then we can gauge away the $\ck_1$-part ($\ck_2$-part
respectively) of $\o_\l$.     To do this, we write $$u_i=
\xi_i+\eta_i\in \ck_1+ \ck_2.$$   Since $\o_0=\sum_i u_i dx_i$ is
flat, both $\sum_i \xi_i dx_i$ and $\sum_i \eta_idx_i$ are flat.
Let $g_1:R^n\to K_1$ and $g_2:R^n\to K_2$ be trivializations of
$\sum_i \xi_i dx_i$ and $\sum_i \eta_i dx_i$ respectively, i.e.,
$$g_1^{-1}dg_1= \sum_i \xi_idx_i, \quad g_2^{-1}dg_2= \sum_i
\eta_i dx_i.$$  Then the gauge transformation of $\o_\l$ by $g_1$ and $g_2$ are
is $$\eqalign{g_1\ast\o_\l &= \sum_i(g_1a_ig_1^{-1} \l +
\eta_i)dx_i,\cr g_2\ast \o_\l &= \sum_i (g_2a_ig_2^{-1} \l+
\xi_i)dx_i,\cr}$$
respectively. 

The $U/(K_1\times K_2)$-{\it system\/} I is the PDE for $g_1:R^n\to
K_1$ and $\eta_1, \cdots,\eta_n:R^n\to \ck_2\cap \ck_\ca^\perp$
such that $$\o_\l^I= \sum_i(g_1a_ig_1^{-1} \l + \eta_i)dx_i$$ is
flat for all $\l\in C$, i.e.,
\refeq[bl]$$\cases{[g_1^{-1}(g_1)_{x_i},
a_j]-[g_1^{-1}(g_1)_{x_j}, a_i] + [a_i,
g_1^{-1}\eta_jg_1]+[g_1^{-1}\eta_ig_1, a_j]=0,&\cr (\eta_j)_{x_i}
- (\eta_i)_{x_j} + [\eta_i, \eta_j]=0.&\cr}$$ Similarly, the
$U/(K_1\times K_2)$-{\it system\/} II is the PDE for $g_2:R^n\to
K_2$ and $\xi_1,\cdots, \xi_n:R^n\to \ck_1\cap\ck_\ca^\perp$ such
that $$\o_\l^{\II}=  \sum_i (g_2a_ig_2^{-1} \l+ \xi_i)dx_i$$ is
flat for all $\l\in C$, i.e.,
\refeq[bm]$$\cases{[g_2^{-1}(g_2)_{x_i},
a_j]-[g_2^{-1}(g_2)_{x_j}, a_i] + [a_i,
g_2^{-1}\xi_jg_2]+[g_2^{-1}\xi_i g_2, a_j]=0,&\cr (\xi_j)_{x_i} -
(\xi_i)_{x_j} + [\xi_i, \xi_j]=0.&\cr}$$

\bs

\newsection $G_{m,n}$-systems.\par

In this section, we assume $$U/K=G_{m,n}=O(m+n)/(O(m)\times
O(n)) \quad {\rm with\,\,} m\geq n.$$ We
write down the $G_{m,n}$-systems I and II explicitly.

Let $\cU =o(m+n)$, and $\sigma : \cU \to \cU $ be the involution
defined by $\sigma (X)=I_{m,n}^{-1}XI_{m,n}$, where
$$I_{m,n}=\pmatrix { I_m & 0 \cr 0 & -I_n }.$$ Let $\ck$  and
$\cp$ denote the $+1$ and $-1$ eigenspaces of $\s$ respectively.
Then  $\cU=\cK+ \cP$ is a Cartan decomposition of $U/K$, where $$
\eqalign{\cK&=o(m)\times o(n) =\left\{\pmatrix { Y_1 & 0 \cr
                                     0 & Y_2 } \bigg|\; Y_1\in o(m), \,
Y_2\in o(n)\right\},\cr
 \cP&=\left\{\pmatrix {0&\xi\cr -\xi^t&0\cr} \, \bigg| \; \xi\in
\cm_{m\times n}\right\}.\cr}$$  Here $\cm_{m\times n}$ is the set
of $m \times n$ matrices.  Note that $$\ca=\left\{\pmatrix {0 & -D
\cr
                        D^t & 0 \cr }\,\bigg| \; D=(d_{ij})\in \cm_{m\times
n}, \, d_{ij}=0 \,\, {\rm if\/}\,\, i\not=j\right\}.$$ is a
maximal abelian subalgebra in $\cP$ and $$\cP\cap \ca^\perp
=\left\{\pmatrix {0 &\xi \cr
                    -\xi^t & 0 } \, \bigg| \; \xi=(\xi_{ij}) \in  \cm_{m\times
n}, \,\, \xi_{ii}=0 \,\, {\rm for\/}\,\,1\leq i\leq n \right\}. $$
Let $$a_i=\pmatrix {0 &  -D_i \cr
                            D_i^t & 0\cr } ,$$
where $D_i\in \cm_{m\times n}$ is the matrix all whose entries are
zero except the ii-th entry is equal to $1$. Then $a_1, \dots ,
a_n$ form a basis of $\ca$. 
The $U/K$-system
\refat{} for this symmetric space is the following PDE for
$\xi=(\xi_{ij}): R^n\to \cm_{m\times n}$ with $\xi_{ii}=0$ for all
$1\leq i\leq n$: 
\refeq[ub]$$\cases{D_i\xi^t_{x_j}-\xi_{x_j}D_i^t
- D_j\xi^t_{x_i} + \xi_{x_i} D_j^t = [D_i\xi^t - \xi D_i^t,
D_j\xi^t- \xi D_j^t], & $i\not=j$,\cr D_i^t\xi_{x_j} -\xi^t_{x_j}
D_i -D_j^t\xi_{x_i} + \xi^t_{x_i} D_j = [D_i^t\xi -\xi^t D_i,
D_j^t\xi -\xi^t D_j], & $i\not=j$.\cr} $$ Its Lax connection \refau{} is
\refeq[av]$$\o_\lambda =\sum_i \left\{ \lambda
          \pmatrix {0  & -D_i \cr D_i^t & 0 }  +
                       \pmatrix {D_i\xi^t - \xi D_i^t &  0 \cr
                                         0  & -\xi^t D_i +D_i^t \xi}
\right\} dx_i. $$

Next we write down the $G_{m,n}$-system I and II
explicitly.   Let 
$$g=\pmatrix{A&0\cr 0&B\cr}\in O(m)\times O(n)$$
be a solution of $$g^{-1}dg =\o_0= \sum_{i=1}^n\pmatrix {D_i\xi^t
- \xi D_i^t &  0 \cr
                                         0  & -\xi^t D_i +D_i^t \xi}
dx_i.$$  Let $$g_1=\pmatrix{A&0\cr 0&I\cr}, \quad g_2=\pmatrix{I&
0\cr 0& B\cr}.$$ Write $$\xi=\pmatrix{F\cr G}, \quad
D_i=\pmatrix{C_i\cr 0\cr}, \quad A=(A_1, A_2),$$ where $F, C_i\in
gl(n)$,  $G\in \cm_{(m-n)\times n}$, $A_1\in \cm_{m\times n}$, and
$A_2\in \cm_{m\times (m-n)}$.  Then $$g_1a_i g_1^{-1} =
\pmatrix{0& -A_1 C_i\cr C_i A_1^t &0\cr}.$$

Let $$\eqalign{\cm_{m\times n}^0 &= \{A_1\in \cm_{m\times n}\n
A_1^tA_1=I\},\cr gl(n)_*&=\left\{(x_{ij})\in
gl(n)\, | \; x_{ii}=0,\ 1\le i\le n\right\}\cr}.$$ 

The $U/K$-{\it system I\/} is the PDE for $(A_1, F):
R^n\to \cm_{m\times n}^0\times gl_\ast(n)$ such that
\refeq[bo]$$\o_\l^I=\sum_{i=1}^n\pmatrix{0 & -A_1C_i\l \cr
                       C_iA_1^t\l  & -F^t C_i +C_i^t F} dx_i$$ is flat for
all $\l\in C$, i.e.,
\refeq[cn]$$\cases{(a_{ij})_{x_k}=f_{jk}a_{ik}, & if $k\not=j$,\cr
(f_{ij})_{x_j} + (f_{ji})_{x_i} + \sum_k f_{ik} f_{jk}=0, & if
$i\not=j$,\cr (f_{ij})_{x_k} = f_{ik}f_{kj}, & if $i, j, k$ are
distinct, \cr}$$ where $A=(a_{ij})$ and $F=(f_{ij})$.
 Note that equation \refcn{} is the condition that the above $\o^I_\l$ is
flat for $\l=1$.    So we have

\refclaim[cx] Proposition.  The following statements are
equivalent for map $(A_1,F):R^n\to \cm_{m\times n}^0\times
gl_\ast(n)$:
\item {(i)} $(A_1,F)$ is a solution of the $G_{m,n}$-system I \refcn{}.
\item {(ii)} $\o_\l^I$ defined by \refbo{}  is flat for all $\l\in C$.
\item {(iii)} $\o^I_\l$ defined by \refbo{} is flat for $\l=1$.

Note that if $a_{ij}\not=0$ for all $1\leq j\leq n$, then  the
first set of equations of \refcn{} implies that $F$ can be
computed from the $i$-th row of $A_1$ by $f_{jk}=(a_{ij})_{x_k}/
a_{ik}$.

\ms
Next we explain the reality conditions.  
Recall that a symmetric space $U/K$ is determined by a conjugate
linear Lie algebra involution $\tau$ and a complex linear Lie
algebra involution $\s$ on the complexified Lie algebra $\cu_\Cx =
\cu\otimes \Cx$ such that
\item {(i)} $\tau$ and $\s$ commute,
\item {(ii)} $\cu$ is the fixed point set of $\tau$, and  $\ck$ and $\cp$ are the
$+1, -1$ eigenspaces of $\s$ on $\cu$ respectively, and
$\cu=\ck+\cp$ is the Cartan decomposition.

\ni We still use $\tau$ and $\s$ to denote the corresponding
involutions on the group $U_\Cx$.

A map $g:C\to U_\Cx$ ($g:C\to \cu_\Cx$ respectively) is said to
satisfy the $U/K$-{\it reality condition\/} if
\refeq[gg]$$\tau(g(\bar\l))=g(\l), \quad \s(g(-\l))= g(\l).$$
A direct computation gives

\refclaim[jh] Proposition.   
\item {(i)} If $A(\l)=\sum_i A_i\l^i:C\to \cu_C$
satisfies the
$U/K$-reality condition, then $A_i\in \ck$ if $i$ is even and is in $\cp$
if $i$ is odd. 
\item {(ii)} The Lax pair $\o_\l$ defined by \refau{} for the $U/K$-system \refat{}
satisfies the $U/K$-reality condition. 

\refpar[gi] Definition.  A {\it frame\/} for a solution $v$ of the
$U/K$-system (I, II respectively) is a trivialization of the corresponding
Lax connection $\o_\l$ ($\o_\l^I, \o_\l^{\II}$ respectively) that satisfies the
$U/K$-reality condition.
\ms

\refpar[ji] Remark. 
\item {(i)} If $g:C\to U_C$ satisfies the $U/K$-reality condition, then $g(0)\in K$. 
\item {(ii)} If  $E(x,\lambda)$ is the trivialization  of
$\o_\l$ defined by \refau{} such that
$E(0,\l)$ satisfies the $U/K$-reality condition, then $E$ also
satisfies the $U/K$-reality condition. \ms

Let $p\in O(m)$ be a constant matrix. The gauge transform of
$\o_\l^I$ by $g=\pmatrix{p&0\cr 0&I\cr}$ is $$g\ast \o_\l^I=
\pmatrix{0&-pA_1C_i\l\cr C_i^t A_1^t p^t\l & -F^t C_i + C_i^t
F\cr},$$ which is flat for all $\l\in C$. Note that the
coefficient of $\l$
 in $g\ast \o_\l$ lies in $\cp$, the constant term lies in $\ck$, and
$$(pA_1)^t(pA_1)=A_1^tp^tpA_1= A_1^tA_1=I.$$  So it follows from
Proposition \refcx{} that

\refclaim[an] Corollary.  Let $(A_1, F)$ be a solution of the
$G_{m,n}$-system I
\refcn{}, and $p\in O(m)$ a constant matrix.  Then $(pA_1, F)$ is also
a solution of \refcn{}. Moreover, if $E^I$ is a frame for $(A_1,F)$,
then $E^Ip^{-1}$ is a frame for $(pA_1, F)$.

\refclaim[fm] Proposition.
\item {(i)} Suppose
$\xi=\pmatrix{F \cr G}$ is a solution of the $G_{m,n}$-system  \refub{}, $\o_\l$ the
corresponding Lax connection, and $g=\pmatrix{A&0\cr 0&B\cr}:R^n\to
O(m)\times O(n)$ satisfies $g^{-1}dg=\o_0$.   Write $A=(A_1,A_2)$
with $A_1\in \cm_{m\times n}$. Then $(A_1, F)$ is a solution of the
$G_{m,n}$-system I  \refcn{}.
\item {(ii)} Conversely, if $(A_1, F)$ is a solution of the $G_{m,n}$-system I \refcn{},
then there exists an $\cm_{(m-n) \times n}$-valued map
$G$ such that $\xi=\pmatrix{F \cr G}$ is a solution of \refub{}.

\proof

(i) follows from the definition of the $U/K$-system I. 

(ii) Choose $A_2$ so that $A=(A_1,A_2)\in O(m)$.
 Let $g=\pmatrix {A^t & 0\cr 0 & I} $. Then the gauge
transformation of $g$ on $\theta _\lambda^I$ is $$\eqalign { g*
&\theta_\lambda^I = g\theta_\lambda^I g^{-1}-dgg^{-1}\cr
              =&\sum_i\left\{ \lambda \pmatrix { 0 & 0 & -C_i \cr 0 & 0 & 0 \cr
                  C_i & 0 & 0 }+
           \pmatrix {  A_1^t\,(A_1)_{x_i} & A_1^t\,(A_2)_{x_i} & 0 \cr
                        A_2^t\,(A_1)_{x_i} & A_2^t\,(A_2)_{x_i} & 0 \cr
                                0 & 0 & -F^tC_i+C_iF }\right\}dx_i. }$$
Although this does not have the same shape as the Lax pair $\o_\l$ of the
$G_{m,n}$-system, we show below that it can be gauged to one.  
From  \refcn{}, we have
\refeq[oa]$$dA_1=A_1\sum(C_iF^t-FC_i)dx_i+\zeta\sum C_idx_i$$
for some
$\zeta : R^n \rightarrow \cm_{m\times n}$. Thus
\refeq[ob]$$A_2^tdA_1=A_2^t\zeta\sum C_idx_i.$$ Since
$A^{-1}dA$ is flat and $$A^{-1}dA= \pmatrix{A_1^tdA_1& A_1^t
dA_2\cr A_2^tdA_1& A_2^tdA_2\cr},$$ we have
\refeq[oc]$$dA_2^t\wedge dA_2 +A_2^t dA_1
\wedge A_1^t dA_2 + A_2^tdA_2\wedge A_2^tdA_2=0.$$ 
By \refob{},
$A_1^tdA_2=(dA_2^tA_1)^t=-(A_2^tdA_1)^t=
-\sum C_i\zeta^tA_2dx_i$. So
it follows from \refoc{} that $A_2^tdA_2$ is flat, and hence there
exists $h : R^n\rightarrow O(m-n)$ such that $h^{-1}dh=A_2^tdA_2$.
Thus if we do a gauge transform by $\hat h=\pmatrix {I & 0 &
0 \cr 0 & h & 0\cr 0 & 0 & I}$ on $g*\theta_\lambda^I$, the
resulting connection 1-form is $$\eqalign{&\hat h\ast(g\ast
\o_\l^I)=\cr &\sum\left\{\lambda\pmatrix{0 & 0 & -C_i\cr 0 & 0 &
0\cr C_i & 0 & 0}+ \pmatrix{A_1^t(A_1)_{x_i} & -C_i \zeta^tA_2h^t
& 0\cr
               hA_2^t\zeta C_i & 0 & 0 \cr
                0 & 0 & -F^tC_i+C_iF}\right\}dx_i.\cr}$$
Set $$G=-hA_2^t\zeta.$$
 From \refoa{}, we have
$A_1^t(A_1)_{x_i}-(C_iF^t -FC_i)=YC_i$, where $Y=A_1^t\zeta$.
Since the left-hand side is skew-symmetric, so is $YC_i$. But
$YC_i=-C_iY$ for all $1\leq i\leq n$ implies that $Y=0$. It follows that
$\hat h
\ast(g\ast
\o_\l^I)$ is the $\o_\l$ defined by \refav{} with
$\xi=\pmatrix{F\cr G}$.  In other words, $\xi$ is a solution of
\refub{}.
\qed

\ms The $G_{m,n}$-system II is the PDE for
$$(F,G,B):R^n\to gl_\ast(n)\times \cm_{(m-n)\times n} \times O(n)$$
such that \refeq[bp]$$\o_\l^{\II}=  \sum_{i=1}^n
\pmatrix{-FC_i+C_iF^t & C_i G^t & -C_i B^t\l\cr -GC_i& 0 & 0\cr
BC_i\l & 0 &0\cr} dx_i$$ is flat for all $\l\in C$, i.e.,
\refeq[cu]$$\cases{(f_{ij})_{x_i} + (f_{ji})_{x_j} + \sum_{k=1}^n
f_{ki}f_{kj} + \sum_{k=1}^{m-n} g_{ki}g_{kj} =0,& if $i\not=j$,\cr
(f_{ij})_{x_k} = f_{ik}f_{kj}, & if $i, j, k$ are distinct,\cr
(b_{ij})_{x_k} = f_{kj} b_{ik},& if $j\not=k$,\cr (g_{ij})_{x_k}=
f_{kj} g_{ik}, & if $j\not=k$.\cr}$$ As a consequence of
Proposition \refcw{} we get

\refclaim[df] Proposition.  
Given $(F,G,B):R^n\to gl_\ast(n)\times \cm_{(m-n)\times n} \times O(n)$, 
the following statements  are equivalent:
\item {(i)} $(F,G,B)$ is a solution of the $G_{m,n}$-system II \refcu{}.
\item {(ii)} $\o_\l^{\II}$ defined by \refbp{} is flat for all $\l\in C$,
\item {(iii)} $\o_\l^{\II}$ defined by \refbp{} is flat for $\l=1$.

Note that if $E$ is a frame for $v=\pmatrix{F\cr G\cr}$, then
$E(x,0)=\pmatrix{A(x)&0\cr 0&B(x)\cr}$ and $E^{\II}=E\pmatrix{A^{-1}&0\cr
0&I_n\cr}$ is a frame for $(F,G,B)$.

If $b_{ji}\not=0$ for all $1\leq i\leq n$, then  the third
equation of \refcu{} implies that $f_{ki}= (b_{ji})_{x_k}/b_{jk}$
if $k\not=i$.  In other words, generically system \refcu{} depends
only on $B$ and $G$.

\refclaim[ig] Corollary.  If $(F,G,B)$ is a solution of the $G_{m,n}$-system II \refcu{}
and $C\in O(n)$ a constant matrix, then $(F,G,CB)$ is also a
solution of \refcu{}.  Moreover, if $E^{\II}$ is a frame for
$(F,G,B)$, then $E^{\II}C^{-1}$ is a frame for $(F,G,CB)$.

\bs Let  $U/K=G_{m,n+1}$ with $m\geq n+1$.  The
rank of $U/K$ is $n+1$.  So the corresponding $U/K$-system has
$(n+1)$ independent variables.  Below we  consider a partial
$U/K$-system of n-variables. Let $(m_1, m_2, m_3, m_4)=
(n,m-n,n,1)$. We partition a matrix $A$ in $o(m+n+1)$ into
$4\times 4$ blocks $A=(A_{ij})$, where $A_{ij}\in \cm_{m_i\times
m_j}$. Let $$a_i=\pmatrix {0&0&-C_i&0\cr 0&0&0&0\cr C_i&0&0&0\cr
0&0&0&0 } , \, {\rm where}\ C_i=\diag(0, \dots , 1, \dots,0)  \
{\rm as \ before}. $$ Then the space $\ca$ spanned by $a_1,
\cdots, a_n$ is an n-dimensional abelian subspace of $\cp$.   Let
$\ck_\ca$ and $\cu_\ca$  denote the centralizer of $\ca$ in $\ck$
and $\cu$ respectively.  Then $\cp\cap\cu_\ca^\perp$ is the space
of elements of the form $$v=\pmatrix {0&0&F&b\cr 0&0&G&0\cr
-F^t&-G^t&0&0\cr -b^t&0&0&0 }, $$ where, $F\in gl(n)_*, \
b\in\cm_{n\times 1}, \ G\in \cm_{(m-n)\times n} $.

The {\it partial $G_{m,n+1}$-system\/} is the
PDE for $(F, G,b):R^n\to gl(n)_\ast\times \cm_{(m-n)\times
n}\times \cm_{n\times 1}$ such that
\refeq[aw]$$\eqalign{\Theta_\lambda= &\sum_i \lambda \pmatrix {0
&0 &-C_i &0\cr 0&0&0&0\cr
            C_i&0&0&0\cr 0&0&0&0}dx_i\cr
 &+\sum_i\pmatrix {-FC_i+C_iF^t&C_iG^t
&0&0\cr -GC_i&0&0&0\cr 0&0&-F^tC_i+C_iF& C_ib\cr 0&0&-b^tC_i&0}
 dx_i}$$
is flat for all $\l\in C$, i.e.,
\refeq[fr]$$\cases{(f_{ij})_{x_i}+(f_{ji})_{x_j} +\sum_k
f_{ki}f_{kj} +\sum_k g_{ki}g_{kj} =0, & if $i\not=j$,\cr
(f_{ij})_{x_k}= f_{ik}f_{kj}, & if $i, j, k$ distinct,\cr
(g_{ij})_{x_k}= g_{ik}f_{kj}, & if $k\not=j$,\cr
 (f_{ij})_{x_j} +
(f_{ji})_{x_i} + \sum_k f_{ik}f_{jk}  + b_i b_j=0, & if $i\not=j$,
\cr (b_i)_{x_j}= f_{ij} b_j.&\cr}$$

The {\it partial $G_{m,n+1}$-system I\/} is the
PDE for maps 
$$(A_1, F, b):R^n\to \cm_{m\times n}^0\times
gl_\ast(n)\times \cm_{n\times 1}$$ such that
\refeq[cs]$$\Theta^I_\l=\sum_i\pmatrix{0& -A_1C_i\l & 0\cr
C_iA_1^t\l & -F^tC_i+C_iF & C_i b\cr 0& -b^tC_i &0\cr} dx_i$$ is
flat for all $\l\in C$, i.e., \refeq[cl]$$\cases{(f_{ij})_{x_j} +
(f_{ji})_{x_i} + \sum_{k=1}^n f_{ik}f_{jk} +
 b_ib_j=0, \quad i\not=j, &\cr
(f_{ij})_{x_k}= f_{ik}f_{kj},  \quad i,j,k  \,\,{\rm distinct},
&\cr (b_i)_{x_j} = f_{ij} b_j, \quad i\not=j, &\cr (a_{ki})_{x_j}=
f_{ij} a_{kj}, \quad i\not=j.& \cr}$$ By Proposition \refcw{} we
have

\refclaim[cr] Proposition. Given a map $(A_1, F, b):R^n\to \cm_{m\times n}^0\times
gl_\ast(n)\times \cm_{n\times 1}$,  the following statements are equivalent:
\item {(i)} $(A_1,F,b)$ is a solution of the partial $G_{m,n+1}$-system I \refcl{},
\item {(ii)} $\Theta_\l^I$ defined by \refcs{} is flat for all $\l\in C$,
\item {(iii)} $\Theta_\l^I$ defined by \refcs{} is flat
for $\l=1$.

\bs

\newsection $G_{m,n}^1$-Systems.\par

In this section, we assume $m\geq n+1$ and
$$U/K=G_{m,n}^1 =O(m+n,1)/(O(m)\times O(n,1)).$$  We write
down the  $G_{m,n}^1$-system I and II explicitly. \ms

Let $\cU =o(m+n,1)=\left\{ X\in gl(m+n+1) \ |\
X^tI_{m+n,1}+I_{m+n,1}X=0 \right\}$ and $\sigma : \cU \to \cU $ be
an involution defined by $\sigma (X)=I_{m,n+1}^{-1}XI_{m,n+1}$,
where $I_{p, q}=\pmatrix { I_p & 0 \cr 0 & -I_q }$. Let $\ck$  and
$\cp$ denote the $+1$ and $-1$ eigenspaces of $\s$ respectively.
Then the Cartan decomposition is \ $\cU=\cK+ \cP$, where $$
\eqalign{\cK&=o(m)\times o(n,1) =\left\{\pmatrix { Y_1 & 0 \cr
                                     0 & Y_2 } \bigg|\; Y_1\in o(m),
Y_2\in o(n,1)\right\},\cr \cP&=\left\{\pmatrix {0&\xi\cr
-J\xi^t&0\cr} \, \bigg| \; \xi\in \cm_{m\times
(n+1)}\right\}.\cr}$$ Here $\cm_{m\times n}$ is the set of $m
\times n$ matrices and $J=I_{n,1}=\diag(1, \cdots, 1, -1)$. It is easy
to see that
$$\ca=\left\{\pmatrix {0 & -DJ \cr
                        D^t & 0 \cr }\,\bigg| \; D=(d_{ij})\in \cm_{m\times
(n+1)}, \, d_{ij}=0 \,\, {\rm if\/}\,\, i\not=j\right\}$$ is a
maximal abelian subalgebra in $\cP$. Let $$a_i=\pmatrix {0 & -D_iJ
\cr
                            D_i^t & 0\cr } ,$$
where $D_i\in \cm_{m\times (n+1)}$ is the matrix all whose entries
are zero except that the ($i,i$)-th entry is equal to $1$. Then $a_1,
\dots , a_{n+1}$ form a basis of $\ca$. 
The $U/K$-system \refat{} for
this symmetric space is the PDE for $\xi=(\xi_{ij}): R^{n+1}\to
\cm_{m\times (n+1)}$ with $\xi_{ii}=0$ for all $1\leq i\leq n+1$ such
that
\refeq[aab]$$\o_\lambda =\sum_i \left\{ \lambda
          \pmatrix {0  & -D_iJ \cr D_i^t & 0 }  +
                       \pmatrix {D_i\xi^t - \xi D_i^t &  0 \cr
                                         0  & -J\xi^t D_iJ +D_i^t \xi}
\right\} dx_i $$ is a family of flat connections on $R^{n+1}$ for
all $\l\in C$, i.e.,
\refeq[aaa]$$\cases{D_i\xi^t_{x_j}-\xi_{x_j}D_i^t - D_j\xi^t_{x_i}
+ \xi_{x_i} D_j^t = [D_i\xi^t - \xi D_i^t, D_j\xi^t- \xi D_j^t],
\quad i\not =j,\cr D_i^t\xi_{x_j} -J\xi^t_{x_j} D_iJ
-D_j^t\xi_{x_i} + J\xi^t_{x_i} D_jJ \cr
  \qquad \qquad  \qquad \qquad = [D_i^t\xi-J\xi^t D_iJ, D_j^t\xi -J\xi^t D_jJ].
\qquad i \not =j.\cr}$$

Let
$$gl(n+1)_*=\left\{(x_{ij})\in gl(n+1)\, | \; x_{ii}=0,\ 1\le i\le
n+1\right\}.$$ 
The $G_{m,n}^1$-system I is the PDE for
$(A_1,F):R^{n+1}\to \cm^0_{m\times (n+1)} \times gl_\ast(n+1)$
such that \refeq[bx]$$\o_\l^I=\sum_{i=1}^n \pmatrix{0& -A_1C_i
J\l\cr C_iA_1^t\l & -JF^tC_iJ + C_i F\cr} dx_i$$ is a flat
connection on $R^{n+1}$ for all $\l\in C$, i.e.,
\refeq[cc]$$\eqalign{&\cases{-(A_1)_{x_j}C_i J + (A_1)_{x_i} C_j J
= -A_1C_iJ \eta_j + A_1C_j J \eta_i,&\cr (\eta_i)_{x_j}
-(\eta_j)_{x_i} =[\eta_i, \eta_j], &\cr} \cr &\quad {\rm where\,}
\quad \eta_i= -JF^tC_iJ + C_i F.\cr}$$

Next we write down the $U/K$-system II. Write $$\xi=\pmatrix{F\cr
G}, \quad D_i=\pmatrix{C_i\cr 0\cr}.$$ Then the
$G_{m,n}^1$-system II is the PDE for
$(F,G,B):R^{n+1}\to gl_\ast(n+1)\times \cm_{(m-n-1)\times (n+1)}
\times O(n,1)$ such that \refeq[aaf]$$\o_\l^{\II}= \sum_{i=1}^{n+1}
\pmatrix{-FC_i+C_iF^t & C_i G^t & -C_i B^tJ\l\cr -GC_i& 0 & 0\cr
BC_i\l & 0 &0\cr} dx_i$$ is flat for all $\l\in C$, i.e.,
\refeq[aad]$$\cases{(f_{ij})_{x_i} + (f_{ji})_{x_j} +
\sum_{k=1}^{n+1} f_{ki}f_{kj} + \sum_{k=1}^{m-n-1} g_{ki}g_{kj}
=0,\quad {\rm if} \ \ i\not=j,\cr (f_{ij})_{x_k} = f_{ik}f_{kj},
\qquad \qquad {\rm if} \ \ i, j, k  \ {\rm are \ distinct},\cr
(b_{ij})_{x_k} = f_{kj} b_{ik}, \qquad \qquad {\rm if}\ \
j\not=k,\cr (g_{ij})_{x_k}= f_{kj} g_{ik}, \qquad \qquad {\rm if}
\ \ j\not=k.\cr}$$ It follows from Proposition \refcw{} that

\refclaim[aae] Proposition.  The following statements are
equivalent for maps $(F,G,B):R^{n+1}\to$ $gl_\ast(n+1)\times
\cm_{(m-n-1)\times (n+1)} \times O(n,1)$ :
\item {(i)} It is a solution of the $G_{m,n}^1$-system II \refaad{}.
\item {(ii)} $\o_\l^{\II}$ defined by \refaaf{} is a flat connection on
$R^{n+1}$ for all $\l\in C$,
\item {(iii)} $\o_\l^{\II}$ defined by \refaaf{} is flat for $\l=1$.

Next we write down the partial $G_{m,n}^1$-systems of $n$ variables. Let
$(m_1, m_2,$
$m_3,m_4)=(n,m-n,n,1)$. We partition a matrix $A$ in $o(m+n,1)$
into $4\times 4$ blocks $A=(A_{ij})$, where $A_{ij}\in
\cm_{m_i\times m_j}$. Let $$a_i=\pmatrix {0&0&-C_i&0\cr 0&0&0&0\cr
C_i&0&0&0\cr 0&0&0&0 }, \, {\rm where}\ C_i=\diag(0, \dots , 1,
\dots,0)  \ {\rm as \ before}. $$ Then the space $\ca$ spanned by
$a_1, \cdots, a_n$ is an $n$-dimensional abelian subspace of $\cp$
(it is not a maximal one).   Let $\ck_\ca$ and $\cu_\ca$  denote
the centralizer of $\ca$ in $\ck$ and $\cu$ respectively.  Then
$\cp\cap\cu_\ca^\perp$ is the space of elements of the form
$$v=\pmatrix {0&0&F&b\cr 0&0&G&0\cr -F^t&-G^t&0&0\cr
b^t&0&0&0\cr},$$ where, $F\in gl(n)_*$, $ b\in\cm_{n\times 1}$,
and $ G\in \cm_{(m-n)\times n} $.

The partial $G_{m,n}^1$-system is the system
for $(F, G, b)$ such that
 \refeq[ka]$$\eqalign{\tau_\lambda= &\sum_i \lambda \pmatrix {0 &0 &-C_i
&0\cr 0&0&0&0\cr
            C_i&0&0&0\cr 0&0&0&0}dx_i\cr
 &+\sum_i\pmatrix {-FC_i+C_iF^t&C_iG^t
&0&0\cr -GC_i&0&0&0\cr 0&0&-F^tC_i+C_iF&C_ib\cr 0&0&b^tC_i&0}
 dx_i,}$$

is flat for all $\l \in C$, i.e.,

\refeq[kb]$$\cases{(f_{ij})_{x_i}+(f_{ji})_{x_j} +\sum_k
f_{ki}f_{kj} +\sum_k g_{ki}g_{kj} =0, & if $i\not=j$,\cr
(f_{ij})_{x_k}= f_{ik}f_{kj}, & if $i, j, k$ distinct,\cr
(g_{ij})_{x_k} = g_{ik} f_{kj}, & if $j\not= k$,\cr (f_{ij})_{x_j}
+ (f_{ji})_{x_i} + \sum_k f_{ik}f_{jk} - b_i b_j=0, & if
$i\not=j$, \cr (b_i)_{x_j}= f_{ij} b_j,&if $i\not=j$. \cr}$$

 The partial $G_{m,n}^1$-system I is the PDE for
$(A_1, F, b):R^n\to gl(n)_\ast\times \cm_{(m-n)\times n}\times
\cm_{n\times 1}$ such that
\refeq[aak]$$\tau^I_\l=\sum_i\pmatrix{0& -A_1C_i\l & 0\cr
C_iA_1^t\l & -F^tC_i+C_iF & C_i b\cr 0& b^tC_i &0\cr} dx_i$$ is
flat for all $\l\in C$.  It follows from Proposition \refcw{} that

\refclaim[aai] Proposition. The following statements are
equivalent for maps $(A_1, F, b):R^n\to \cm_{m\times n}^0\times
gl_\ast(n)\times \cm_{n\times 1}$:
\item {(i)} It is a solution of the partial $G_{m,n}^1$-system I:
\refeq[aaj]$$\cases{(f_{ij})_{x_j} + (f_{ji})_{x_i} + \sum_{k=1}^n
f_{ik}f_{jk} - b_ib_j=0, \quad i\not=j, &\cr (f_{ij})_{x_k}=
f_{ik}f_{kj},  \quad i,j,k  \,\,{\rm distinct}, &\cr (b_i)_{x_j} =
f_{ij} b_j, \quad i\not=j, &\cr (a_{ki})_{x_j}= f_{ij} a_{kj},
\quad i\not=j.& \cr}$$
\item {(ii)} $\tau_\l^I$ defined by \refaak{} is a flat connection on $R^n$
for all $\l\in C$.
\item {(ii)} $\tau_\l^I$ defined by \refaak{} is a flat connection on $R^n$
for $\l=1$.

\bs


\newsection Moving frame method for submanifolds.\par

In this section, we  review some elementary theory of submanifolds
in Euclidean space (cf. Chapter 2 of [PT]). The basic local
invariants of submanifolds in $R^N$ are the first, second
fundamental forms and the induced normal connection.  These
invariants satisfy the Gauss, Codazzi and Ricci equations.  If we
use the method of moving frames, then these equations arise as the
condition for a $o(N)$-valued 1-form to be flat.  Natural
geometric conditions on submanifolds in $R^N$ often lead to
interesting PDE.  If in addition
these submanifolds come in  a family that depends holomorphically
on a parameter, then the corresponding PDE has a Lax connection.
Therefore we can use techniques from soliton theory to study these
submanifolds.

Let $X:M^n\to R^{n+m}$ be an immersed submanifold.  Henceforth we
agree on the following index conventions unless otherwise stated:
$$1\le i, j, k\le n ,\quad n+1\le\alpha,\beta, \gamma \le
n+m,\quad 1\le A, B, C \le n+m.$$   Let $e_A$ be a local
orthonormal frame on $M$ such that $e_\a$ are normal field, and
$$E=(e_1, \cdots, e_{n+m}),$$ i.e., the
$A$-th column of $E$ is $e_A$.  Thus $E\in O(m+n)$. Let $\w_i$ be the
local orthonormal frame of $T^*M$ dual to $e_i$. Then $$dX=\sum_i
\w_ie_i,$$ and the first fundamental form is $$I=\sum_i \w_i^2.$$
Set $\w_{AB}= \li e_A, de_B\ri$, where $\li , \ri$ is the inner
product on $R^{n+m}$, or equivalently, $$de_B= \sum_{A}
\w_{AB}e_A = -\sum_A \w_{BA} e_A.$$ Write this in matrix form to
get  $$dE= E(\w_{AB}), $$ i.e.,
 $$(\w_{AB})= E^{-1} dE.$$
In other words, $(\w_{AB})$ is a flat $o(n+m)$-valued 1-form on $M$,
and $\w_{AB}$   satisfies the Maurer-Cartan equation:
\refeq[bz]$$d\w_{AB}=-\sum_C \w_{AC}\wedge\w_{CB}.$$ The structure
equation is \refeq[af]$$d\omega_i=
-\sum_{j}\omega_{ij}\wedge\omega_j, \quad \w_{ij}+\w_{ji}=0.$$
Here $(\w_{ij})$ is the Levi-Civita connection 1-form for the
induced metric $I$, and it can be computed in terms of $\w_1,
\cdots, \w_n$ by solving the structure equation \refaf{}.

The second fundamental form of $M$ is $$ \II = \sum_{i,\a}
\w_{i\a}\w_i e_\a.$$ Given $\xi\in\nu(M)_x$, the {\it shape
operator\/} $A_\xi$ is the self-adjoint operator defined by $<\II,
\xi>$, i.e., $$<\II(v_1, v_2), \xi> = <A_\xi(v_1), v_2>.$$ The
eigenvalues and eigendirections of $A_\xi$ are the {\it principal
curvatures\/} and {\it principal directions\/} of $M$ with respect to
$\xi$ respectively.

The induced connection $\K^\perp$ on the normal bundle $\nu(M)$ is
defined by $$\K^\perp \xi= (d\xi)^\perp$$ for any normal field
$\xi$, where $(d\xi)^\perp$ is the normal component of $d\xi$.  In
particular, $$\K^\perp e_\a= -\sum_\b \w_{\a\b} e_\b.$$ The normal
curvature is $$\W^\perp_{\a\b} = d\w_{\a\b}+ \sum_\gamma
\w_{\a\g}\wedge\w_{\g\b}.$$

Equation \refbz{} gives the fundamental equations for $M$: The
$ij$-th entry gives the Gauss equation
\refeq[ag]$$\eqalign{d\omega_{ij} &=-\sum_k\omega_{ik}\wedge
\omega_{kj} - \sum_\alpha\omega_{i\alpha}\wedge\omega_{\alpha
j},\cr &=-\sum_k \w_{ik}\wedge \w_{kj} + \W_{ij},\cr}$$
where $\W_{ij}= \sum_\a \w_{i\a} \wedge
\w_{j\a}$ is the Riemann curvature tensor of I.
 The $i\a$-th
entry gives  the Codazzi equations
\refeq[ah]$$d\omega_{i\alpha}=-\sum_j\omega_{ij}\wedge
\omega_{j\alpha} -\sum_\b \w_{i\b}\wedge\w_{\b\a},$$ and the
$\a\b$-th entry gives the Ricci's equations:
$$d\w_{\a\b}+\sum_\gamma \w_{\a\g}\wedge\w_{\g\b}
=\W^\perp_{\a\b}= -\sum_i \w_{\a i}\wedge \w_{i \b}.$$

\refclaim[bj] Fundamental Theorem of Submanifolds in Euclidean
space.  Let $(M,g)$ be an n-dimensional Riemannian manifold, $\xi$
a rank $m$ orthogonal vector bundle over $M$, $\K_0$ an
$O(m)$-connection on $\xi$, and $b$ a smooth section of
$S^2(T^*M)\otimes \xi $.  If $g, b$ and $\K_0$ satisfy the
Gauss-Codazzi-Ricci equations, then there exist a local isometric
immersion of $M$ into $R^{n+m}$ and a bundle isomorphism between
$\xi$ and $\nu(M)$ such that $g, b$ and $\K_0$ are the first,
second fundamental forms and induced normal connection
respectively.  This immersion is unique up to rigid motions.

A normal vector field $\eta$ is {\it parallel\/} if $\K^\perp
\eta=0$.  The normal bundle $\nu(M)$ is {\it flat\/} if the
induced normal connection $\K^\perp$ is flat, i.e.,
\refeq[bk]$$\W^\perp_{\a\b}= d\w_{\a\b} +\sum_\g
\w_{\a\g}\wedge\w_{\g\b} =0.$$ If $\nu(M)$ is flat, then it
follows from \refbk{} that there exists a local parallel normal frame.
So we may assume that $(e_{n+1}, \cdots, e_{n+m})$ is parallel,
i.e., $$\w_{\a\b}=0.$$   Then the Ricci equation is
\refeq[ai]$$0=d\w_{\a\b}=\sum_i \w_{i\a}\wedge\w_{i \b}.$$ In
particular, this implies that $[A_{e_\a}, A_{e_\b}]=0$ for all
$\a, \b$. Hence the family $\{A_v\n v\in \nu(M)_p\}$ of shape
operators of $M$ at $p$ is a family of commuting self-adjoint
operators on $TM_p$, and generically there is a smooth common
eigenframe.

\ms

Henceforth, we assume that $\nu(M)$ is flat,  and $e_\a$ is
parallel, i.e., $\w_{\a\b}=0$.
So equation
\refah{} becomes 
\refeq[jy]$$d\w_{i\a}= -\sum_j \w_{ij}\wedge \w_{j\a},
\quad \w_{ij}+\w_{ji}=0.$$ 

Next we recall the well-known Cartan
Lemma.

\refclaim[abb] Cartan Lemma.  If $\tau_1, \cdots, \tau_n$ are
linearly independent 1-forms on an n-dimensional manifold, then
there exists a unique $o(n)$-valued 1-form $(\tau_{ij})$ such that
\refeq[ax]$$d\tau_i = -\sum_{j=1}^n \tau_{ij}\wedge\tau_j, \quad
\tau_{ij}+\tau_{ji}=0, \quad 1\leq i\leq n.$$

A direct computation gives

\refclaim[ab] Corollary.  If $\tau_i= b_i dx_i$, then the solution
$(\tau_{ij})$ for system \refax{} is given by $$\tau_{ij}=
{(b_i)_{x_j}\over b_j} dx_i - {(b_j)_{x_i}\over b_i} dx_j.$$

\refpar[ay] Definition.  Let $M^n$ be a submanifold in a
Riemannian manifold $N^{n+m}$.  The normal bundle $\nu(M)$ is
called {\it non-degenerate\/} if the space  of shape
operators $\{A_v\n v\in \nu(M)_p\}$ has dimension $n$ for all $p\in
M$.

\ms

When $M$ is a surface in $R^3$, it is well-known that if $p\in M$
is not umbilic then there exist a local coordinates $x_1, x_2$
such that ${\p\over \p x_1}, {\p\over \p x_2}$ are eigenvectors
for the shape operator $A_{e_3}$.  In other words, the two
fundamental forms are of the form $$I=a_1^2 dx_1^2 + a_2^2 dx_2^2,
\quad \II= b_1 dx_1^2 +b_2 dx_2^2$$ for some functions $a_i, b_i$.
Such coordinates  are called {\it line of curvature
coordinates\/}.  We generalize this notion to submanifolds below.

\refpar[al] Definition.  Let $M$ be an n-dimensional submanifold
of $R^m$ with flat normal bundle, and
$\{e_\a\}=\{e_{n+1}, \cdots, e_m\}$ a local orthonormal parallel
normal frame.    Local coordinates $x_1,\cdots, x_n$ are called
{\it line of curvature coordinates\/}  with respect to $\{e_\a\}$
if $\{{\p\over \p x_i}\}$ is a common orthogonal eigenbasis for
the shape operators $A_v$ for all $v\in \nu(M)$, or equivalently, if
the two fundamental forms of $M$ are of the form $$\eqalign{I&=
\sum_{i=1}^n a_i^2dx_i^2, \cr \II&=\sum_{i=1, j=1}^{n, m-n}
b_{ji}dx_i^2 e_{n+j}\cr}$$ for some smooth functions $a_i$ and
$b_{ij}$.  If in addition $a_1^2+ \cdots + a_n^2=1$ ($a_1^2+\cdots
+ a_{n-1}^2 -a_n^2 = \pm 1$ resp.), then $x_1, \cdots, x_n$ are
called {\it spherical (hyperbolic resp.) line of curvature
coordinates\/}.

\ms

Let $M$ be a submanifold of $R^{n+m}$ with flat normal bundle.
Although generically there exist orthonormal tangent frame $e_i$
that are common eigenbasis for all shape operators, we can not
always find coordinates $x_1, \cdots, x_n$ so that ${\p\over \p
x_i}$ is parallel to $e_i$ for all $i$. However,  n-dimensional
submanifolds in $R^{2n-1}$ with constant sectional curvature are
known to have flat and non-degenerate normal bundle, and admit
spherical line of curvature coordinates (cf. [Ca], [M]).

\bs


\newsection Submanifolds associated to $G_{m,n}$-systems.\par

In this section, we describe submanifolds associated to the
$G_{m,n}$-systems I and II. 

 Recall that the Lax connection of
the $G_{m,n}$- system I \refcn{} is $\o_\l^I$ defined by \refbo{}.  Let $\o_{ij}$
denote the $ij$-th entry of $\o_1^I$ (i.e., $\o_\l^I$ at $\l=1$).
Then $$\o_{ij}=\cases{0, \qquad {\rm if\,\ } 1\leq i,j\leq m, & \cr
a_{j,i-m}dx_{i-m}, \qquad {\rm if\, \ }  m+1\leq i\leq m+n, 1\leq
j\leq m, &\cr   -f_{j-m, i-m}dx_{j-m} + f_{i-m, j-m} dx_{i-m},
\qquad {\rm if\ \ } m+1\leq i, j\leq m+n. & \cr}$$ Set
$$\eqalign{&\w_i=\o_{m+i,1}= a_{1i} dx_i, \quad \w_{ij}=
\o_{m+i,m+j}= -f_{ji} dx_j + f_{ij} dx_i, \,\, {\rm if\,\,} 1\leq
i, j\leq n,\cr & \w_{n+i,n+j}= \o_{ij} =0, \quad {\rm if\/}\quad
1\leq i, j\leq  m,\cr & \w_{i,n+j}=\o_{m+i, j}= a_{ji} dx_i, \quad
{\rm if\/}\quad 1\leq i\leq n, 1\leq j\leq m.\cr}$$ Since
$\o_1^I=(\o_{ij})$ is flat, we have $$\eqalign{d\w_i&= d\o_{m+i,
1} \cr &= -\sum_{j=1}^m \o_{m+i,j}\wedge \o_{j 1} -\sum_{j=1}^n
\o_{m+i, m+j} \wedge\o_{m+j,1}\cr &= -\sum_{j=1}^n \o_{m+i, m+j}
\wedge\o_{m+j,1} = - \sum_{j=1}^n \w_{ij}\wedge \w_j.\cr}$$ If
$\sum_{i=1}^n \w_i^2$ is non-degenerate, then $(\w_{ij})$ is the
Levi-Civita connection 1-form for the metric $\sum_{i=1}^n
\w_i^2$.  Let $E$ be a trivialization of $\o_1^I$, i.e.,
$E^{-1}dE= \o_1^I$. Let $e_{n+i}$ denote the $i$-th column of $E$
for $1\leq i\leq m$, and let $e_i$ denote the $(m+i)$-th column of
$E$, i.e., $E=(e_{n+1}, \cdots, e_{n+m}, e_1, \cdots, e_n)$.  It
follows from $dE=E\o_1^I$ that $e_{n+1}$ is an n-dimensional
submanifold in $S^{n+m-1}\subset R^{n+m}$ such that  $e_{n+1},
e_{n+2}, \cdots, e_{n+m}$ is a parallel normal frame and $(x_1,
\cdots, x_n)$ is a spherical line of curvature coordinates such that
the two fundamental forms are $$I=\sum_{i=1}^n a_{1i}^2 dx_i^2,
\quad \II=\sum_{i=1, \a=n+1}^{n, m+n} a_{1i} a_{\a-n, i} dx_i^2
e_\a.$$ We state this and the converse in the Theorem below.

\refclaim[aa] Theorem. (Flat n-submanifolds in $S^{m+n-1}$).  Let
 $X$ be a local isometric immersion of flat n-dimensional submanifold in
$S^{n+m-1}$ with flat, non-degenerate normal
bundle, and $m\geq n$.  Let $e_{n+1}=X$, and fix a local parallel
normal frame   $e_{n+2}, \cdots, e_{n+m}$. Then:
\item {(i)} There exist spherical  line of curvature coordinates
$(x_1, ..., x_n)$ and a smooth $\cm_{m\times n}$-valued map $A_1$
 such that $A_1^tA_1=I$ and the first and second
fundamental forms of $X$ are given by \refeq[ak]$$
I=\sum_{i=1}^na_{1i}^2\,dx_i^2,\qquad
           \II=\sum_{i=1}^n\sum_{j=2}^{m} a_{1i}a_{ji}dx_i^2\,e_{n+j}.$$
\item {(ii)} Let $f_{ij}=(a_{1i})_{x_j}/a_{1j}$ for $i\not=j$, $f_{ii}=0$
for all $1\leq i, j\leq n$, and  $F=(f_{ij})\in gl_\ast(n)$.
 Then the Gauss, Codazzi and Ricci equations
for the immersion $X$ is the $G_{m,n}$-system I
\refcn{} for  $(A_1, F)$.  In other words, $(A_1, F)$ is a
solution of \refcn{}.
\item {(iii)} Let $e_i={1\over a_{1i}} X_{x_i}$, $g=(e_{n+1}, \cdots, e_{m+n},
e_1, \cdots, e_n)\in O(m+n)$. Then
\refeq[cq]$$g^{-1}dg=\sum_i\pmatrix{0& -A_1C_i\cr C_i A_1^t& -F^t
C_i +C_i F\cr}dx_i,$$ which is equal to the Lax connection \refbo{}
$\o_\l^I$  of equation \refcn{} at $\l=1$.
\item {(iv)} Conversely, if $(A_1, F)$ is a solution of \refcn{}, then
system \refcq{} is solvable.  Let $g$ be a solution of \refcq{},
and $X$ the first column $g$. If all entries of the first row of
$A_1$ are non-zero, then $X$ is an isometric immersion of flat
$n$-submanifolds in $S^{n+m-1}$ with flat and non-degenerate
normal bundle such that the two fundamental forms are as in (i),
where $A_1=(a_{ij})$.

\proof

(i) can be proved the same way as for isometric immersions of
n-submanifolds in $R^{2n-1}$ with constant sectional curvature
$-1$ (cf. [Ca], [M]).

From (i), we have $$\omega_i=a_{1i}dx_i, \quad \w_{i, n+j}=
a_{ji}dx_i, \quad \w_{\a\b}=0.$$ By Corollary \refab{}, we get
$\omega_{ij} =f_{ij}dx_i-f_{ji}dx_j$.  So $g^{-1}dg= (\w_{AB})$,
and (ii), (iii) follow.

(iv) follows from the Fundamental Theorem of submanifolds in
Euclidean space.  \qed

\ms
 When $m=n$,  the above theorem was proved by Tenenblat in [Ten], where equation
\refcn{} is called the {\it generalized wave equation}. \ms

\ss
 If we use a different parallel frame $\tilde
e_{n+2}, \cdots, \tilde e_{n+m}$ for the immersion $X$ in the
above Theorem, then there exists a constant matrix $p=(p_{ij})\in
O(m-1)$ such that $\tilde e_{n+i}= \sum_{j=2}^m p_{ij}e_{n+j}$ for
$2\leq i\leq m$.  The solution of \refcn{} given by $X$ and
parallel frame $\tilde e_\a$ is $(\tilde pA_1, F)$, where $\tilde
p=\pmatrix{1&0\cr 0&p\cr}$ and $p=(p_{ij})$.

\ms

Suppose $(A,F)$ is a solution of \refcn{}, and $p\in O(m)$.  By
Corollary \refan{}, $(pA, F)$ is also a solution of \refcn{}.  As
a consequence of Theorem \refaa{} (iv), we get

\refclaim[ap] Corollary.  Let $X, e_{n+2}, \cdots, e_{n+m}$, and
$(A_1, F)$ be as in Theorem \refaa{}, and $c=(c_1, \cdots, c_m)$ a
unit vector.  If all components of $cA_1$ never vanish, then
$Y=c_1X+c_2 e_{n+2} +\cdots + c_me_{n+m}$ is again an immersion of
a flat n-submanifold in $S^{n+m-1}$ with flat and non-degenerate
normal bundle.

\ms


\refclaim[bt] Theorem.  (Flat n-submanifolds in $R^{n+m}$). Let $n\leq m$, and 
$X$ a local isometric immersion of flat  n-dimensional submanifold in $R^{m+n}$ with flat
and non-degenerate normal bundle.  Fix a local parallel normal frame $e_{n+1},
\cdots, e_{n+m}$. Then:
\item{(i)} There exist a line of curvature coordinates $(x_1, ..., x_n)$, a
smooth $\cm_{m\times n}$-valued map $A_1=(a_{ij})$ and a smooth $\cm_{n\times
1}$-valued map
$b=(b_1,..., b_n)$ such that $A_1^tA_1=I$ and the
first and second fundamental forms of $X$ are given by
\refeq[bu]$$I=\sum_i b_i^2dx_i^2,\ \ \II=\sum_{i=1}^n\sum_{j=1}^m
b_ia_{ji}dx_i^2 e_{n+j}.$$
\item{(ii)} Let $f_{ij}=( b_i)_{ x_j}/b_j$ for $i\ne j, f_{ii}=0,$ and
$F=(f_{ij})$. Then $(A_1,F)$ is a solution of the $G_{m,n}$-system I \refcn{}. 
Moreover, if $$e_i := {1\over b_i} X_{x_i}, \quad g:= (e_{n+1}, \cdots,
e_{n+m}, e_1, \cdots, e_n),$$ then $$g^{-1}dg =\sum_i \pmatrix{0&
-A_1C_i\cr C_i A_1^t & -F^t C_i + C_i F\cr}dx_i,$$ which is  the Lax connection
$\o_\l^I$ defined by \refbo{} for equation
\refcn{} at $\l=1$.
\item{(iii)} Conversely, if $(A_1,F)$ is a solution of \refcn{} and $b_1,
\cdots, b_n$ satisfy 
\refeq[dw]$$(b_i)_{x_j}= f_{ij} b_j, \quad
1\leq i\not=j\leq n,$$ then there exists a local isometric
immersion of $R^n$ in $R^{n+m}$ such that the two fundamental
forms are given  by \refbu{}, where $A_1=(a_{ij})$ and
$F=(f_{ij})$.
\item{(iv)} $X$ lies in a hypersphere of radius 1 if and only if $b=vA_1$
for some constant unit vector $v\in \cm_{1\times m}$.
\enditem

\proof Statements (i), (ii) and (iii) follows from an argument similar
to those  for Theorem \refaa{}. To prove (iv), we assume $\|
X-X_0\|=1$ for some constant vector $ X_0$. Then $X-X_0$ is a
parallel normal field, which implies that there exists constant
unit vector $v=(v_1, \cdots, v_m)$ such that $X-X_0= \sum_i v_i
e_{n+i}$.  Hence $$d(X-X_0)=dX =\sum_i v_i de_{n+i} = \sum_i v_i
\w_{j,n+i}e_j= \sum_i v_ia_{ij} dx_j e_j.$$ But $dX=\sum_j b_j
dx_je_j$.  So $b=vA_1$.  The converse can be proved by reversing
the argument.   \qed

 \ms
\refpar[dx] Remark.   System
\refdw{} was studied by Darboux [Da]. It can be solved if and only
if \refeq[eg]$$(f_{ij})_{x_k}= f_{ik}f_{kj}, \quad i, j, k \,\,
{\rm distinct\/}.$$ This condition  is obtained by equating
$(b_i)_{x_jx_k}= (b_i)_{x_kx_j}$: $$\eqalign{(b_i)_{x_jx_k}&=
(f_{ij}b_j)_{x_k}= (f_{ij})_{x_k} b_j + f_{ij}f_{jk} b_k\cr
&=(f_{ik}b_k)_{x_j} = (f_{ik})_{x_j} b_k + f_{ik}f_{kj} b_j.\cr}$$
Moreover, the space of solutions of \refdw{} depends on $n$
arbitrary smooth functions of one variable.  In fact, if
$\xi(t)=(\xi_1(t), \cdots, \xi_n(t))$ is a curve from $(-\e,\e)$
to $R^n$ such that $\xi_i'(t)$ never vanishes for all $1\leq i\leq
n$, then given any smooth maps $b_1^0, \cdots, b_n^0:(-\e, \e)\to
R$ there exists a unique solution $(b_1, \cdots, b_n)$ of \refdw{}
such that $b_i(\xi(t))= b_i^0(t)$.   The compatibility condition
\refeg{} is the third equation of \refcn{}.  So $(b_1, \cdots,
b_n)$ in Theorem \refbt{} (iii) always exists.

\ms
We can use a proof similar to that of the previous theorem to deduce:

\refclaim[cb] Theorem. (Local isometric immersion of $S^n$ in $S^{n+m}$).
Suppose
$X$ is a local isometric immersion of $S^n$ in $S^{n+m}$
with flat and non-degenerate normal bundle, and $\{e_\alpha\}$ is a parallel normal
frame.
 Then:
\item{(i)} There exist  a local coordinate system $(x_1,\dots,x_n)$, a smooth
$\cm^0_{m\times n}$-valued map
$A_1=(a_{ij})$  and a smooth $\cm_{n\times 1}$-valued map $b=(b_1,\cdots,
b_n)^t$ such that
 the two fundamental forms of $X$ are given by
\refeq[cd]$$I=\sum_i b_i^2\,dx_i^2,\ \ \II=\sum_{i=1, j=1}^{n,
m} a_{ji} b_idx_i^2 e_{n+j}.$$
\item{(ii)} Let $f_{ij}=(b_i)_{x_j}/b_j$ if $i\not=j$, $f_{ii}=0$ for all
$1\leq i\leq n$,  and $F=(f_{ij})$. Then $(A_1,F,b)$ is a solution
of the partial $G_{m,n+1}$-system I \refcl{}.
\item {(iii)} Let $e_i= {1\over b_i}X_{x_i}$,  and
$g=(e_{n+1}, \cdots, e_{n+m}, e_1, \cdots, e_n, X)$.    Then
\refeq[cy]$$ g^{-1}dg =\sum_i\pmatrix{0& -A_1C_i & 0\cr C_iA_1^t &
-F^tC_i+C_iF & C_i b\cr 0& -b^tC_i &0\cr} dx_i,$$ which is  the Lax connection
$\Theta_\l$ defined in \refcs{}  for system \refcl{} at $\l=1$.
\item {(iv)} If $(A_1,F, b)$ is a solution of \refcl{}, then system  \refcy{} is solvable. 
Let $g$ be a solution of \refcy{}, and let $X$ denote  the last column of $g$.
Then $X$ is a local isometric immersion of $S^n$ in $S^{n+m}$.

\bs Next we study submanifolds associated to  the
$G_{m,n}$- system II \refcu{}.  We will show that
each solution of \refcu{} gives rise to an $O(n)$-family of
submanifolds with common line of curvature coordinates.  In order to
simplify the notation, we make the following definition:

\refpar[aq] Definition.  Let $m>n$,  $\co$ a domain in $R^n$, and $X_i:\co\to R^m$
an immersion with flat and non-degenerate normal
bundle for $1\leq i\leq n$.   $(X_1, \cdots, X_n)$ is called a
 {\it  n-tuple in $R^m$ of type $O(n)$\/}
($O(n-1,1)$ resp.) if  
\i {(i)} the normal plane of $X_i(x)$ is parallel to the normal plane of $X_j(x)$ for all
$1\leq i,j\leq n$ and $x\in \co$, 
\i {(ii)} there exists a common parallel normal frame
$\{e_\a\}_{\a=n+1}^m$,
\i {(iii)} $x\in \co$ is a spherical (hyperbolic resp.) line of
curvature coordinate system (cf. Definition
\refal{}) with respect to $e_\a$ for each $X_i$ such that the fundamental forms for
$X_j$ are
\refeq[ar]$$\eqalign{I_j&= \sum_{i=1}^n b_{ji}^2 dx_i^2,\cr \II_j&=
\sum_{i=1, k=1}^{n,m-n} b_{ji} g_{k i} dx_i^2 e_{n+k}\cr}$$
for some $O(n)$-valued (($O(n-1,1)$- resp.)  map $B=(b_{ij})$ and a
$\cm_{(m-n)\times n}$-valued map
$G=(g_{ij})$.

\ss
We note that an $n$-tuple in $R^m$ of type $O(n)$ or $O(n-1,1)$ is an n-tuple of 
parametrized $n$-dimensional submanifolds in $R^m$ and the parametrization is a
spherical or hyperbolic line of curvature coordinates.

\refclaim[ac] Theorem.  Let $(X_1, \cdots, X_n)$ be an  n-tuple  in
$R^m$ of type
$O(n)$, $e_{n+1}, \cdots$, $e_m$  common parallel normal frame, and $(x_1, \cdots,
x_n)$ a common spherical line of curvature coordinates for all $X_j$'s  such
that the two fundamental forms $I_j, \II_j$ for $X_j$ are given by
\refar{}.  Set $f_{ij}=(b_{1j})_{x_i}/b_{1i}$ if $i\not=j$,
$f_{ii}=0$, and $F=(f_{ij})$. If all entries of $G$ are nonzero,
then $(F,G,B)$ is a solution of \refcu{}, the $G_{m,n}$-system II.

\proof  It follows from the definition of 
n-tuples that $$\w^{(j)}_1= b_{j1}dx_1,\,\, \w_2^{(j)} = b_{j2}
dx_2,\,\, \cdots, \,\,\w^{(j)}_n = b_{jn} dx_n$$ is a dual 1-frame
for $X_j$, and $\w_{i,n+k}^{(j)}=g_{ki}dx_i$ for each $X_j$. Note
that $\w_{i,n+k}^{(j)}$ is independent of $j$.   By Corollary
\refab{}, the Levi-Civita connection 1-form for the metric $I_j$
is $$\w_{ik}^{(j)} = -f^{(j)}_{ik}dx_k + f^{(j)}_{ki} dx_i, \quad
{\rm where\,\,} f^{(j)}_{ik}={(b_{jk})_{x_i}\over b_{ji}}.$$ Since
$$d\w^{(j)}_{i,n+k}=-\sum_{r=1}^n \w^{(j)}_{ir}\wedge
\w^{(j)}_{r,n+k}$$ for $1\leq k\leq m-n$ and $g_{k1}, \cdots,
g_{kn}$ are non-zero,  Cartan's Lemma and Corollary \refab{} imply
that $$f^{(j)}_{ir}= {(g_{kr})_{x_i}\over g_{ki}},$$ which is
independent of $j$.  Hence $$\w^{(j)}_{ij}= \w^{(1)}_{ij} =
-f_{ij}dx_j+f_{ji}dx_i.$$ The structure equation,
Gauss-Codazzi-Ricci equations for $X_1, \cdots, X_n$ imply that
$(F,G,B)$ is a solution of \refcu{}.  \qed

\ms The converse is also true:

\refclaim[dc] Theorem. (n-tuples in $R^m$ of
type $O(n)$). Suppose
$(F,G,B):R^n\to gl_\ast(n)\times
\cm_{(m-n)\times n}\times O(n)$ is a solution of the
$G_{m,n}$-system II \refcu{}.  Let
$$F=(f_{ij}),\quad  G=(g_{ij}), \quad B=(b_{ij}).$$ Then:
\item {(i)}
$$\tau= \sum_{i=1}^n \pmatrix{-FC_i+C_iF^t & C_iG^t\cr -GC_i &
0\cr}dx_i$$ is a  flat $o(m)$-valued connection 1-form.  Hence
there exists $A:R^n\to O(m)$ such that \refeq[az]$$A^{-1}dA =\tau
= \sum_{i=1}^n \pmatrix{-FC_i+C_iF^t & C_iG^t\cr -GC_i &
0\cr}dx_i.$$
\item {(ii)} Write $A=(A_1,A_2)$ with $A_1\in \cm_{m\times n}$ and $A_2\in
\cm_{m\times (m-n)}$. Then $$\sum_{i=1}^n A_1C_i B^t dx_i$$ is
exact.  So there exists a map $X:R^n\to \cm_{m\times n}$ such that
\refeq[di]$$dX= -\sum_{i=1}^n A_1C_i B^t dx_i.$$
\item {(iii)} Suppose all the entries of $B$  are
non-zero. Let  $X_j:R^n\to R^m$ denote the  $j$-th  column of $X$
(solution of \refdi{}), and $e_i$ denote the $i$-th column of $A$.
Then $(X_1, \cdots, X_n)$ is an  n-tuple in
$R^m$ of type $O(n)$.  In fact,
\itemitem {(1)} $e_1, \cdots, e_n$ are tangent to $X_j$ for all $1\leq
j\leq n$; so the tangent planes of $X_1, \cdots, X_n$ are
parallel, \ii {(2)} $\{e_{n+1}, \cdots, e_m\}$ is a parallel normal
frame for each $X_j$, \ii {(3)} the two fundamental forms for the
immersion $X_j$ are \refeq[cz]$$\eqalign{I_j &= \sum_{i=1}^n
b_{ji}^2 dx_i^2,\cr \II_j &= -\sum_{i=1, k=1}^{n, m-n} g_{ki}b_{ji}
dx_i^2 e_{n+k}.\cr}$$

\proof By Proposition \refdf{}, $\o_\l^{\II}$ defined by \refbp{}
is flat for all $\l\in C$.  In particular, $$\o_0^{\II} =
\pmatrix{\tau &0\cr 0&0\cr}$$ is flat.

The gauge transformation of $\o_\l^{\II}$ by $\pmatrix{A&0\cr
0&I\cr}$ is \refeq[dj]$$\pmatrix{A&0\cr 0&I\cr}\ast \o_\l^{\II} =
\sum_{i=1}^n \pmatrix{0&-\l A_1C_iB^t\cr  \l BC_i A_1^t & 0\cr},$$
which is flat for all $\l\in C$.   It follows from Proposition
\refdg{} that $\sum_{i=1}^n A_1C_i B^t dx_i$ is exact.  This
proves (ii).

Equate the $j$-th column of equation \refdi{} to get $$dX_j=
-\sum_{k=1}^n b_{jk}e_kdx_k.$$  So $I_j= \sum_k b_{jk}^2 dx_k^2$.
The rest of (iii) follows from the fact that $A^{-1} dA= \tau$.
 \qed

\refclaim[ad] Corollary.  Let $(X_1, \cdots, X_n)$ be an  n-tuple in $R^m$ of type
$O(n)$, and $p\in O(m)$ and
$q\in O(n)$ constant matrices.  Then $p(X_1, \cdots, X_n)q$ is
also an  n-tuple in $R^m$ of type $O(n)$.

\proof  We call a local orthonormal frame $A=(e_1, \cdots, e_m)$
an adapted frame for the  n-tuple $(X_1,
\cdots, X_n)$ of type $O(n)$  if  $e_1, \cdots, e_n$ are common principal
curvature directions and $e_{n+1}, \cdots, e_m$ are a common
parallel normal frame for each $X_j$.  By Theorem \refac{}, there
exist $G$ and $B$ such that $(F,G,B)$ is a solution of system
\refcu{}.  A direct computation shows that $Xq$ is an  n-tuple with $Aq$ as an adapted
frame and the corresponding solution of \refcu{} is the same $(F,G,B)$.

Note that if $A$ is a solution of \refaz{}, then so is  $pA$.   By
Theorem \refdc{} (ii), $pX$ is an  n-tuple of type $O(n)$ in $R^m$.
\qed

\ms The immersion $X_j$ in Theorem \refdc{} can be obtained either
by solving the system \refdi{} using integration or by an analogue of Sym's formula
(cf. [Sy]) below.

\refclaim[iv] Proposition.  Suppose $\o(x,\l)= \a_1(x)\l + \a_0(x)$ is a
flat $\cg$-valued connection 1-form on $x\in R^n$, and $E(x,\l)$ is a
trivialization of $\o(x,\l)$, i.e., $E^{-1}dE=\o$.  Set 
$$Y(x)= {\p E\over \p \l}(x,0)E^{-1}(x,0).$$
Then $dY= E(x,0)\a_1(x) E(x,0)^{-1}$. 

\proof A direct computation gives
$$\eqalign{dY&=\left({\p\over \p \l} (dE)\right)E^{-1} - {\p E\over
\p\l} E^{-1}dE E^{-1}\bigg|_{\l=0}\cr
&= \left({\p\over \p\l} (E\o)\right)E^{-1}\bigg|_{\l=0} -{\p E\over \p
\l}(x,0)\o(x,0) E(x,0)^{-1}\cr
&= E(x,0)\a_1(x) E(x,0)^{-1}. \qed\cr}$$

\refclaim[pa] Corollary. Let $E(x,\l)$ be a frame for a solution
$\xi$ of the $G_{m,n}$-system \refub{}, and $$Y(x)={\p E\over
\p
\l}(x,0)E^{-1}(x,0).$$ Then:
\item {(i)} $Y=\pmatrix{0&X\cr -X^t& 0\cr}$ for some $X\in \cm_{m\times n}$.
\item {(ii)} $X=(X_1, \cdots, X_n)$ is an n-tuple in $R^m$ of type $O(n)$. 
\item {(iii)}
 $dX= -A_1\d B^t$, where $\d=\diag(dx_1, \cdots, dx_n)$. In other
words, $X$ satisfies \refdi{}

\proof Since $E$ satisfies the reality condition, $E(x,0)\in O(m)\times O(n)$.
Write $E(x,0)=\pmatrix{A(x)&0\cr 0&B(x)\cr}$ and set 
$$\d=\diag(dx_1, \cdots, dx_n), \quad \b=\pmatrix{\d\cr 0\cr}.$$
  It
follows from Proposition \refiv{} and the fact that the Lax connection
 of \refub{} is 
$\o(x,\l)=\a_0\l + \a_1$ with
$$\a_1=\pmatrix{0& -A_1\d
B^t\cr B\d A_1^t&0\cr},$$ that 
$$dY = E(x,0)\pmatrix{0& -\b\cr \b & 0\cr} E(x,0)^{-1} 
= \pmatrix{0& -A_1\d B^t\cr B \d A_1^t &0\cr}.$$
So $dX = - A_1 \d B^t$, i.e., $X$ solves \refdi{}.  
\qed

\ms

We end this section by studying the $G_{m,2}$-system II \refcu{}.

\refclaim[ds] Proposition.  The $G_{m,2}$-system
II \refcu{} is the Gauss-Codazzi equations for a surface in $R^m$
that admits spherical line of curvature coordinates.

\proof Suppose $M$ is a surface in $R^m$, which admits a spherical
line of curvature coordinate system $(x_1, x_2)$ with respect to a
parallel normal frame $e_3, \cdots, e_m$. Then there exist  a
function $u$ and a $\cm_{(m-2)\times 2}$ -valued map $G=(g_{ij})$
such that $$\eqalign{I&= \cos^2 u \, dx_1^2 + \sin^2 u \, dx_2^2,
\cr \II&= \sum_{j=1}^{m-2} (g_{j1}\cos u \,dx_1^2 + g_{j2} \sin u
\,dx_2^2)e_{2+j}. \cr}$$  The Gauss-Codazzi-Ricci equations are
the following system for $u, g_{j1}, g_{j2}$:
\refeq[db]$$\cases{u_{x_1x_1}-u_{x_2x_2} +
\sum_{i=1}^{m-2}g_{i1}g_{i2}=0,&\cr (g_{i2})_{x_1} =g_{i1}
u_{x_1},&\cr (g_{i1})_{x_2}= - g_{i2} u_{x_2}.&\cr}$$ Let
\refeq[dk]$$B=\pmatrix{\cos u &\sin u\cr -\sin u &\cos u\cr},
\quad F=\pmatrix{0&u_{x_1}\cr -u_{x_2}& 0\cr}, \quad G=(g_{ij}).$$
Then $(F,G,B)$ is a solution of the $G_{m,2}$
system II \refcu{}.

Conversely, if $(F,G,B)$ is a solution of the $G_{m,2}$-system II
\refcu{}, then  we may assume $B=\pmatrix{\cos u & \sin u\cr -\sin u
&\cos u\cr}$.  If $\sin u\cos u\not=0$, then by the third equation of
\refcu{}, we have $f_{12}= u_{x_1}$ and
$f_{21} = -u_{x_2}$, i.e., $F=\pmatrix{0& u_{x_1}\cr -u_{x_2} &
0\cr}$. Let $g_{ij}$ denote the $ij$-th entry of $G$.  Then system
\refcu{} is system \refdb{}. \qed

\refclaim[dy] Corollary. Let $X_1:\co\to R^m$ be an immersion and $(x,y)\in \co$ the
spherical line of curvature coordinates with respect to a
parallel normal frame $e_3, \cdots, e_m$. Then there exists an immersion
surface $X_2$ unique up to translation such that $(X_1,X_2)$ is a
 2-tuple in $R^m$ of type $O(2)$. Moreover,
the fundamental forms of $X_1, X_2$ are respectively given by
\refeq[da]$$\eqalign{&I_1=
\cos^2 u\, dx_1^2 + \sin^2 u\, dx_2^2,\quad \II_1= \sum_{j=1}^m(
g_{j1}\cos u \,dx_1^2 + g_{j2} \sin u \,dx_2^2) e_{2+j},\cr
&I_2=\sin^2 u \, dx_1^2 + \cos^2 u\, dx_2^2,\quad \II_2=
\sum_{j=1}^{m-2}(-g_{j1}\sin u \, dx_1^2  + g_{j2} \cos u\,
dx_2^2)e_{2+j}.\cr}$$

   It follows from the Gauss equation \refag{}
that $$K_1=-K_2 = {\sum_{j=1}^{m-2} g_{j1}g_{j2}\over \sin u\cos
u}.$$   So we have

\refclaim[aj] Corollary.  Let $\co$ be an open subset of $R^2$, $(X_1,X_2):\co\to
R^m\times R^m$  a  2-tuple in $R^m$ of type
$O(2)$, and $(x_1, x_2)\in \co$ spherical line of curvature coordinates.  Then
$$K_2(x)= -K_1(x),$$ where $K_1, K_2$ are respectively the Gaussian curvature of
$X_1, X_2$.

\refpar[ih] Definition.  If $(M_1, M_2)$ is a  2-tuple
 in $R^m$ of type $O(2)$ ($O(1,1)$ resp.)  with a parallel normal frame $e_3,
\cdots, e_m$ and spherical line of curvature coordinates, then we
call $M_2$ a {\it C-dual\/} of $M_1$.  Note that any two C-duals
of $M_1$ are differed by a translation.

\ms

\refpar[dz] Example. Recall that given a surface $M$ in $R^3$ with curvature
$-1$, there exist spherical line of curvature coordinates $x_1, x_2$ such
that the two fundamental forms are $$\eqalign{I&= \cos^2 u \,
dx_1^2 + \sin^2 u\, dx_2^2,\cr \II&= \sin u\cos u \,
(dx_1^2-dx_2^2),\cr}$$ and $u$ satisfies
$$u_{x_1x_1}-u_{x_2x_2} = \sin u\cos u. \eqno({\rm SGE})$$
 This implies that $(u, \sin u, -\cos u)$ is a solution of \refdb{}.  
Let $X(x_1,x_2)$ denote the immersion of $M$.  Then $(X,e_3)$ is a 2-tuple in $R^3$
of type $O(2)$, where
$e_3$ is the unit normal  of $M$, which is an open subset of $S^2$. 


\bs

\newsection Submanifolds associated to $G_{m,n}^1$-
systems.\par

In this section, we  describe
submanifolds whose fundamental equations are given by
$G_{m,n}^1$-systems.

Let $R^{k, 1}$ denote the Lorentz space equipped with the
non-degenerate bilinear form of index one: $$\langle x, y\rangle
_1=x_1y_1+\dots +x_ky_k-x_{k+1}y_{k+1}.$$  The moving frame
computation for submanifolds in $R^{k,1}$ can be carried out in a
similar way as for submanifolds in $R^{k+1}$ except that the
Levi-Civita connection $1$-form $(\w_{ij})$ of the flat Lorentzian
metric $\li , \ri_1$ is $o(k,1)$-valued.

Let $H^k$ denote the k-dimensional simply connected space form of
sectional curvature -1 (i.e., a {\it hyperbolic
k-space\/}). It is well-known that $$\{x\in R^{k,1}\n \langle x,
x\rangle_1= -1, x_{k+1}>0 \}$$ with the induced metric is isometric to $H^k$.
We need the following Proposition later, which can be proved using
a direct computation (cf. Chap. 2 of [PT]).

\refclaim[dl] Proposition.  Let  $v_0\in R^{k,1}$ be a constant
non-zero vector, $c\in R$ a constant, and $N_{v,c}$ the
hypersurface defined by $$N_{v_0,c}=\{x\in H^k\n  \langle x,
v_0\rangle_1= c\}.$$ Then $N_{v,c}$  is  a totally umbilic
hypersurface of $H^k$ with constant sectional curvature $${-\li
v_0, v_0\ri_1\over c^2+ \li v_0, v_0\ri_1}.$$

We use methods similar to those  in section 6 to find submanifolds whose
fundamental equations are the various $G_{m,n}^1$-systems.
 For the $G_{m,n}^1$-system I,  let $(\o_{ij})$ denote the connection 1-form
$\o^I_\l$ defined by \refbx{} at $\l=1$. Then $(\o_{ij})$ is a
$o(m+n,1)$-valued 1-form and 
$$\cases{ \o_{m+i, m+j} =  -\e_i\e_j
f_{ji} dx_j + f_{ij} dx_i,& $1\leq i, j\leq n+1$,\cr \o_{m+i, \a}=
a_{\a i} dx_i, & $1\leq i\leq n+1, 1\leq \a\leq m$, \cr \o_{\a\b}=
0, & $1\leq \a, \b\leq m$,\cr}$$ 
where $\e_1=\cdots =\e_n =-\e_{n+1}= 1$. 
Let $$\cases{\w_{ij}=
\o_{m+i,m+j},& $1\leq i, j\leq n+1$,\cr \w_{i\a} = \o_{m+i, \a},&
$1\leq i\leq n+1, 1\leq \a\leq m$, \cr \w_{\a\b}= \o_{\a\b}, &
$1\leq \a, \b\leq m$.\cr}$$ Then the flatness of $(\o_{ij})$ is
exactly the Gauss-Codazzi-Ricci equations for
 (n+1)-dimensional flat, Lorentzian submanifolds in $R^{n+m,1}$.  So we get

\refclaim[by] Theorem. (Local isometric immersions of $R^{n,1}$ in $R^{m+n,
1}$).  Let  $X$ be a local isometric immersion of $R^{n,1}$ in $R^{n+m,1}$ with flat and
non-degenerate normal bundle.  Fix a local parallel normal
frame
$e_{n+2},
\cdots, e_{n+m+1}$. Then:
\item{(i)} There exist a line of curvature coordinates $(x_1, ..., x_{n+1})$ and a
$\cm_{m \times (n+1)}$-valued map $A_1=(a_{ij})$  and a map
$b=(b_1,..., b_{n+1})$ such that $A_1^tA_1=I$ and
the first and second fundamental forms of $X$ are given by
\refeq[ca]$$I=\sum_{i=1}^{n+1} \e_ib_i^2dx_i^2,\ \
\II=\sum_{i=1}^{n+1}\sum_{j=2}^{m-n} b_ia_{ji}dx_i^2 e_{n+j},$$
where $\e_1=\cdots =\e_n = -\e_{n+1}=1$.
\item{(ii)} Let $f_{ij}=( b_i)_{ x_j}/b_j$ for $i\ne j, f_{ii}=0,$ and
$F=(f_{ij})$. Then $(A_1,F)$ is a solution of the $G_{m,n}^1$-system I \refcc{}.
\item{(iii)} Conversely, if $(A_1,F)$ is a solution of \refcc{} and $b_1,
\cdots, b_{n+1}$ satisfies $(b_i)_{x_j}= f_{ij} b_j$ for all
$i\not=j$, then there exists a local isometric immersion of
$R^{n,1}$ in $R^{m+n,1}$ such that the two fundamental forms are
given  by \refca{}, where $A_1=(a_{ij})$ and $F=(f_{ij})$.

Next we consider the partial $G_{m,n}^1$-system I of n variables, system
\refaaj{}.  Let $\o_{ij}$ denote the $ij$-th entry of $\tau^I_\l$,
 defined by \refaak{}, at $\l=1$.  So we have
$$\cases{\o_{\a\b}= 0, & $1\leq \a, \b\leq m$,\cr \o_{m+i, m+j}=
-f_{ji}dx_j + f_{ij}dx_i, & $1\leq i, j\leq n$, \cr \o_{m+i, \a}=
a_{\a i} dx_i, & $1\leq i\leq n, 1\leq \a\leq m$, \cr \o_{m+i,
m+n+1}= b_idx_i, & $1\leq i\leq n$.\cr}$$ Let $$\cases{\w_i =
\o_{m+i, m+n+1}, & $1\leq i\leq n$,\cr \w_{ij}= \o_{m+i, m+j}, &
if $1\leq i, j\leq n$,\cr \w_{i,n+j}= \o_{m+i, j},& if $1\leq
i\leq n, \quad1\leq j\leq m$,\cr \w_{n+i, n+j}=\o_{ij},& if $1\leq
i, j\leq m$. \cr}$$ Then system \refaaj{} (given by the flatness
of $\tau_1^I$) is exactly the Gauss-Codazzi-Ricci equations for
local isometric immersions of $H^n$ in $H^{n+m}\subset R^{n+m,1}$ with flat
and non-degenerate normal bundle.   We summarize this below.

\refclaim[bbb] Theorem. (Local isometric immersion of $H^n$ in $H^{n+m}$).  Let
$X$ be a local isometric immersion of $H^n$ in $H^{n+m}$ with flat
and non-degenerate normal bundle, and $\{e_\alpha\}$ a parallel normal frame.
 Then:
\item{(i)} There exist  a line of curvature coordinate system $(x_1,\dots,x_n)$, a 
$\cm_{m\times n}$-valued map 
$A_1=(a_{ij})$ satisfying $A_1^tA_1=I$,  and a $\cm_{n\times 1}$-valued map
$b=(b_1,\cdots, b_n)^t$ such that
 the two fundamental forms of $X$ are given by
\refeq[bbc]$$I=\sum_i b_i^2\,dx_i^2,\ \ \II=\sum_{i=1, j=1}^{n,
m-n} a_{ji} b_idx_i^2 e_{n+j}.$$
\item{(ii)} Let $f_{ij}=(b_i)_{x_j}/b_j$ if $i\not=j$, $f_{ii}=0$ for all
$1\leq i\leq n$,  and $F=(f_{ij})$. Then $(A_1,F,b)$ is a solution
of the partial $G_{m,n}^1$-system \refaaj{}.
\item {(iii)} Let $e_i= {1\over b_i}X_{x_i}$,  and
$g=(e_{n+1}, \cdots, e_{n+m}, e_1, \cdots, e_n, X)$.    Then
\refeq[bbd]$$ g^{-1}dg =\sum_i\pmatrix{0& -A_1C_i & 0\cr C_iA_1^t&
-F^tC_i+C_iF & C_i b\cr 0& b^tC_i &0\cr} dx_i,$$ which is equal to
the Lax connection $\tau^I_\l$ defined in \refaak{}  for system
\refaaj{} at $\l=1$.
\item {(iv)} If $(A_1,F, b)$ is a solution of \refaaj{}, then system
\refbbd{} is solvable.
\item {(v)}  Let $g$ be a solution of
\refbbd{}. Then the last column of $g$ is an isometric immersion of
a constant sectional curvature $-1$ submanifold of $H^{n+m}$.
\item{(vi)} $M^n$ lies in a totally umbilic hypersurface of $H^{n+m}$ if
and only if there is a constant vector $w\in \cm_{m\times 1}$ such
that
 $b=A_1^tw$. Moreover, $M^n$ lies in a flat totally umbilic hypersurface if and
only if  $\Vert w\Vert=1$.

\proof Statements (i)-(v) follows from standard submanifold
theory. To prove (vi), note that every umbilic hypersurface of
$H^{n+m}$ is the intersection of $H^{n+m}$ with a hyperplane (cf. [PT]),
i.e., it is of the form $$N_{v,c}=\{x\in H^{n+m}\n \li x,
v\ri_1=c\}$$ for some constant $v\in R^{n+m,1}$ and $c\in R$. By
Proposition \refdl{}, $N_{v,c}$ has constant sectional curvature
$$-{\li v, v\ri_1\over c^2+ \li v, v\ri_1}.$$
 Suppose the image of $X$ lies in $N_{v,c}$. Then
$<dX, v>_1=0$, which implies that $v$ is normal to $X$.  But
$v$ is constant, so $v$ is a parallel normal vector field of $M$
as a submanifold of $R^{n+m,1}$. Hence
\refeq[cf]$$v=\sum_{k=1}^{m} v_k e_{n+k} + v_{m+1}X$$ for some
constants $v_1, \cdots, v_{m+1}$. It follows from $\li X, v\ri_1 =
c$ and $\li X, X\ri_1 =-1$ that $v_{m+1}=-c$. Differentiate
\refcf{} to get $$ c\w_i = cb_idx_i= \sum_{k=1}^m
v_k\omega_{i,n+k}=\sum_{k=1}^m v_k a_{ki}dx_i,$$ 
so 
\refeq[eh]$$cb_i =\sum v_k a_{ki}.$$ Let $\hat v=(v_1, \cdots, v_m)^t$, and
$w=\hat v/c$. Then \refeh{} implies that $b=A_1^tw$. Since $\li v,
v\ri_1= \N \hat v\N^2 -c^2$, it follows from the formula for the
sectional curvature that $N_{v,c}$ is flat if and only if $\li v,
v\ri_1=0$, i.e., $\N \hat v\N^2=c^2$, or equivalently, $\N
w\N=1$. \qed

\ms  For $m=n$,  the above theorem was proved by Terng in [Te2].

\ms

\refclaim[cca] Theorem. ($(n+1)$-tuples in $R^m$ of type
$O(n,1)$).   Suppose
$m>n+1$, and
$$(F,G,B):R^{n+1}\to gl_\ast(n+1)\times \cm_{(m-n-1) \times (n+1)}\times O(n,1)$$ is
a solution of the $G_{m,n}^1$-system II
\refaad{}.  Let $F=(f_{ij})$,  $G=(g_{ij})$, and $B=(b_{ij})$.
Then:
\item {(i)}
\refeq[do]$$\w= \sum_{i=1}^{n+1} \pmatrix{-FC_i+C_iF^t & C_iG^t\cr
-GC_i & 0\cr}dx_i$$ is a  flat $o(m)$-valued connection 1-form.
Hence there exists $A:R^{n+1}\to O(m)$ such that
\refeq[dn]$$A^{-1}dA =\w = \sum_{i=1}^{n+1} \pmatrix{-FC_i+C_iF^t
& C_iG^t\cr -GC_i & 0\cr}dx_i.$$
\item {(ii)} Let $A=(A_1,A_2)$ with $A_1\in \cm_{m\times (n+1)}$. Then
$\sum_{i=1}^{n+1} A_1C_i JB^{-1} dx_i$ is exact, so there exists
a map $X:R^{n+1}\to \cm_{m\times (n+1)}$ such that
\refeq[ccb]$$dX= -\sum_{i=1}^{n+1} A_1C_i B^tJ dx_i.$$
\item {(iii)} Suppose all the entries of the $j$-th row of $B$  are
non-zero. Let  $X_j:R^{n+1}\to R^m$ denote the  $j$-th  column of
$X$, and $e_i$ denote the $i$-th column of $A$.  Then $(X_1,
\cdots, X_{n+1})$ is a  $(n+1)$-tuple in
$R^m$ of type $O(n,1)$, and the two fundamental forms for the
immersion $X_j$ are \refeq[ccc]$$\eqalign{I_j &= \sum_{i=1}^{n+1}
b_{ji}^2 dx_i^2,\cr \II_j &= - \sum_{i=1, k=1}^{n+1, m-n-1}
g_{ki}b_{ji} dx_i^2 e_{n+k+1}.\cr}$$
 \item {(iv)} Conversely, suppose  $(X_1, \cdots, X_{n+1})$ is an  $(n+1)$-tuple in
$R^m$ of type $O(n,1)$ such that the two fundamental forms are of the form
\refccc{} with respect to a common parallel normal frame $e_{n+2}, \cdots, e_m$ for
some
$B=(b_{ij})$ and $G=(g_{ij})$. Then $(F,G,B)$ is a solution of
\refaad{}, where $$F=(f_{ij}), \quad f_{ij}={ (b_{1j})_{x_i}\over
b_{1i}}.$$

\proof The proof is similar to that of Theorem \refdc{}.  We only give the proof of (iii)
here.   Let $\w_{ij}$ denote the $ij$-th entry of $\w$ (defined by
\refdo{}),   i.e., \refeq[ec]$$\w_{ij}=\cases{-f_{ij}dx_j +
f_{ji}dx_i, & $1\leq i,j\leq n+1$,\cr g_{j-n-1, i}dx_i, & $1\leq
i\leq n+1$ and $n+2\leq j\leq m$,\cr 0, & $n+2 \leq i, j\leq
m$.\cr}$$ Let $e_i$ denote the i-th column of $A$ (defined by
\refdn{}), i.e., $A=(e_1, \cdots, e_m)$. Since $dA=A\w$, we have
\refeq[eb]$$de_i=\sum_{j=1}^m \w_{ji} e_j, \quad 1\leq i\leq m.$$
By Theorem \refcca{} (ii), \refeq[ea]$$dX_r = -\sum_{i=1}^{n+1}
b_{ri} \e_r e_i dx_i,$$ where $X_r$ is the $r$-th column of $X$.
So $\{e_1, \cdots, e_{n+1}\}$ is a common orthonormal tangent
frame, and $$\w_1^{(r)}= -\e_r b_{r1}dx_1, \cdots, \w_{n+1}^{(r)}
= -\e_r b_{r, n+1}dx_{n+1}$$ is the dual coframe for $X_r$.  Since
$$\w_{i, n+1+j}^{(r)} = \li de_{n+1+j}, e_i\ri_1= \w_{i,n+j+1}=
g_{ji} dx_i,$$ $\{e_{n+2}, \cdots, e_m\}$ is a common parallel
normal frame and $(x_1, \cdots, x_{n+1})$ is a hyperbolic line of
curvature coordinate system for each $X_r$. \qed

\refclaim[iw] Corollary.  Let $(X_1, \cdots, X_{n+1})$ be a 
$(n+1)$-tuple in $R^m$ of type $O(n,1)$, and $p\in O(m)$, $q\in O(n,1)$ constant
matrices.  Then $p(X_1, \cdots, X_{n+1})q$ is also a 
$(n+1)$-tuple in $R^m$ of type $O(n,1)$.

\refclaim[jj] Corollary.  Let $v$ be a solution of the $G_{m,n}^1$-system \refaaa{}, 
$E$ a frame for $v$, and  $Y=({\p E\over \p \l} E^{-1})(x,0)$. Then:
\item {(i)} $Y=\pmatrix{0&X\cr -JX^t&0\cr}$ for some
$\cm_{m\times (n+1)}$-valued map
$X$ and $J=I_{n,1}=\diag(1, \cdots, 1, -1)$.  
\item {(ii)} $X=(X_1, \cdots, X_n)$ is an (n+1)-tuple in $R^m$ of type $O(n,1)$.  

\bs


\newsection $G_{m,1}^1$-systems and isothermic surfaces.\par

We study the relation between the $G_{m,1}^1$-system and isothermic
surfaces in $R^m$.  

\refclaim[dp] Proposition.  The $G_{m,1}^1$-system II \refaad{} is the
Gauss-Codazzi-Ricci equations for surfaces in $R^m$ admitting
hyperbolic line of curvature coordinates.

\proof Let $\co$ be a domain in $R^2$, $X:\co\to R^m$ an immersion 
with flat and non-degenerate normal bundle, and $(x_1,x_2)\in \co$ a hyperbolic line of
curvature coordinate system with respect to the parallel normal
frame $\{e_3,\cdots, e_m\}$. Then there exist $u$ and $g_{ki}$ for
$1\leq k\leq m-2$ and $i=1, 2$ such that the two fundamental forms
are \refeq[cg]$$\eqalign{I&=\cosh^2 u\, dx_1^2 + \sinh^2 u\,
dx_2^2, \cr \II&=\sum_{k=1}^{m-2} (g_{k1}\cosh u\, dx_1^2 + g_{k2}
\sinh u\, dx_2^2)e_{k+2}.\cr}$$  The Gauss-Codazzi-Ricci equations
for $X$ are \refeq[dq]$$\cases{u_{x_1x_1} + u_{x_2x_2} +
\sum_{k=1}^{m-2} g_{k1}g_{k2} =0, &\cr  (g_{k1})_{x_2}= u_{x_2}
g_{k2}, &\cr (g_{k2})_{x_1} = u_{x_1} g_{k1}.&\cr}$$ This is
exactly the $G_{m,1}^1$-system \refaad{} for
$(F,G,B)$, where \refeq[dr]$$F=\pmatrix{0&u_{x_1}\cr
u_{x_2}&0\cr}, \quad G=(g_{ki}),\quad B=\pmatrix{\cosh u&\sinh
u\cr \sinh u & \cosh u\cr} \in O(1,1).$$

Conversely, if $(F,G,B)$ is a solution of the
$G_{m,1}^1$-system \refaad{}, then since $B\in O(1,1)$ we
may assume $$B=\pmatrix{\cosh u& \sinh u\cr\sinh u & \cosh
u\cr}.$$  The third equation of \refaad{} implies that $f_{12}=
u_{x_1}$ and $f_{21}= u_{x_2}$, i.e., $(F,G,B)$ is of the form
\refdr{}.  Write equation \refaad{} for $(F,G,B)$ in terms of $u$
and $g_{ki}$ to get equation \refdq{}. This completes the proof.
\qed

\refclaim[ee] Corollary.  Let $\co$ be a domain of $R^2$, and $X_1:\co\to R^m$ an
immersion with flat normal bundle and $(x,y)\in \co$ a hyperbolic line of curvature
coordinate system with respect to a parallel normal frame $e_3, \cdots, e_m$. Then
there exists an immersion $X_2$ unique up to translation such that
$(X_1,X_2)$ is a  2-tuple in $R^m$ of type
$O(1,1)$. Moreover, the fundamental forms of $X_1, X_2$ are given respectively by
\refeq[da]$$\eqalign{&\cases{I_1= \cosh^2 u\, dx_1^2 + \sinh^2 u\,
dx_2^2,&\cr \II_1= \sum_{j=1}^{m-2}( g_{j1}\cosh u \,dx_1^2 +
g_{j2} \sinh u \,dx_2^2) e_{2+j},&\cr}\cr &\cases{I_2=\sinh^2 u \,
dx_1^2 + \cosh^2 u\, dx_2^2,&\cr \II_2=
\sum_{j=1}^{m-2}(-g_{j1}\sinh u \, dx_1^2  - g_{j2} \cosh u\,
dx_2^2)e_{2+j}.&\cr}\cr}$$ 

Let $K_1, K_2$ denote the Gaussian curvature of $X_1$ and $X_2$
respectively.  The Gauss equation implies that $$K_1= K_2 =
{\sum_{k=1}^{m-2} g_{k1}g_{k2}\over \sinh u\cosh u}.$$ Hence we get

\refclaim[ajj] Corollary.  Let $(X_1, X_2):\co\to R^m$ be a  2-tuple
 in $R^m$ of type $O(1,1)$ such that $(x_1, x_2)\in \co$ is a
hyperbolic line of curvature coordinate system. Then
$$K_2(x)= K_1(x)$$ where $K_1, K_2$ are respectively the Gaussian
curvature of $X_1, X_2$.

\ms

Next we describe the relation between 
2-tuple in $R^m$ of type $O(1,1)$ and isothermic surfaces in
$R^m$.  The following definition,
 given by Burstall in [Bu], generalizes the classical notion of isothermic
surfaces in $R^3$ ([Da]).

\refpar[dt] Definition. Let $\co$ be a domain in $R^2$. An immersion $X:\co\to R^m$
is called {\it isothermic\/} if it has flat normal bundle and the two fundamental forms are
of the form
$$I=e^{2u}(dx_1^2 + dx_2^2), \quad
\II=\sum_{k=1}^{m-2}e^{u}(g_{k1}dx_1^2 + g_{k2}dx_2^2)e_{2+k}$$
with respect to some parallel normal frame $e_3, \cdots, e_m$, or equivalent $(x_1,
x_2)\in \co$ is conformal and line of curvature coordinate system for $X$.

\ms 

A direct computation gives

\refclaim[el] Proposition.  The Gauss-Codazzi-Ricci equation for
isothermic surfaces in $R^m$ is \refdq{}.

It was first noted by Cie\'sli\'nski, Goldstein and Sym in [CGS] that equation for
isothermic surfaces in $R^3$ has a Lax pair.  There have been many papers ([Ci],
[BHPP], [Bu], [HP], [HMN]) using techniques from
solition theory to study isothermic surfaces in 3-space.

Next we give a  simple relation between isothermic surfaces in $R^m$ and 2-tuples
in $R^m$  of type $O(1,1)$.

\refclaim[ed] Proposition.  Suppose  $(X_1, X_2)$ is a  2-tuple in $R^m$ of type
$O(1,1)$.  Let $Y_1= X_1-X_2$ and $Y_2=X_1+X_2$. Then both $Y_1$ and $Y_2$
are isothermic.

\proof Let $(u, (g_{ij}))$ be the solution of \refdq{} associated
to $(X_1, X_2)$.   We use the same notations as in Theorem
\refcca{}.   Write $X=(X_1,X_2)$. By \refea{}, we get
$$\cases{dX_1= -\cosh u\ dx_1 e_1 - \sinh u\ dx_2 e_2,&\cr dX_2=
\sinh u \ dx_1e_1 +\cosh u \ dx_2 e_2. &\cr }$$ Note that
$$F=\pmatrix{0& u_{x_1}\cr u_{x_2}&0\cr}, \quad \w_{12}=
u_{x_2}dx_1- u_{x_1}dx_2, \quad \w_{i, 2+j}= g_{ji}dx_i.$$
 Use \refeb{} to get
$$\cases{\nabla e_1= \w_{21}e_2 = (u_{x_2}dx_1 - u_{x_1} dx_2)e_2,
&\cr \nabla e_2= \w_{12}e_1 = -(u_{x_2}dx_1 - u_{x_1} dx_2) e_1,
&\cr de_{k+2}= \w_{1,k+2}e_1 + \w_{2,k+2}e_2 = g_{k1}dx_1 e_1 +
g_{k2} dx_2 e_2.&\cr}$$ Compute directly to get $$dY_1= dX_1-dX_2
= -(\cosh u+\sinh u) (dx_1e_1+dx_2e_2).$$ So the induced metric
for $Y_1$ is $e^{2u}(dx_1^2 + dx_2^2)$. Similarly, we can get
explicit formula for $dY_2$ and $de_k$ for $3\leq k\leq m$. This
computation implies that the two fundamental forms for $Y_1$ and
$Y_2$ are \refeq[ej]$$\cases{I_1= e^{2u}(dx_1^2 + dx_2^2), &\cr
\II_1 = -e^u\sum_{k=1}^{m-2}(g_{k1}dx^2_1+ g_{k2}dx_2^2)
e_{2+k},&\cr}$$ 
\refeq[ei]$$\cases{I_2=
e^{-2u}(dx_1^2+dx_2^2),&\cr \II_2=
e^{-u}\sum_{k=1}^{m-2}(-g_{k1}dx^2_1 + g_{k2}dx^2_2)e_{2+k}.
&\cr}$$ Hence $Y_1$ and $Y_2$ are isothermic surfaces. \qed

\ms

We call $(Y_1, Y_2)$ an {\it isothermic pair\/}.  If $(u,
g_{k1}, g_{k2})$ is the solution of \refdq{} corresponding to
$Y_1$, then $(-u, g_{k1}, -g_{k2})$ is the solution of \refdq{}
corresponding to $Y_2$.  
 For $m=3$,
the transformation from $Y_1$ to $Y_2$ is the classical Christoffel
transformation.  Burstall generalizes this result to arbitrary $m$ in
[Bu]. In fact, he developed in [Bu] a beautiful theory of isothermic surfaces in $R^m$
and explained its relation to conformal geometry, Clifford algebra and the curved flat
system associated to $G_{m,1}^1$.   His result motivated us
to study the $G_{m,n}$-system II and $G_{m,n}^1$-system II for general $m, n$. 

Next we describe the
$o(m+1,1)$ model used by Burstall, Jeromin, Pedit and Pinkall
[BHPP] for $m=3$ and by Burstall [Bu] for general $m\geq 3$. Let
$$\eqalign{&\s_{k,1}=\pmatrix{I_{k-1}&0\cr 0&
\s_{1,1}\cr}, \qquad r_{k,1}=\pmatrix{I_{k-1}&0\cr 0& r_{1,1}\cr},
\cr &{\rm where\,\,} \s_{1,1}=\pmatrix{0&1\cr 1&0\cr}, \quad
r_{1,1}={1\over \sqrt{2}}\pmatrix{1&1\cr -1&1\cr}.}$$ Note that
the bilinear form defined by $\s_{k,1}$ and $I_{k,1}$ are
isometric, and $$r_{k,1}^t \s_{k,1} r_{k,1} = I_{k,1}.$$
  Let
$$O_\s(k, 1)=\{\xi\in GL(k+1)\n \xi^t\s\xi= \s\}.$$ Then
$\phi_{k,1}(g)= r_{k,1}^tgr_{k,1}$ is a group isomorphism from
$O(k,1)$ to $O_\s(k,1)$. If we replace $\o_\l^{\II}$ by $\W_\l^{\II}
= (\phi_{m+1,1})_\ast(\o_\l^{\II})$ in Theorem \refcca{}, then the
integration in Theorem \refcca{} (ii) gives an isothermic pair
$Y=(Y_1, Y_2)$.

\ms

  The following proposition, which states that a CMC surface in $R^3$ without umbilic
points admits isothermic coordinates,  is well-known (for a proof see [PT]).

\refclaim[eo] Proposition.  
Let $M$ be an immersed  surface
without umbilic points in space form $N^3(c)$, and $k_1, k_2$ the principal
curvature functions.  If  $M$ has constant mean curvature $H$ and $k_1>k_2$, then
there exists a conformal line of curvature coordinates $x,y$ such that 
$$I= {2\over k_1-k_2} \left(dx^2+dy^2\right), \quad \II=\left(1+{H\over
k_1-k_2}\right)dx^2 +\left(-1+{H\over k_1-k_2}\right)dy^2.$$
Moreover, if we write $e^{2u}={2\over k_1-k_2}$, then the Gauss-Codazzi
equation for
$M$ is 
\refeq[dua]$$u_{xx}+u_{yy}= e^{-2u} - \left({H^2\over 4}+c\right)e^{2u}.$$

\refclaim[ke] Corollary.  With the same notation as in Proposition \refeo{}.  
If $c= -H^2/4$, then equation
\refdua{} becomes
\refeq[du]$$u_{xx}+u_{yy}= e^{-2u}.$$
Moreover,  $(u, e^{-u}, -e^{-u})$ is a solution of \refdq{} if and
only if $u$ is a solution of \refdu{}. 

\ms

Proposition \refeo{} has a natural generalization to surfaces in $R^m$ observed by
Burstall [Bu].  First we recall a  definition due to Chen [C]:

\refpar[kc] Definition.  A surface $M$ in $R^m$ with flat normal bundle is called a {\it
generalized H-surface\/} if there exists a unit parallel normal field $v$ such that $\li H,
v\ri =$ constant, where $H$ is the mean curvature vector field.  

\ss
Note that a generalized $H$-surface in $R^3$ is a CMC surface.  

\ms
The following analogue of Proposition \refeo{} can be proved exactly the same way
for generalized H-surfaces  in $R^m$.

\refclaim[jz] Proposition ([Bu]).   
Let $M$ be a generalized  $H$-surface in $R^m$, $H$
the mean curvature vector, and $(e_3, \cdots, e_m)$ a parallel normal field such
that $\li H, e_3\ri $ $= c $ is a constant.  If the shape operator
$A_{e_3}$ has two distinct eigenvalues $k_1>k_2$, then there exists isothermic
coordinates $x,y$ on $M$ such that 
$$\eqalign{I&={2\over k_1-k_2} \left(dx^2+dy^2\right), \cr
\II_3&=\li \II, e_3\ri=\left(1+{c\over k_1-k_2}\right)dx^2 +\left(-1 + {c\over
k_1-k_2}\right)dy^2,\cr
\II_\a&=\li \II, e_\a\ri = g_{\a 1} dx^2 + g_{\a 2}dy^2, \quad 3< \a\leq m\cr}$$
for some $g_{\a i}$.  

We call the isothermic coordinates $x, y$ in the above Proposition {\it the
canonical isothermic coordinates\/} for the generalized $H$-surface $M$. 

\ms
The following Proposition generalizes the classical Bonnet transformation for CMC
surfaces in
$R^3$. The proof follows from a direct computation as in the classical case. 

\refclaim[kd] Proposition ([Bu]).  Let $\co$ be a domain in $R^2$, $Y:\co\to
R^m$ an isothermic immersion of a generalized
$H$-surface in $R^m$,  $(x_1,x_2)\in \co$ the canonical isothermic coordinates, and
$e_3$ the parallel unit normal field such that $\li H, e_3\ri= 2c$ is a constant. Then:
\i {(i)}  If $c\not=0$, then $(Y, Y- {1\over c} e_3)$ is an isothermic pair.
\i {(ii)} If $c=0$, then $(Y,e_3)$ is an isothermic pair.

\ms 

Next we study the class of solutions of the $G_{3,1}^1$-system II corresponding to
isothermic immersions of minimal surface in $R^3$.  
If $Y$ is an isothermic immersion of a minimal surface in $R^3$ such that 
$$I=e^{2u}(dx^2+ dy^2), \quad \II= dx^2-dy^2,$$ then 
$(Y,e_3)$ is an isothermic pair, where $e_3$ is the unit normal of $Y$.  It follows from
Propositions \refdp{}, \refkd{} and Corollary \refee{} that $((e_3-Y)/2,
(e_3+Y)/2)$ is a 2-tuple in
$R^3$ of type $O(1,1)$, and the corresponding solution
$(F,G,B)$ of the $G_{3,1}^1$-system II \refaad{}  is given by
\refeq[cj]$$F=\pmatrix{0&u_x\cr u_y&0\cr}, \quad G=(e^{-u},
-e^{-u}), \quad B=\pmatrix{\cosh u&\sinh u\cr \sinh u& \cosh
u\cr}.$$ 
The converse is a consequence of the definition of $G_{3,1}^1$-system
II, Corollary \refjj{}, Proposition \refed{} and Corollary \refdu{}:

\refclaim[ew] Proposition. If $u$ is a solution of \refdu{}, then 
$$v=\pmatrix{0& u_x\cr u_y& 0\cr e^{-u}& -e^{-u}\cr}$$ is a solution of the
$G_{3,1}^1$-system \refaaa{}.  Moreover, let $E$ be a frame  for $v$, i.e.,
$$E^{-1}dE = \o_\l =\pmatrix{\d v^t - v\d^t& -\d J \l\cr
\d^t \l & - Jv^t\d J + \d^t v\cr},$$ where $\d=\pmatrix{dx_1& 0\cr 0
& dx_2\cr 0&0\cr}$ and $J=\diag(1, -1)$.
Then
\item {(i)}$({\p E\over \p \l} E^{-1})(x,y,0)= \pmatrix{0& X(x,y)\cr -JX^t(x,y)&
0\cr}$ for some
$\cm_{3\times 2}$-valued map $X$,
\item {(ii)} $X=(X_1, X_2)$ is a 2-tuple in $R^3$ of type $O(1,1)$,
\item {(iii)} let $Y_1=X_1-X_2$ and $Y_2= X_1+X_2$, then $(Y_1, Y_2)$ is an
isothermic pair such that $Y_1$ is minimal and $Y_2$ is the unit normal of $Y_1$,
\item {(iv)} $(x,y)$ is the canonical isothermic coordinate system for minimal surface
$Y_1$.

\ms

It follows from Corollary \refke{} that a solution of \refdu{} gives a minimal
surfaces in $R^3$ and a surfaces in
$N^3(-1)=H^3$ with mean curvature $2$ unique up to ambient isometries. 
Proposition \refew{} gives a construction of the minimal immersion  in  $R^3$
from the frame of the $G_{3,1}^1$-system.
Next we show that the immersion of the corresponding CMC surface in $H^3$ can
also be read easily from the frame of the $G_{3,1}^1$-system II. First we recall that if
$E^0(x,y,\l)$ is a frame of the solution $v=\pmatrix{F\cr G\cr}$ of the
$G_{3,1}^1$-system, then $E^0(x,y,0)=\pmatrix{A(x,y)&0\cr 0&B(x,y)\cr}$ for
some $A, B$ and 
$$E(x,y,\l)= E^0(x,y,\l)\pmatrix{I&0\cr 0&B^{-1}(x,y)\cr}$$ is a frame
for the solution $(F,G,B)$ of the $G_{3,1}^1$-system II.

\refclaim[ex] Proposition.  Let $u$ be a solution of \refdu{}, and $(F,G,B)$ defined
by \refcj{}.  Then $(F,G,B)$ is a solution of the $G_{3,1}^1$-system II \refaad{}.
Moreover, let $E(x,y,\l)$ be a trivialization of the Lax pair $\o_\l^{\II}$ 
defined by
\refaaf{} for $(F,G,B)$, i.e., $E^{-1} dE= \o_\l^{\II}$, and $e_i(x,y)$ the i-th
column of $E(x,y,\sqrt{2})$.  Then:
\item {(i)} $v_0=\sqrt{2}\ e_3+e_4+e_5$ is independent of $(x,y)$ and $\li v_0,
v_0\ri_1=2$.
\item {(ii)} $Y=e_3 + \sqrt{2}\  e_5$ lies in the totally umbilic
hypersurface
$$N_{v_0,0}= \{p\in R^{4,1}\n \li p, p\ri_1= -1, \li p, v_0\ri_1=0\}$$ 
of $H^4$ and $N_{v_0,0}$ is isometric to $H^3$, where $H^n$ is the hyperbolic
space form $N^n(-1)$.
\item {(iii)} $Y$ has constant mean curvature $2$ in $H^3$.

\proof Let $\w_{ij}$ denote the ij-th entry of $\o_\l^{\II}$ at $\l=\sqrt{2}$. 
Since $dE= E\o^{\II}_\l$, we have
\refeq[cii]$$de_j= \sum_{i=1}^5 \w_{ij} e_i, \quad \w_{ij} +
\e_i\e_j \w_{ji}=0$$ for $1\leq j\leq 5$, where $\e_1=\cdots
=\e_4=-\e_5=1$.
Read $\w_{ij}$ from \refaaf{} to get 
\refeq[ch]$$\eqalign{&\w_{12} = u_y\ dx - u_x\ dy,\cr &\w_{13} =
e^{-u}dx, \quad \w_{14}= -\sqrt{2}\  \cosh u\  dx, \quad \w_{15} = \sqrt{2}\ \sinh u
\ dx,\cr &\w_{23}= -e^{-u}\ dy, \quad \w_{24}= - \sqrt{2}\ \sinh u\  dy, \quad
\w_{25} =  \sqrt{2}\ \cosh u \ dy, \cr & \w_{ij}=0, \quad {\rm if\,\,} 3\leq
i, j\leq 5.\cr}$$  

To prove (i), we  compute directly
$$\eqalign{dv_0&=\sqrt{2}\ de_3+ de_4+ de_5,\cr & =
\sqrt{2}\ \w_{13}e_1 +  \sqrt{2}\ \w_{23}e_2 + \w_{14}e_1 + \w_{24}e_2 +
\w_{15}e_1 + \w_{25}e_2 \cr &=(\sqrt{2} \ \w_{13}+\w_{14} +  \w_{15})e_1
+ (\sqrt{2} \ \w_{23} + \w_{24} + \w_{25})e_2.\cr }$$ Substitute \refch{} to
the above equation to conclude that $dv_0=0$. 

By Proposition \refdl{}, $N_{v_0,0}$ has constant sectional curvature
$-1$.  A direct computation gives (ii).

Use \refcii{}, \refch{} and a direct computation to get
$$dY = {e^u\over\sqrt{2}}\ (dx\ e_1 + dy\
e_2).$$ So $e_1, e_2$ are tangent to $Y$ and the dual 1-frame is 
$$\w_1= {1\over \sqrt{2}}\  e^u dx, \quad \w_2= {1\over \sqrt{2}}\ e^u dy.$$ 
The normal plane of
$N_{v_0,0}$ in
$R^{4,1}$ at $Y$ is spanned by $Y$ and
$v_0=\sqrt{2}\ e_3+e_4+e_5$.  It is easy to write down an $O(4,1)$-frame
$(e_1, e_2,
\tilde e_3, \tilde e_4, \tilde e_5)$ on $Y$ such that $\tilde e_3$ is normal to $Y$
in $N_{v_0,0}$ and $\tilde e_5=Y$: 
$$\eqalign{\tilde e_3&= e_3 - {1\over \sqrt{2}}\  e_4 + {1\over \sqrt{2}}\  e_5,
\cr
\tilde e_4&= (e_3+ {1\over \sqrt{2}}\ e_4 + {1\over \sqrt{2}}\ e_5),\cr
\tilde e_5&= Y = e_3+ \sqrt{2}\ e_5.\cr}$$
Let $\tilde e_i= e_i$ for $1\leq i\leq 2$, and $\tilde \w_{ij}= \li d\tilde e_i, \tilde
e_j\ri_1$. A direct computation gives
$$\tilde \w_{13}= (e^{-u} + e^u) dx, \quad \tilde \w_{23} = (-e^{-u} + e^u)
dy.$$ 
So  $Y$ has constant mean curvature $2$ in $H^3$. 
\qed

\bs


\newsection Loop group action for $G_{m,n}$-systems.\par

In this section,  we  give an explicit construction of the dressing action of a rational map
with two simple poles on the space of solutions of the $G_{m,n}$-systems.

First we review the construction of the dressing action.  Let $U/K$ be a symmetric
space, and $$\eqalign{ G_+&= \{ g: \Cx \to U_{\Cx}   \n g \, \,{\rm is\,
holomorphic, \, satisfies \, the\,} \cr & \quad U/K{\rm -reality\,
condition\, \refgg{}} \} ,\cr G_- &= \{ g:S^2\to U_{\Cx}\n
  g(\infty )=I , \, g {\rm\, is \, meromorphic,\, and}\cr
& \quad{\rm  satisfies \, the\, } U/K{\rm -reality\, condition \,
\refgg{}} \} .\cr}$$

Let $g\in G_-$, $v$ a solution of the $U/K$-system \refat{}, and
$E(x,\l)$ the frame of $v$ with $E(0,\l)=I$, i.e., $$E^{-1}dE=
\o_\l, \quad E(0,\l)=I,$$ where $\o_\l$ is defined by \refau{}.  Since $\o_\l$ is
holomorphic for
$\l\in
\Cx$, $E(x,\cdot)\in G_+$. It is known that (cf. [TU2])
$g(\l)E(x,\l)$ can be factored as 
\refeq[fa]$$g(\l) E(x,\l) =
\tilde E(x,\l) \tilde g(x,\l)$$ such that
 $\tilde E(x,\cdot)\in G_+$ and $\tilde
g(x,\l)\in G_-$ for  $x$ in a neighborhood of the origin.
Moreover, this factorization \reffa{} can be obtained explicitly
using residue calculus, and we have:
 \item {(i)} $\tilde E$ is the frame for
a new solution $\tilde v$ of the $U/K$-system with $\tilde
E(0,\l)=I$.
\item {(ii)} Write $\tilde g(x,\l) = I + \l^{-1}m_1(x) + \l^{-2}m_2(x) +
\cdots$. Then  $\tilde v=v- p_1(m_1)$ is a new solution of the
$U/K$-system, where $p_1$ is the projection onto
$\cp\cap\ca^\perp$.
\item {(iii)}  $g\,\sharp\, v= \tilde v$ defines an action of $G_-$ on the
space of germs of solutions of the $U/K$-system.
\item {(iv)} If $U$ is compact and $v$ is a smooth solution decaying at
infinity, then the factorization \reffa{} can be carried out for
all $x\in R^n$ and the solution $g\sharp v$ is globally defined.
\item {(v)} Let $\tilde \o$ denote the Lax connection of the solution $\tilde v$.
Then \refeq[gj]$$\tilde g \o - d\tilde g = \tilde \o \tilde g,$$
and this gives a system of compatible ordinary differential equations
for $\tilde v$.

\ms

\refpar[gk] Remark. Let $E$ be a frame (cf. Definition \refgi{})
of the solution $v$ of the $U/K$-system \refat{}. Then
\item {(i)} $g(\l)E(x,\l)$ is a frame of $v$ if and only if
$g$ satisfies the $U/K$-reality condition,
\item {(ii)}   if $g_0\in U_\Cx$ is a constant, then $g_0E(x,\l)$ is a
frame of $v$ if and only if $g_0\in K$.

\ms Next we find certain simple element in
$G_-$ explicitly. The $G_{m,n}$- reality condition is
\refeq[fb]$$  \cases{ \overline{g(\bar{\lambda})}=g({\lambda
}),\cr  I_{m,n} g({\lambda }) I_{m,n}^{-1} = g({-\lambda }),\cr
g({\lambda })g({\lambda})^{t} = I .}$$

\refpar[ii] Remark. If $g(\l)$ satisfies the $G_{m,n}$- reality condition
\reffb{}, then
\item {(i)} $g(0)\in O(m)\times O(n)$ and $g(\l)\in O(m+n,C)$ for all
$\l\in C$,
\item {(ii)} $g(\l)$ satisfies the $U(n)$-reality condition, i.e.,
$$g(\bar\l)^* g(\l)=I.$$

It is known (cf. [U1]) that the group of rational maps $g:S^2\to
GL(n,C)$ satisfying the $U(n)$-reality condition
$\overline{g(\bar\l)}^tg(\l)= I$ is generated by the set of 
\refeq[fc]$$h_{z,\pi}(\l)=
\pi + {\l-z\over \l-\bar z} (I-\pi),$$ where $z\in C$ and $\pi$ is a Hermitian projection of
$C^n$.   Although $h_{z,\pi}$ defined by \reffc{} does not satisfy
the $G_{m,n}$- reality condition \reffb{}, a product of suitable choices of
two such elements satisfies \reffb{}.  To construct this element, let $R^m=\cm_{m\times
1}$,
$W$ and $Z$ unit vectors in $\reals^m$  and $\reals^n$
respectively,  and $\pi $   the Hermitian projection of $C^{n+m}$
onto $ \Cx\pmatrix{ W \cr iZ \cr} $, the complex linear subspace
spanned by $\pmatrix{W\cr iZ\cr}$.   So \refeq[za]$$ \pi= {1\over
2} \pmatrix{ W W^t & -i WZ^t  \cr
 iZW^t & ZZ^t \cr}.$$
Note that $\bar\pi$ is the Hermitian projection onto $\Cx
\pmatrix{-W\cr iZ\cr}$, which is perpendicular to $\pmatrix{W\cr
iZ}$. This implies $$\pi\bar{\pi}=\bar{\pi}\pi = 0.$$ Let $s\in
\reals, s\not= 0 $. Define \refeq[zb]$$g_{s,\pi}(\lambda )  =
\left( \pi +  {\lambda -is\over \lambda +is} (I - \pi ) \right)
\left( \bar{\pi} + {\lambda +is\over \lambda -is}  (I - \bar{\pi}
) \right).$$ Substitute \refza{} to \refzb{} to get
\refeq[hg]$$g_{s,\pi}(\l)= I + {2s\over \l^2 +
s^2}\pmatrix{-sWW^t& \l WZ^t\cr -\l ZW^t& -s ZZ^t\cr}.$$ A direct
computation implies that $g_{s,\pi}$ satisfies the $G_{m,n}$-reality
condition \reffb{}.  So $g_{s,\pi}\in G_-$.  Note also that
$$g_{s,\pi}(0)=\pmatrix{ I -2   WW^t & 0  \cr 0 & I
 -2 ZZ^t \cr}.$$

Below we give an explicit construction of the dressing action of $g_{s,\pi}$ on the
space of solutions of the $G_{m,n}$-system \refub{}. 

\refclaim[zc]Theorem.  Let $\xi: R^n\to \cm_{m\times n}$ be a
solution of the $G_{m,n}$-system \refub{}, and
$E(x,\l)$ a frame of $\xi$ such that $E(x,\l)$ is holomorphic for
$\l\in C$. Let $W$ and $Z$ be unit vectors in $R^m, R^n$
respectively, $\pi$ the Hermitian projection of $\Cx^{n+m}$ onto
$\Cx \pmatrix{W\cr iZ\cr}$, and $g_{s, \pi}$ the map defined by
\refzb{}.  Let $\tilde{\pi}(x)$  denote the Hermitian projection onto
 $\Cx\pmatrix{\tilde{W} \cr i\tilde{Z}\cr}(x)$, where
\refeq[gl]$$\pmatrix{\tilde{W} \cr
i\tilde{Z}\cr}(x)=  E(x,{-is}) ^{-1} \pmatrix{W\cr iZ \cr}.$$
Let $\hat W=\tilde W/\N \tilde W\N$ and  $\hat Z=\tilde
Z/\N \tilde Z\N$,  
\refeq[fd]$$\tilde E(x,\l)= g_{s,\pi}(\l)
E(x,\l) g_{s, \tilde \pi(x)}(\l)^{-1},$$ 
\refeq[fj]$$\tilde \xi =
\xi - 2s (\hat W\hat Z^t)_\ast,$$ 
where $(\eta_\ast)_{ij}= \eta_{ij}$ if $i\not=j$, and
$(\eta_\ast)_{ii}=0$ if $1\leq i\leq n$.
 Then $\tilde \xi$ is a solution of \refub{}, $\tilde E$ is a
frame for $\tilde \xi$ and $\tilde E(x,\l)$ is holomorphic in $\l\in C$.

When $m=n$, Theorem \refzc{} was proved in [Zh].  

\ms
To prove Theorem \refzc{}, we need the following lemma.

\refclaim[fo] Lemma.  With the same assumption as in Theorem
\refzc{}, the following statements are true:
\item {(i)} $\tilde W(x)\in R^m$ and $\tilde Z(x)\in R^n$.
\item {(ii)}  $\N \tilde
W(x)\N = \N \tilde Z(x)\N$ for all $x$, and $g_{s,\tilde \pi(x)}$
satisfies the $G_{m,n}$-reality condition \reffb{}, i.e., it belongs to
$G_-$.
\item {(iii)} $\tilde E(x,\l)$ is holomorphic in  $\lambda \in \Cx$.

\proof

(i) Let $A^\ast$ denote $\bar A^t$.     Since $E(x,\l)$ satisfies
the $G_{m,n}$- reality condition \reffb{}, \refeq[gb]$$I_{m,n}E(x,{-is})^{-1}
I_{m,n}^{-1} = E(x, is)^{-1} = E(x, -is)^\ast.$$
 Write $E(x, -is)^{-1}= \pmatrix{\eta_1& \eta_2\cr
\eta_3&\eta_4\cr}$ with $\eta_1\in gl(m,C)$ and $\eta_4\in
gl(n,C)$. It follows from \refgb{} that $\eta_1, \eta_4$ are real
and $\eta_2, \eta_3$ are pure imaginary.  This implies that
$\tilde W, \tilde Z$ are real. \ss

(ii) Let $\li v_1, v_2\ri = v_1^t v_2$ denote the standard
bi-linear form on $\Cx^{n+m}$. Then  $\li u, u\ri=0$, where
$u=\pmatrix{W\cr iZ\cr}$. Since $E$ satisfies the reality
condition \reffb{}, $E(x,\l)\in O(m+n,C)$.  So $$\li
E(x,-is)^{-1}(u), E(x,-is)^{-1}(u)\ri = \li u, u\ri =0.$$ This
implies that $\N \tilde W(x)\N = \N \tilde Z(x)\N$ and $g_{s,
\tilde \pi(x)}$ satisfies the reality condition \reffb{}.

(iii) If $g$ satisfies \reffb{}, then $g(\l)^{-1}= g(\bar\l)^*$.
So \refeq[fy]$$\eqalign{&\tilde E(x,\l) = g_{s,\pi}(\l) E(x,\l)
g_{x, \tilde \pi(x)}(\bar \l)^\ast =\cr & (\pi+{\l-is\over
\l+is}\pi^\perp)(\bar\pi + {\l+is\over
\l-is}\bar\pi^\perp)E(x,\l)(\overline{\tilde \pi} + {\l-is\over
\l+is}\overline{\tilde \pi}^\perp)( \tilde \pi+{\l+is\over
\l-is}\tilde \pi^\perp),\cr}$$ where $\pi^\perp= I-\pi$ and
$\tilde \pi^\perp= I-\tilde \pi$. The right hand side only has
poles at $is, -is$ of order 2.
First note that the coefficient of ${1\over (\l-is)^2}$ on the
right hand side of \reffy{} is $-4s^2\pi\bar\pi^\perp
E(x,is)\overline{\tilde \pi}\tilde\pi^\perp$, which is equal to $
-4s^2\pi E(x,is)\overline{\tilde \pi}$ because
\refeq[fz]$$\pi\bar\pi = \bar\pi \pi =0, \quad \tilde \pi
\overline{\tilde \pi}= \overline{\tilde \pi}\pi=0.$$  Claim that
\refeq[gd]$$\pi E(x,is)\overline{\tilde \pi} =0.$$ To see this, we
note that the image of $\pi E(x,is)\overline{\tilde \pi}$ is
spanned by $$\eqalign{&\pi E(x, is) \pmatrix{\tilde W\cr -i \tilde
Z\cr} = \pi E(x,is) I_{m,n}\pmatrix{\tilde W\cr i\tilde Z\cr},
\quad {\rm by\,} \refgb{} \cr &= \pi
I_{m,n}E(x,-is)\pmatrix{\tilde W\cr i\tilde Z\cr} = \pi I_{m,n}
\pmatrix{W\cr iZ\cr} = \pi \pmatrix{W\cr -iZ\cr} = 0.\cr}$$  This
proves the claim.  So the coefficient of ${1\over (\l-is)^2}$ on
the right hand side of \reffy{} is zero.

The coefficient of ${1\over \l-is}$ on the right hand side of
\reffy{} is $4is\pi E(x,is)\overline{\tilde \pi}$,
 which is zero because  of
\refgd{}.  So $\tilde E(x,\l)$ is holomorphic at $\l=is$. Similar
computation implies that $\tilde E(x,\l)$ is holomorphic at $\l=
-is$. Hence $E(x,\l)$ is holomorphic for $\l\in \Cx$. \qed \ms

\refpar[] Proof of Theorem \refzc{}.

A direct computation gives \refeq[fg]$$\tilde E^{-1} d\tilde E=
\tilde g E^{-1} dE \tilde g^{-1} - d\tilde g\tilde g^{-1},$$ where
$\tilde g(x, \l)=g_{s, \tilde\pi(x)}(\l)$. But $E^{-1} dE= \sum
(a_i\l + [a_i,v])dx_i$, $\tilde E(x,\l)$ is holomorphic in $\l\in
\Cx$ and $\tilde g(x,\l)$ is holomorphic at $\l=\infty$. So
$\tilde E^{-1} d\tilde E$ must be of the form $\sum (a_i \l +
\eta_i) dx_i$. Write $$\tilde g(x,\l)= I + \l^{-1} m_1(x) +
\cdots.$$ A direct computation shows that \refeq[ge]$$m_1(x)=
2s\pmatrix{0& \hat W\hat Z^t\cr -\hat Z\hat W^t& 0\cr}, \quad
m_1(x)\in \cp.$$  Multiply \reffg{} by $\tilde g$ on the right to
get $$\eqalign{&\left(\sum_i (a_i\l + \eta_i)dx_i
\right)(1+m_1\l^{-1} + \cdots) \cr &\,\, =(I+m_1\l^{-1} +\cdots)
\left(\sum_i (a_i\l + [a_i, v])dx_i\right) -(dm_1 \l^{-1}
+\cdots).\cr}$$ Equate the constant term of the above equation to
get $$\eta_i = [a_i, v-m_1] = [a_i, v-p_1(m_1)],$$ where $p_1$ is
the projection from $\cp$ onto $\cp\cap\ca^\perp$. Write $\tilde
v=\pmatrix{0&\tilde \xi\cr -\tilde \xi^t & 0\cr}$.  Theorem
follows from \refge{}. \qed
 \ms

Since $E(x,\l)$ in Theorem \refzc{} is not assumed to satisfy the initial condition
$E(0,\l)=I$, the resulting new solution $\tilde \xi$ is not necessarily equal to 
the dressing action of $g_{s,\pi}$ on $\xi$.  But they are related as follows:

\refclaim[fx] Corollary.  Suppose $E$ is a frame of the solution
$\xi$ of the $G_{m,n}$-system \refub{} such that $E(x,\l)$ is holomorphic for $\l\in
\Cx$.
\item {(i)} If $E(0,\l)=I$, then  $\tilde
\xi$ obtained in Theorem \refzc{} is  $g_{s,\pi}\sharp \xi$ and
$\tilde E$ is the frame of $\tilde \xi$ with $\tilde E(0,\l)=I$.
\item {(ii)} Let $g_+(\l)= E(0,\l)$, and $\tilde \xi$ the new solution of \refub{} 
obtained in Theorem \refzc{}. Then $g_+\in G_+$ and  $\tilde
\xi=\tilde g_-\sharp \xi$, where $\tilde g_-$ is obtained by
factoring $g_{s,\pi}g_+= \tilde g_+ \tilde g_-$ with $\tilde
g_\pm\in G_\pm$.

The functions $\tilde W, \tilde Z$ in Theorem \refzc{} also
satisfy a system of compatible first order differential equations.
This follows from taking the differential of the defining equation
\refgl{} of $\tilde W, \tilde Z$: \refeq[fi]$$\eqalign{ &
d\pmatrix{\tilde{W} \cr i\tilde{Z} \cr} =d E_{-is}^{-1} \pmatrix{W
\cr iZ \cr}  = - E_{-is}^{-1} dE_{-is} \pmatrix{\tilde{W}\cr
i\tilde{Z} \cr}  =- \o_{-is} \pmatrix{ \tilde{W} \cr i\tilde{Z}
\cr}, \cr }$$ where $\o_\l$ is defined by \refav{} and
$E_\l(x)=E(x,\l)$.  We write system \reffi{} explicitly:
\refeq[fh]$$\cases{\tilde W_{x_i}=(\xi D_i^t-D_i\xi^t)\tilde W +
sD_i\tilde Z, &\cr \tilde Z_{x_i}= sD_i^t\tilde W + (\xi^tD_i -
D_i^t\xi)\tilde Z,&\cr}$$ or 
\refeq[gs]$$  \cases{(\tilde{w}_i
)_{x_i } =-\sum_{j\not= i }^n f_{ji}\tilde{w}_j -\sum_{j=1}^{m-n}
g_{ji}\tilde{w}_{j+n} +s\tilde{z}_i   & $i\leq n$, \cr
(\tilde{w}_i )_{x_j} = f_{ij}  \tilde{w}_j & $i\leq n , j\not= i$,
\cr (\tilde{w}_i )_{x_j} =   g_{ij}\tilde{w}_j   & $i>n$ ,\cr
(\tilde{z}_i )_{x_j} =   f_{ji} \tilde{z}_j & $j\not= i$ ,\cr
(\tilde{z}_i )_{x_i}=-\sum_{j\not= i}^n  f_{ij} \tilde{z}_j +
s\tilde{w}_i .&}$$ 
The converse is also true:

\refclaim[zf] Proposition. Given $\xi:R^n\to \cm_{m\times n}^0$ and
 a real number $s\not=0$, we have:
\item {(i)} System \reffh{} is solvable for $\ti W, \ti Z$ if and only if $\xi$ is a solution of
the $G_{m,n}$-system \refub{}.
\item {(ii)} Let $\xi$ be a solution of \refub{}, and $\ti W, \ti Z$ solution of
\reffh{} with initial condition $\ti W(0)= W, \ti Z(0)= Z$ for some unit vectors
$W\in R^m, Z\in R^n$. Set
$$\tilde \xi = \xi - 2s {(\tilde W\tilde Z^t)_\ast\over \N \tilde
W\N \N \tilde Z\N},$$  
where $(Y_\ast)_{ij}=
Y_{ij}$ if $i\not=j$, and $(Y_\ast)_{ii}=0$ if $1\leq i\leq n$.
Then  $\ti \xi$ 
is a solution of \refub{}, and is equal to $g_{s,\pi}\sharp \xi$, where $\pi$ is the
projecton onto $C(W, iZ)^t$.

\proof System \reffh{} is the same as system \reffi{}, i.e.,
\refeq[fp]$$d\pmatrix{\tilde W\cr i\tilde Z\cr}=
-\o_{-is}\pmatrix{\tilde W\cr i\tilde Z\cr},$$ where
$$\o_{-is}=\sum_{j=1}^n \pmatrix{D_j\xi^t-\xi D_j^t& isD_j\cr
-isD_j& -\xi^tD_j + D_j^t\xi\cr} \, dx_j.$$ But \reffp{} is
solvable if and only if $\o_{-is}$ is flat. This is equivalent to
$\xi$ being a solution of \refub{}.  So there exists a unique
solution $\tilde W, \tilde Z$ of \reffh{} such that $\tilde W(0)=
W$ and $\tilde Z(0)=Z$.   Let $E(x,\l)$ be the frame of $\xi$ with
$E(0,\l)=I$. Then $E_{-is}^{-1}\pmatrix{W\cr i Z\cr}$ is also a
solution of \reffh{} with the same initial condition.  Hence they
must be equal and Proposition follows from Theorem \refzc{}. \qed

\ms

Since the Lax connection of the $G_{m,n}$-system I, II are gauge equivalent
to the Lax connection of the $G_{m,n}$-system, we can write down the action
of $g_{s,\pi}$ on solutions of these two systems easily from the
action on solutions of the $G_{m,n}$-system.

Use the same notations as in Theorem \refzc{}.   Let
\refeq[go]$$\eqalign{\tilde E^\na(x,\l)&=g_{s,\pi}^{-1}\tilde
E(x,\l)= E(x,\l) g_{s,\tilde\pi(x)}^{-1}\cr &=
E(x,\l)\pmatrix{I_m- {2s^2\over \l^2+s^2}\hat W\hat W^t&
-{2s\l\over \l^2+s^2} \hat W\hat Z^t\cr {2s\l\over \l^2+s^2} \hat
Z\hat W^t & I_n- {2s^2\over \l^2+s^2} \hat Z\hat Z^t\cr}.\cr}$$
Since both $g_{s,\pi}$ and $\tilde E$ satisfy the $G_{m,n}$- reality
condition \reffb{}, so does $\tilde E^\na$.  Hence $\tilde E^\na$ is
a frame of $\tilde \xi$.  Note that $\tilde E^\na(x,\cdot)$ is not
in $G_+$.   The $G_{m,n}$- reality condition implies that both $E(x,0)$ and
$\tilde E^\na(x,0)$ are in  $O(m)\times O(n)$.  Write
$$E(x,0)=\pmatrix{A(x)&0\cr 0&B(x)\cr}, \quad \tilde
E^\na(x,0)=\pmatrix{\tilde A^\na(x)&0\cr 0& \tilde B^\na(x)\cr}.$$
It follows from \refgo{} that we
have \refeq[ga]$$\cases{\tilde A^\na=A(I-2\hat W\hat W^t),&\cr
\tilde B^\na= B(I-2\hat Z\hat Z^t).&\cr}$$ Write
$$\xi=\pmatrix{F\cr G\cr}, \quad \tilde \xi=\pmatrix{\tilde F\cr
\tilde G\cr},\quad A=(A_1,A_2),\quad \tilde A^\na=(\tilde A^\na_1,
\tilde A^\na_2),$$ where $A_1, \tilde A^\na_1\in \cm_{m\times n}$
and $A_2, \tilde A^\na_2\in \cm_{m\times(m-n)}$.  Rewrite \reffj{}
as $$\pmatrix{\tilde F\cr \tilde G\cr}= \pmatrix{F\cr G\cr} -
2s\left(\hat W\hat Z^t\right)_\ast.$$
So $(A_1, F)$ and
$(\tilde A^\na_1, \tilde F)$ are solutions of the
$G_{m,n}$-system I \refcn{}, and $(F,G,B)$ and
$(\tilde F, \tilde G, \tilde B^\na)$ are solutions of the
$G_{m,n}$-system II \refcu{}.  Recall that
\item {(i)} $\o_\l^I$ defined by \refbo{} is the Lax connection of the solution $(A_1,
F)$ of the $G_{m,n}$-system I \refcn{},
\item {(ii)}  $\o_\l^{\II}$ defined by \refbp{} is the Lax connection of solution
$(F,  G,  B)$ of the $G_{m,n}$-system II \refcu{},
\item {(iii)}  $\o_\l^I$ is  the gauge transformation of
$\o_\l$ by $\pmatrix{A&0\cr 0&I_n\cr}$, and $\o_\l^{\II}$ is the
gauge transformation of $\o_\l$ by $\pmatrix{I_m&0\cr 0&B\cr}$.

\ni So the frames  $E$ of $\xi$, $\tilde E^\na$ of $\tilde \xi$,
$E^I$ of $(A_1, F)$, $\tilde {E^\na}^I$ of $(\tilde A^\na_1,
\tilde F)$ and $E^{\II}$ of $(F,G,B)$ and $\tilde {E^\na}^{\II}$ of
$(\tilde F, \tilde G, \tilde B^\na)$ are
 related by
$$\eqalign{E^I(x,\l)&=E(x,\l)\pmatrix{A^t&0\cr 0& I_n\cr }, \cr
\tilde {E^\na}^I(x,\l) &= \tilde E^\na(x,\l)\pmatrix{\tilde
{A^\na}^t&0\cr 0& I_n\cr},\cr \tilde {E^\na}^{\II}(x,\l)&= \tilde
{E^\na}(x,\l)\pmatrix{I_m&0\cr 0& \tilde{B^\na}^t\cr}.\cr}$$ Use
\refgo{} to get \refeq[fq]$$\eqalign{  \tilde
{E^\na}^I(x,\l)&=E^I(x,\l)\left(I-{2\over
\l^2+s^2}\pmatrix{\l^2A\hat W\hat W^t A^t& s\l A\hat W\hat Z^t\cr
s\l\hat Z\hat W^t  A^t& s^2 \hat Z\hat Z^t\cr}\right),\cr
&=E^I(x,\l)\left(I-{2\over \l^2+s^2}\pmatrix{\l A\hat W\cr s\hat
Z\cr}\pmatrix{\l\hat W^tA^t & s\hat Z^t}\right), \cr}$$ and
\refeq[gr]$$\tilde {E^\na}^{\II}(x,\l)= E^{\II}(x,\l) \left(I -
{2\over \l^2+s^2}\pmatrix{s^2\hat W\hat W^t & -s\l \hat W\hat Z^t
B^t\cr -s\l B\hat Z\hat W^t & \l^2 B\hat Z\hat Z^t
B^t\cr}\right).$$ 

It follows from
Corollary \reffx{} that the new solution $\tilde \xi$ obtained in
Theorem \refzc{} depends on $g_{s,\pi}$ and the frame $E$.
Henceforth, we will use the following notations:
$$\eqalign{&(\tilde \xi, \tilde E^\na)=g_{s,\pi}\cdot (\xi, E),
\cr & \tilde {A^\na} = g_{s,\pi}\cdot A, \quad \tilde {B^\na} =
g_{s,\pi}\cdot B,\cr &  (\tilde A^\na_1, \tilde F, \tilde
{E^\na}^I)=g_{s,\pi}\cdot (A_1, F, E^I),\cr & (\tilde F, \tilde G,
\tilde B^\na, \tilde {E^\na}^{\II})=g_{s,\pi}\cdot (F,G, B,
E^{\II}).\cr}$$

\ms

Use exactly the same argument as for Theorem \refzc{} to get the
action of $g_{s,\pi}$ on solutions of the partial
$G_{m,n+1}$-system \reffr{}.   We summarize the
results for this case below.

\refclaim[zk]Theorem.  Let $(F , G  , b )$  be a solution of the partial
$G_{m,n+1}$-system 
\reffr{}, $\Theta_{\lambda}$ (defined by \refaw{}) its Lax connection,
and  $E(x,\l)$ a frame of $(F,G,b)$. Let $W\in R^m$ and $Z\in
R^{n+1}$ be unit vectors, $\pi$ the Hermitian projection onto
$\Cx\pmatrix{W\cr iZ\cr}$, and $g_{s,\pi}$ defined by \refzb{}.
Let \refeq[fu]$$\pmatrix{\tilde{W} \cr i\tilde{Z}\cr}(x)= E(x,-is)
^{-1}\pmatrix{ W\cr iZ \cr},$$ $\tilde{\pi}(x)$  the Hermitian
projection onto
 $\Cx\pmatrix{ \tilde{W} \cr i\tilde{Z} \cr} (x)$,
$$\tilde{E} (x,\l)=
g_{s,\pi}(\l)E(x,\l)g_{s,{\tilde\pi(x)}}^{-1}(\l),$$
\refeq[fs]$$\pmatrix{\tilde F&\tilde b\cr \tilde G&0\cr} =
\pmatrix{F& b\cr G &0\cr} - 2s(\hat W\hat Z^t)_\ast,$$ where $\hat
W=\tilde W/\N \tilde W\N$, $\hat Z=\tilde Z/\N \tilde Z\N$.  Then:
\item {(i)} $(\tilde F, \tilde G,
\tilde b)$ is a solution of the partial $G_{m,n+1}$-system \reffr{} and $\tilde E$ is a
frame.
\item {(ii)}  $(\tilde W, \tilde Z)$ is a solution of
\refeq[ft]$$\cases{(\tilde{W})_{x_j } =-\pmatrix{-FC_j+ C_jF^t&
C_j G^t\cr -GC_j& 0\cr}\tilde W +sC_j \tilde Z, &\cr (\tilde
Z)_{x_j} = -\pmatrix{-F^tC_j+C_jF & C_j b\cr -b^t C_j& 0\cr}\tilde
Z + s (C_j, 0)\tilde W.&\cr}$$
\item {(iii)} System \refft{} is
solvable for $(\tilde W, \tilde Z)$ if and only if $(F,G,b)$ is a
solution of \reffr{}.
\item {(iv)} If $(F,G,b)$ is a solution of \reffr{} and $(\tilde W,\tilde
Z)$ is a solution of \refft{} such that $\tilde W(0), \tilde Z(0)$
are unit real vectors, then $(\tilde F, \tilde G, \tilde b)$ defined
by formula \reffs{} is a solution of \reffr{}.

 Let $(F,G,b)$ be a solution of the partial $G_{m,n+1}$-system \reffr{}, and $E$ a frame
of
$(F,G,b)$.  Since $E(x,0)\in O(m)\times O(n+1)$, there exist
$A(x)\in O(m)$, and $B(x) \in O(n+1)$ such that
$$E(x,0)=\pmatrix{A(x)& 0\cr 0& B(x)\cr}.$$ Write $A=(A_1,A_2)$
with $A_1\in \cm_{m\times n}$ and $A_2\in \cm_{m\times(m-n)}$.
Then $(A_1, F,b)$ is a solution of the partial $G_{m,n+1}$-system I \refcl{} and
$$E^I(x,\l)= E(x,\l)\pmatrix{A^t&0\cr 0&I_n}$$ is a frame  of $(A_1,F,b)$. Let
$(\tilde F, \tilde G, \tilde b)$, $\hat W, \hat Z$, and $\tilde E$
be as in Theorem \refzk{}, and $$\tilde E^\na= g_{s,\pi}^{-1}\tilde E,
\quad \tilde A^\na=A(I-2\hat W\hat W^t).$$  Then
$(\tilde{A^\na}_1 ,\tilde{F},\tilde{b})$ is a new solution of the
partial $G_{m,n+1}$-system I, where $\tilde
A^\na=(\tilde A_1^\na, \tilde A_2^\na)$. Moreover, $$\eqalign
{\tilde {E^\na}^I(x,\lambda)= E^I(x,\lambda)
  \left(I-{2\over \l^2+s^2}\pmatrix {\l\ A \hat{W}\cr
s \hat{Z}} \pmatrix { \l  \hat{W}^t A^t & s \hat{Z}^t}\right)  }$$
is a frame of $(\tilde F, \tilde G, \tilde b)$.  We will use
the following notations:$$\eqalign{& (\tilde F, \tilde G, \tilde b, \tilde
E^\na)= g_{s,\pi}\cdot (F, G, b, E), \cr &( \tilde A_1^\na,\tilde F,\tilde b, \tilde
{E^\na}^I) = g_{s,\pi}\cdot (A_1, F, b, E^I).\cr}$$

\bs

\newsection Ribaucour Transformations for $G_{m,n}$-systems.\par

The main goal of this section is to give geometric interpretations
of the action of $g_{s,\pi}$ on the spaces of solutions of the various
$G_{m,n}$-systems constructed in section 9.

We first review some notions in classical differential geometry (cf. [Da], [Bi1], [Bi2]).
Given a two parameter family of spheres $S(x,y)$ in $R^3$, $$S(x,y):
\quad p(x,y) + r(x,y) w, \quad w\in S^2,$$ generically there exist
two enveloped surfaces $M, \tilde M$, i.e., they are tangent to
these spheres. To see this, we fix a parametrization $w(u,v)$ on
$S^2$. To find an enveloped surface is to find $u(x,y), v(x,y)$ so
that $w(u(x,y), v(x,y))$ is normal to the surface $$X(x,y)=
p(x,y)+ r(x,y)w(u(x,y), v(x,y))$$ at $X(x,y)$. So $u,v$ need to
satisfy $X_x\cdot w=0$ and $X_y\cdot w=0$. Or equivalently,
\refeq[gn]$$p_x\cdot w + r_x=0, \quad p_y \cdot w + r_y=0.$$ This
means that $w(u(x,y), v(x,y))$ must lie in the intersection of the
two circles defined by \refgn{} on $S^2$. Since generically two
such circles intersect at exactly two points for each $(x,y)$, we
obtain two surfaces $M, \tilde M$ and  a map $\ell:M\to \tilde M$
so that $X(x,y)$ and $\tilde X(x,y)$ lie in the same sphere
$S(x,y)$.  Note that the map $\ell$ is characterized by the
property that the normal line at $q$ to $M$ and the normal line at
$\ell(q)=\tilde q$ to $\tilde M$ meet at equidistance $r(q)$.   We
call such map $\ell$  a {\it sphere congruence\/}.   A sphere congruence $\ell:M\to \tilde
M$ is called a {\it Ribaucour transformation\/} if it preserves line of curvature directions
and the lines $p+te$ and  $\ell(p)+t \ell_\ast(e)$ meet at equidistance for any principal
direction $e\in TM_p$ and $p\in M$.  

\ms
 Natural generalizations of sphere congruence and Ribaucour transformation to
submanifolds in space forms are given in [DT].  For $x\in N^m(c)$ and $v\in
TN^m(c)_x$, let
$\g_{x,v}(t)= \exp(tv)$ denote the geodesic.  

\refpar[zn] Definition ([DT]).  Let $M^n$ and $\tilde M^n$ be
submanifolds of the space form $N^m(c)$.  A {\it sphere
congruence\/} is a vector bundle isomorphism $P:\nu(M)\to
\nu(\tilde M)$ that covers a diffeomorphism $\ell:M\to \tilde M$
with the following properties:
\item{(a)} If $\xi$ is a parallel normal vector field of $M$,
then $P\circ\xi\circ\ell^{-1}$ is a parallel normal field of
$\tilde M$.
\item{(b)} For any nonzero vector $\xi\in \nu_x(M)$, the geodesics $\g_{x,\xi}$ and
$\g_{\ell(x), P(\xi)}$ intersect at a point that is equidistant from $x$ and $\ell(x)$
(the distance depends on $x$).

\refpar[ht] Definition ([DT]). A sphere congruence $P:\nu(M)\to
\nu(\tilde M)$ that covers $\ell:M\to \tilde M$ is called a {\it
Ribaucour Transformation\/} if  it satisfies the following
additional properties:
\item {(i)}  If $e$ is an eigenvector of the shape
operator $A_{\xi}$ of $M$, then $\ell_{\ast} (e)$ is an
eigenvector of the shape operator $A_{P(\xi)}$ of $\tilde M$.
\item {(ii)} the geodesics $\g_{x,e}$ and $\g_{\ell(x), \ell_\ast(e)}$
intersect at a point equidistant to $x$ and $\ell(x)$. 

\ms

Below we show that the action of $g_{s,\pi}$ on solutions of the $G_{m,n}$-system I
correspond to a Ribaucour transformation for flat submanifolds in $S^{n+m-1}$.  

\refclaim[zo]Theorem.  Let
$E^I$ be a frame of the solution $(A_1, F)$ of the $G_{m,n}$-system I \refcn{},
$g_{s,\pi}$ given by \refzb{}, and $$(\tilde A_1^\na,\tilde F,
\tilde {E^\na}^I)= g_{s,\pi}\cdot (A_1, F, E^I)$$ as in section 9.  Write
\refeq[gq]$$\eqalign{E^I(x,1)&=(X(x), e_{n+2}(x), \cdots,
e_{n+m}(x), e_1(x), \cdots, e_n(x)),\cr \tilde {E^\na}^I(x,1) &=
(\tilde X(x), \tilde e_{n+2}(x), \cdots, \tilde e_{n+m}(x), \tilde
e_1(x), \cdots, \tilde e_n(x)).\cr}$$ Then:
\item {(1)} Both $X$ and $\tilde X$ are immersions of flat n-dimensional
submanifolds of $S^{n+m-1}$ with  flat, non-degenerate normal
bundle, $x_1, \cdots, x_n$  line of curvature coordinates,
$\{e_\alpha\}_{\a=n+2}^{n+m}$ and $\{\tilde e_\a\}_{\a=n+2}^{n+m}$
are parallel normal frames for $X$ and $\tilde X$ respectively.
\item {(2)}  $(A_1, F)$ and $(\tilde {A^\na_1}, \tilde F)$ are solutions of
\refcn{} corresponding to $X$ and $\tilde X$ as in Theorem
\refaa{} respectively.
\item {(3)} The bundle morphism $P:\nu(M)\to \nu(\tilde M)$ defined by
$P(e_\a(x))= \tilde e_\a(x)$ for $n+2\leq \a\neq n+m$  is a
Ribaucour Transformation covering the map $X(x)\mapsto \tilde
X(x)$.

\proof

(1) and (2) follow from   Theorems  \refaa{} and \refzc{}. \ss

(3) Let  $A_2, G$ and $A =(A_1 ,A_2 )$ be given as in Proposition
\reffm{}, and $\hat{W},\hat{Z}$ as in  Theorem \refzc{}. Let
$$\gamma = (\gamma_{n+1},\ldots ,\gamma_{m+n}, \gamma_{1} ,\ldots
,\gamma_{n})=\pmatrix {\sin\rho\, \hat{W}^t A^t \quad
 \cos\rho\, \hat{Z}^t},$$ where $\sin \rho = 1/\sqrt{1+s^2}$ and
$\cos \rho = s/\sqrt{1+s^2}$.   Substitute $\l=1$ in \reffq{} to get
\refeq[jm]$$\tilde {E^\na}^I(x,1) = E^I(x,1) \left( I-2\pmatrix{\sin\rho\,
A\hat W\cr \cos\rho\, \hat Z\cr}\pmatrix{\sin \rho\, \hat W^t A^t
& \cos\rho\, \hat Z^t\cr} \right).$$ 
Substitute \refgq{} to the above equation to get formulas for each column of $\tilde
{E^\sharp}^I$:
$$\eqalign{\tilde X&= X(1-2\g_{n+1}^2)- 2\g_{n+1}\left(\sum_{j=2}^m
\g_{n+j}e_{n+j} +\sum_{j=1}^n \g_j e_j\right),\cr
\tilde e_i &= - 2\g_{n+1} \g_i X - 2\g_i\left(\sum_{j=2}^m \g_{n+j}e_{n+j} +
\sum_{j=1}^m \g_je_j\right) + e_i,\cr
\tilde e_{n+i}&= - 2\g_{n+1}\g_{n+i} X - 2 \g_{n+i}\left(\sum_{j=2}^m
\g_{n+j}e_{n+j} +
\sum_{j=1}^m \g_j e_j\right) + e_{n+i}.\cr}$$
So we have
$$\eqalign{\g_i \tilde X - \g_{n+1} \tilde e_i &= \g_i X - \g_{n+1} e_i,\cr
\g_{n+i} \tilde X -\g_{n+1} \tilde e_i &= \g_{n+i} X - \g_{n+1}e_i.\cr}$$
Let $$\psi_i=\arctan(\gamma_{n+1}/\gamma_{i}), \quad {\rm and\,\,}\quad
\phi_i=\arctan(\g_{n+1}/\g_{n+i}).$$ Then the above equations are
\refeq[jo]$$\eqalign{\cos\psi_i\,
X-\sin\psi_i\,e_i &=\cos\psi_i\,\tilde X-\sin\psi_i\,\tilde e_i,\cr
\cos\phi_i X -\sin \phi_i e_{n+i} &= \cos \phi_i \tilde X -\sin \phi_i \tilde e_{n+i}.\cr}$$
Geometrically, this means  that the geodesic of $S^{n+m-1}$ at
$X(x)$ in the direction $e_i (x)$ intersects the geodesic of
$S^{n+m-1}$ at $\tilde X(x)$ in the direction $\tilde e_i (x)$ at
a point equidistant to $X(x)$ and $\tilde X(x)$.  So $P$ is a
Ribaucour Transformation. \qed

\ms

Note that the distance $\phi_i$ and
$\psi_i$'s in \refjo{} in the proof of Theorem \refzo{} satisfy 
$${\sum_{i=1}^n \cot^2\psi_i(x)\over 1+ \sum_{i=2}^m \cot^2\phi_i(x)} = s^2,$$
for all $x$.

The following theorem gives the geometric transformation for flat submanifolds in
$R^{n+m}$ corresponding to the action of $g_{s,\pi}$ on the space of solutions of the
$G_{m,n}$-system I.   The proof  is similar to that of Theorem
\refzo{}.

\refclaim[zq] Theorem.  Let $X$ be a local isometric immersion of $R^n$ in $R^{n+m}$
with flat and non-degenerate normal bundle, and $(x_1,
\cdots, x_n)$ a line of curvature coordinate system. Let $(A_1, F)$ be the corresponding
solution of the $G_{m,n}$-system I \refcn{}, 
$b, g$ as in Theorem \refbt{} (i) and (ii), and $E^I$ a frame of $(A_1, F)$ with
$E^I(x,1)=g(x)$.  Let $ A, G$ be as in Proposition \reffm{},
and $E=E^I\pmatrix{A&0\cr 0& I_n\cr}$ a frame for the solution
$\xi=\pmatrix{F\cr G\cr}$ of \refub{}.  Let $g_{s,\pi}$, $\tilde W, \tilde Z$,
$\hat W=\tilde W/\N \tilde W\N$, $\hat Z=\tilde Z/\N \tilde Z\N$ be as in Theorem
\refzc{} for the solution
$\xi$, $E$ a frame of $\xi$, and $$(\tilde F, \tilde G, \tilde E^\na)=
g_{s,\pi}\cdot (F,G,E), \quad (\tilde {A^\na_1}^I, \tilde F,
\tilde {E^\na}^I)= g_{s,\pi}\cdot (A_1, F, E^I)$$ as in section 9.  Write
$$\eqalign{E^I(x,1) &= (e_{n+1}(x), \cdots, e_{n+m}(x), e_1(x),
\cdots, e_n(x)),\cr \tilde {E^\na}^I(x,1) &= (\tilde e_{n+1}(x),
\cdots, \tilde e_{n+m}(x), \tilde e_1(x), \cdots,\tilde
e_n(x)).\cr}$$
Set 
$$\eqalign{&\gamma = (\gamma_{n+1},\ldots ,\gamma_{m+n},
\gamma_{1} ,\ldots ,\gamma_{n}):=\pmatrix {\sin\rho\ \hat{W}^t A^t
& \cos\rho \hat{Z}^t}, \cr
& \eta(x) =E^I(x,1) \gamma(x)^t, \quad {\rm where\ \ } \rho= \cot^{-1} s.\cr}$$
Then:
\item {(1)} There exists $\phi$ such that
$\phi_{x_i}=-s  b_i  \tilde{z}_i$ for all $i=1,\ldots ,n$, where
$\tilde z_i$ is the i-th coordinate of $\tilde Z$.
\item {(2)} Let $\tilde X =X+{2\cos\rho\, \phi \over \vert \vert \tilde{W}
\vert \vert } \ \eta$.  Then $\tilde X$
is again an immersed flat submanifold in
$R^{n+m}$ with flat and non-degenerate normal bundle.
\item {(3)} The bundle morphism $P(e_\a(x))= \tilde e_\a(x)$ for $n+1\leq \a\leq
n+m$ is a Ribaucour Transformation covering the map $X(x)\mapsto
\tilde X(x)$.
\item {(4)} The solution of \refcn{} corresponding to $\tilde X$ is $(\tilde
{A^\na_1}, \tilde F)$.

\proof

(1) By \refgs{}, we have $(\tilde z_j)_{x_i}= f_{ij} \tilde z_j$.
Theorem \refbt{} implies that $(b_i)_{x_j}= f_{ij}b_j$ for
$i\not= j$.  A direct computation gives $(b_i\tilde
z_i)_{x_j} = (b_j\tilde z_j)_{x_i}$.  So $\phi$ exists.

\ss

(2) It follows from the definition of $\eta$ and equation \refgs{}
for $\tilde{w}_i ,\tilde{z}_j$ that
$$d\eta = -s \left( \sum_{i=1}^n \hat{z}_i
 \hat{w}_i dx_i \right) \eta +{s\over \cos \rho}
 \sum_{i=1}^n \hat{w}_i dx_i e_i$$ and $$ d\left( {1\over \vert \vert
\tilde{W} \vert \vert}\right) = -{s\over \vert \vert \tilde{W}
\vert \vert} \left( \sum_{i=1}^n \hat{z}_i \hat{w}_i dx_i \right)
.$$ A straightforward calculation gives $$d \tilde{X}  = dX +
{2\cos \rho \over s}d\left( {\phi \over \vert \vert \tilde{W}
\vert \vert}\eta \right) =\sum_i \tilde{b}_i dx_i \tilde{e}_i,$$
where $\tilde b_i=b_i+ {2\phi \hat w_i\over\vert \vert \tilde W
\vert \vert}$. Hence $\tilde{X}$ is a submanifold with the stated
properties. \ss

(3)  Set $\l=1$ in equation \reffq{} to get $$\tilde{e}_i  =e_i
-2\gamma_i \eta.$$ So we have $$\tilde{X}+ {\cos \rho \phi \over
\gamma_i s \vert \vert \tilde{W} \vert \vert} \tilde{e}_i =X+
{\cos \rho \phi \over \gamma_i s \vert \vert \tilde{W} \vert
\vert} e_i,$$ which implies that $P$ is a Ribaucour Transformation
covering $X\mapsto \tilde X$.

\ss (4) follows from  (2) and (3). \qed

\ss
The following theorem gives the geometric transformation for local isometric immersions
of $S^n$ in $S^{m+n}$ corresponding to the action of
$g_{s,\pi}$ on the space of solutions of the partial
$G_{m,n+1}$-system I.

\refclaim[zs] Theorem.  Let $X$ be a local isometric immersion of $S^n$ in $S^{n+m}$
with flat and non-degenerate normal bundle, and
$(x_1, \cdots, x_n)$ a line of curvature coordinate system. Let $(A_1, F, b)$ be the
solution of the partial $G_{m,n+1}$-system I \refcl{} corresponding to $X$, 
$$g=(e_{n+1}, \cdots, e_{n+m}, e_1,
\cdots, e_n, X)$$  as in Theorem \refcb{}, $E^I$ the frame of $(A_1, F, b)$ such
that $E^I(x,1)= g(x)$, and $$(\tilde {A_1^\na}, \tilde F, \tilde b,
\tilde {E^\na}^I)= g_{s,\pi}\cdot (A_1, F, b, E^I)$$ as in section 9.  Write
$$\tilde{E^\na}^I(x,\l) = (\tilde e_{n+1}(x), \cdots, \tilde
e_{n+m}(x), \tilde e_1(x),\cdots, \tilde e_n(x), \tilde X(x)).$$
Then:
\item {(i)} $\tilde X$ is a local isometric immersion of $S^n$ in $S^{n+m}$ with flat
and non-degenerate normal bundle, $x$ is line of curvature
coordinates, and $\{\tilde e_\a\}_{\a=n+1}^{n+m}$ is a parallel
normal frame.
\item {(ii)} The solution of \refcl{} corresponding to $\tilde X$ is
$(\tilde {A_1^\na}, \tilde F, \tilde b)$.
\item {(iii)} The bundle morphism $P(e_\a(x))= \tilde e_\a(x)$ for $n+1\leq
\a\leq n+m$ is a Ribaucour Transformation that covers the map
$X(x)\mapsto \tilde X(x)$.

\ms
Next we give the geometric transformation for n-tuples in $R^m$ of type $O(n)$
corresponding to the action of
$g_{s,\pi}$ on the space of solutions of the $G_{m,n}$-system II. 

 \refclaim[zv]Theorem.  Let $\xi=\pmatrix{F\cr G\cr}$ be a solution of the
$G_{m,n}$-system \refub{}, $E$ a frame of $\xi$, $E(x,0)=\pmatrix{A(x)&0\cr 0&
B(x)\cr}$,  $(F,G,B)$ the corresponding solution of the $G_{m,n}$-system II \refcu{},
and
$$(\tilde F, \tilde G, \tilde {B^\na}, \tilde {E^\na}^{\II}) =
g_{s,\pi}\cdot (F,G,B, E^{\II}), \quad \tilde {A^\na} =
g_{s,\pi}\cdot A.$$ 
 Let $e_i$ and $\tilde e_i$ denote the
$i$-th columns of $A$ and $\tilde A^\na$ respectively. 
Then:
\item {(i)}  $$\eqalign{&{\p E\over\p \l}(x,0)
E^{-1}(x,0) = \pmatrix{0&X(x)\cr -X(x)^t& 0\cr},\cr &{\p \tilde
{E^\na}\over\p \l }(x,0)\tilde {E^\na}^{-1}(x,0) =
\pmatrix{0&\tilde X(x)\cr -\tilde X^t(x)& 0\cr}\cr}$$ for some $X$
and $\tilde X$. 
\item {(ii)} $X=(X_1, \cdots, X_n)$ and $\tilde X=(\tilde X_1, \cdots, \tilde
X_n)$ are  n-tuples in $\reals^m$ of type
$O(n)$ such that $\{e_\a\}_{\a=n+1}^m$ and $\{\tilde
e_\a\}_{\a=n+1}^m$ are parallel normal frame for $X_j$ and $\tilde
X_j$ respectively for all $1\leq j\leq n$.
\item {(iii)} The solutions of the $G_{m,n}$-system II \refcu{} corresponding to $X$ and
$\tilde X$ as given in Theorem \refdc{} are $(F, G, B)$ and $(\tilde F,
\tilde G, \tilde {B^\na})$ respectively.
\item {(iv)} The bundle morphism $P(e_\a(x))= \tilde e_\a(x)$ for $n+1\leq
\a\leq m$ is a Ribaucour Transformation covering the map
$X_j(x)\mapsto \tilde X_j(x)$ for each $1\leq j\leq n$.
\item {(v)} There exist maps $\phi_{ij}$ such that $X_j+\phi_{ij} e_i = \tilde X_j +
\phi_{ij} \tilde e_i$ for $1\leq j\leq n$ and $1\leq i\leq m$. 

\proof It follows from Theorem \refdc{}, Proposition \refiv{}, Corollary \refpa{} and
formula \refgo{} that $$\eqalign{ & \pmatrix{0&\tilde X\cr -\tilde
X^t&0\cr} = {\partial \tilde {E^\na} \over
\partial\lambda}(x,0)\tilde{E^\na}^{-1}(x,0) \cr
& ={ \partial E\over \partial \lambda}(x,0) E(x,0)^{-1}  + {2\over
s} E(x,0) \pmatrix{ 0 & \hat{W}\hat{Z}^t \cr -\hat{Z} \hat{W}^t &
0 \cr}E(x,0)^{-1} \cr & = { \partial E \over \partial
\lambda}(x,0) E(x,0)^{-1}  + {2\over s} \pmatrix{ 0 &  A \hat{W}
\hat{Z}^t B^{-1} \cr -B \hat{Z} \hat{W}^t A^{-1} &  0 \cr}\cr &=
\pmatrix{0&X\cr -X^t&0\cr} + {2\over s} \pmatrix{ 0 &  A \hat{W}
\hat{Z}^t B^{-1} \cr -B \hat{Z} \hat{W}^t A^{-1} &  0 \cr}.\cr}$$
So $$\tilde{X}=X+{2\over s} A \hat{W}\hat{Z}^t B^t.$$   Let $\eta
= \sum_{j=1}^m \hat{w}_j e_j$. Then $$\tilde{X}_i = X_i +{2 \over
s} \sum_{j=1}^n \hat{z}_j b_{ji} \ \eta.$$ By \refga{}, we get
$\tilde{e}_i = e_i -2 \hat{w}_i \eta$. Hence for each $1\leq l\leq
m$ and $1\leq i\leq n$ we have $$ X_j +\phi_{ij} e_i = \tilde X_j
+ \phi_{ij} \tilde e_i, \quad {\rm where\,\,} \phi_{ij} = {
\sum_{l=1}^n \tilde{z}_l b_{lj} \over s \tilde{w}_i }.$$    This
finishes the proof. \qed

\ss

\refpar[zx] Example. Recall that when $n=2$ and $m=3$, the $G_{3,2}$- system
II \refcu{} is \refdb{}, which is the Gauss-Codazzi equation for surfaces in $R^3$
parametrized by spherical line of curvature coordinates.   Let
$(u,r_1,r_2 )$ be a solution of 
\refdb{}, and
$(F,G,B)$ defined by \refdk{} the corresponding solution of \refcu{}. 
Let $\hat W,\hat Z$ as in Theorem \refzc{} for the solution
$\xi=\pmatrix{F\cr G\cr}$ and
$E$ a frame, and $(\tilde F, \tilde G, \tilde B^\na)=
g_{s,\pi}\cdot (F,G,B)$.  Then 
$$\tilde F=\pmatrix{0&\tilde
u_{x}\cr -\tilde u_{y}&0\cr}, \quad \tilde G=(\tilde r_1, \tilde r_2),
\quad \tilde B^\na= \pmatrix{\cos \tilde u& \sin\tilde u\cr \sin
\tilde u& -\cos\tilde u\cr}$$ for some solution $(\tilde u, \tilde
r_1, \tilde r_2)$ of \refdb{}.  To see this, we write 
$$\hat Z^t= (\hat z_1,
\hat z_2)= (-\sin \a, \cos\a)$$ 
for some function $\a$.  It
follows from \refga{} that 
$$\eqalign{ \tilde{B^\na}& =B (I-2
\hat{Z}\hat{Z}^t) = \pmatrix{\cos u & \sin u \cr -\sin u & \cos u
\cr} \pmatrix{\cos 2\alpha & \sin 2\alpha \cr \sin 2\alpha & -\cos
2\alpha \cr} \cr & = \pmatrix{\cos (-u+2\alpha ) & \sin
(-u+2\alpha ) \cr \sin (-u+2\alpha ) & -\cos (-u+2\alpha )
\cr}.\cr}$$
  Let $X=(X_1, X_2)$ and $\tilde X=(\tilde X_1, \tilde
X_2)$ be the 2-tuples in $R^2$ of type $O(2)$ corresponding to $(F,G,B)$ and $(\tilde
F, \tilde G, \tilde B^\na)$ as in Theorem \refzv{}. Then  the first
and second fundamental forms of $\tilde{X}_1$ and $\tilde{X}_2$
are given by 
$$\eqalign{  I_1 & = \cos^2 (-u+2\alpha ) dx^2 +
\sin^2 (-u+2\alpha ) dy^2 ,\cr \II_1 & = -\tilde{r}_1 \cos
(-u+2\alpha ) dx^2 - \tilde{r}_2 \sin (-u+2\alpha ) dy^2 , \cr I_2
& = \sin^2 (-u+2\alpha )dx^2 + \cos^2 (-u+2\alpha )dy^2 ,\cr \II_2
& = -\tilde{r}_1 \sin (-u+2\alpha ) dx^2 + \tilde{r}_2 \cos
(-u+2\alpha ) dy^2 ,\quad {\rm where\,}\cr & \tilde{r}_i = r_i
-2s\hat{w}_3 \hat{z}_i = r_i -{2s\tilde{w}_3 \tilde{z}_i \over
\tilde{z}_1^2 +\tilde{z}_1^2 }.\cr}$$ 

\ms

\refpar[il] Example.  Let $(u,r_1 ,r_2
)=(0,0,0)$ be the trivial solution of \refdb{}.  It is easy to see that the 2-tuple in 
$R^3$ of type $O(2)$ 
corresponding to the trivial solution is
$$X =\pmatrix{-x&0\cr 0& -y\cr 0&0\cr}$$
and  
$$E(x,y,\l)=\pmatrix{\cos \lambda x & 0&0& -\sin \lambda x & 0 \cr
0&\cos\lambda y & 0&0& -\sin \lambda y  \cr 0&0&1&0&0 \cr
 \sin \lambda x & 0&0& \cos \lambda  x &0 \cr
 0& \sin \lambda y & 0&0&  \cos \lambda y  \cr} $$
is a frame for the trivial solution. 

Below we write down explicitly the 2-tuple
in $R^3$ of type $O(2)$ constructed by applying the Ribaucour transformation to the
trivial solution.  It follows from Theorems \refzv{} (i) and \refzc{} that 
\refeq[im]$$\pmatrix{\tilde{w}_1 \cr \tilde{w}_2 \cr \tilde{w}_3 \cr
\tilde{z}_1 \cr \tilde{z}_2 \cr}  = \pmatrix{ \cosh (sx) w_1
+\sinh (sx) z_1 \cr \cosh (sy) w_2 +\sinh (sy)z_2 \cr w_3 \cr
\cosh (sx)z_1 +\sinh (sx)w_1 \cr \cosh (sy)z_2 +\sinh (sy)w_2
\cr}.$$ 
Use Theorem \refzv{} and a direct computation to get:
\ms
\item {(i)} If $\vert z_1 \vert \geq \vert w_1 \vert$ and  $\vert z_2
\vert \geq \vert w_2 \vert$, then 
$$ \eqalign{\tilde{X}_1&= -\pmatrix{ x \cr
0\cr 0 \cr} +  {2 a \cosh (sx) \over s (a^2 \cosh^2 (sx ) +b^2
\cosh^2 (sy ))} \pmatrix{ a\sinh (sx) \cr b \sinh (sy ) \cr w_3
\cr },\cr
\tilde{X}_2&= -\pmatrix{ x \cr
0\cr 0 \cr} +  {2 b \cosh (sy) \over s (a^2 \cosh^2 (sx ) +b^2
\cosh^2 (sy ))} \pmatrix{ a\sinh (sx) \cr b \sinh (sy ) \cr w_3
\cr },\cr}$$
 where
$a= \sqrt{z_1^2 -w_1^2 }$, $b=\sqrt{z_2^2 -w_2^2 }$, and
$a^2+b^2=w_3^2 <1$.   See Figure 1. 

\ms

\item {(ii)} If $\vert w_1 \vert \geq \vert z_1 \vert$ and  $\vert
z_2 \vert \geq \vert w_2 \vert$, then
    $$\eqalign{\tilde X_1&= -\pmatrix{ x \cr 0\cr 0\cr}
+ {2 a \sinh (sx )\over s (a^2 \sinh^2 (sx ) +b^2 \cosh^2 (sy ))}
\pmatrix{ a\cosh (sx ) \cr b \sinh (sy ) \cr w_3 \cr},\cr
\tilde X_2&= -\pmatrix{ x \cr 0\cr 0\cr}
+ {2 b \cosh (sy) \over s (a^2 \sinh^2 (sx ) +b^2 \cosh^2 (sy ))}
\pmatrix{ a\cosh (sx ) \cr b \sinh (sy ) \cr w_3 \cr},\cr}$$ where
$a= \sqrt{w_1^2 -z_1^2 }$, $b=\sqrt{z_2^2 -w_2^2 }$, and 
 $0\leq b^2-a^2 = w_3^2 <1$.  See Figure 2. 

\ms
\item {(iii)} The first and second fundamental forms of $\tilde{X}_1$
are $$\eqalign{ I_1 & = \cos^2 (2\alpha ) dx^2 + \sin^2 (2\alpha )
dy^2 ,\cr \II_1 & =   -{2sw_3 \over \tau} \sin (\alpha )\cos
(2\alpha ) dx^2 + {2s w_3 \over \tau} \cos (\alpha )\sin (2\alpha
) dy^2 , \cr }$$ where $\tilde z_1, \tilde z_2$ are given by \refim{}, 
$\tan \a = -\tilde z_1/\tilde z_2$, and 
$$\tau =\{(\cosh (sx)z_1 -\sinh (sx)w_1
)^2 + (\cosh (sy)z_2 -\sinh (sy )w_2 )^2\}^{1\over 2}.$$

\ms

We have seen in Example \refdz{} that
\item {(i)} solutions of SGE correspond to $K=-1$ surfaces of $\reals^3$ up
to rigid motion,
 \item {(ii)} $(u, \sin u, -\cos u)$ is a solution of \refdb{}
if and only if $u$ is a solution of the SGE,
\item {(iii)} if $(X_1, X_2)$ is the 2-tuple corresponding
to $(u, \sin u, -\cos u)$, then $X_1$ has Gaussian curvature $-1$,
 $u$ is the corresponding solution of the SGE, and $X_2$ is the
unit normal of $X_1$. 
\ss\ni
  Below we give a condition on $s, \pi$ so that
the action of $g_{s,\pi}$ on $(X,e_3)$ also represents a 2-tuple in $R^3$ of type
$O(2)$ corresponding to a $K=-1$ surface.

\refclaim[zy] Corollary.  Let $(u ,\sin u ,-\cos u)$ be a solution
of \refdb{}, $Z=(z_1, z_2)^t$, $W=(w_1, w_2, w_3)^t$ real constant unit vectors, 
and  constant $s$ such that 
$$sw_3 =  z_1\sin u
(0,0) - z_2\cos u(0,0).$$  Let $\pi$ denote the Hermitian projection onto $(w_1,
w_2, w_3, iz_1, iz_2)^t$.  Then the new solution $g_{s,\pi}\cdot
(u, \sin u, -\cos u)$ of \refdb{} is of the form $(\tilde u, \sin
\tilde u, -\cos\tilde u)$ for some $\tilde u$.

\proof Let $\tilde Z, \tilde W, \hat Z, \hat W$ be as in Theorem
\refzc{} for the solution $\xi=\pmatrix{F\cr G\cr}$ of the $G_{3,2}$-system, where $F$
and
$G$ are defined by \refdk{}.  Since $(\tilde W, \tilde Z)$
satisfies \refgs{}, we have $$\cases{d\tilde w_3=\sin u \ \tilde
w_1 dx_1-\cos u \ \tilde w_2 dx_2, \cr d\tilde z_1 =-\tilde z_1 du
+ s \tilde w_1 dx_1, \cr d\tilde z_2 =\tilde z_1 du + s \tilde w_2
dx_2.}$$ Thus $s d \tilde w_3=\sin u \ (d\tilde z_1 +\tilde z_2
du) -\cos u \ (d\tilde z_2
 -\tilde z_1 du) =d(\sin u \ \tilde z_1 - \cos u \ \tilde z_2)$.
It follows that
 $s\hat{w}_3 =\sin u \hat{z}_1 -\cos u \hat{z}_2$. \qed

\ms

\refpar[gt] Example.  $u=0$ is the trivial solution of the SGE.  It follows from
Example \refdz{} that $(0,0,-1)$ is a solution of \refdb{}.  Use Corollary \refzy{}
and Example \refzx{} to compute $g_{s,\pi}\cdot (0,0,-1)$ to get three parameter
family of solutions of the SGE.  Some of these examples (see Figure 3) are:

\ms
\item {(1)}  If $0<s =\sin c \leq 1$, then
 \refeq[in]$$\tilde{u}=2\tan^{-1}\left( {\sin  c \sin  (y\cos
c  )\over \cos  c \cosh  (x\sin c )}\right).$$ This is the breather
solution for the SGE, and the corresponding curvature $-1$ surface in $R^3$ is
given by $$\tilde{X}=  -\pmatrix{x \cr 0 \cr 0 \cr}  +r\pmatrix{
\sinh (x\sin c )\cr \cos  y\cos (y \cos c  ) +{\rm sec }  c \sin
y\sin (y \cos c ) \cr \sin y \cos (y\cos c ) -{\rm sec } c \cos y
\sin (y\cos c ) \cr },$$ where $$ r=   {2 \cosh (x\sin c )\over
\sin c \{\cosh^2 (x \sin c  ) +\tan^2 c \sin^2 (y\cos c )\} }.$$
\ms 

\item {(2)}  If $s=\cosh c >1 $, then $$\tilde{u}=2\tan^{-1}\left(
{\cosh c \cosh (y\sinh c  )\over \sinh c \sinh (x\cosh c )}\right),
$$ and the corresponding $-1$ curvature surface in $R^3$ is $$
\tilde{X}=    -\pmatrix{  x \cr 0 \cr 0 \cr} + r \pmatrix{ \cosh
(x\cosh c
 ) \cr \cos y \sinh (y\sinh c )+{\rm sech}c \sin y \cosh (y\sinh c ) \cr \sin
y \sinh (y\sinh c )-{\rm sech} c \cos y \cosh (y \sinh c ) \cr},$$
where $$  r= {2 \sinh (x\cosh c )\over \cosh c \{ \sinh^2 (x \cosh
c  ) +\coth^2 c\, \sinh^2 (y\sinh  c  )\} }.$$

\bs


\newsection Loop group actions for $G_{m,n}^1$-Systems.\par

\ms

In this section,  we  construct the action of a rational map with two simple poles on the
space of solutions of the $G_{m,n}^1$-system explicitly. Since
the calculation and proofs are similar to those for the $G_{m,n}$-system in
section 9, we  only state the results.

 The {\it $G_{m,n}^1$-reality condition} is 
\refeq[gw]$$\cases{
\overline{g (\bar{\lambda})}=g(\lambda ),&\cr  I_{m,n+1} g(\lambda
)I_{m,n+1} = g(-\lambda ),&\cr g(\lambda )^tI_{n+m,1}g({\lambda})^{t}
=I_{n+m,1}.&\cr}$$ 
Let 
$$\eqalign{G_+&=\{g:C\to GL_C(n+m+1)\big|
\, g \, {\rm is\, holomorphic \, and \, satisfies\,
\refgw{}}\},\cr G_-&=\{g:S^2\to GL(n+m+1,C)\big| \, g\, {\rm is\,
rational\, and \, satisfies\, \refgw{}}\}.
\cr}$$
Let $\Cx^{n+m+1}$ be equipped with the  bi-linear form:
$$\li u,\ v\ri_1=\sum_{i=1}^{n+m} \bar u_i v_i - \bar u_{n+m+1}
v_{n+m+1}.$$ 
Let $W =(w_1 ,\ldots ,w_m)^t$ and $Z= (z_1 ,\ldots
,z_{n+1}  )^t $ be unit vectors in $R^m$ and the Lorentz space
$R^{n,1}$ respectively,  and $\pi $   the orthogonal projection of
$\Cx^{n+m+1}$ onto the span of $ \pmatrix{ W \cr iZ\cr}$ with
respect to $\li \ , \ \ri_1$.  So 
\refeq[zz]$$ \pi= {1\over 2} \pmatrix{ W
W^t & -i WZ^t  \cr
 iZW^t & ZZ^t \cr} I_{n+m,1}.$$
Since $Z, W$ are real vectors,  $\pi \bar{\pi}=\bar{\pi}\pi = 0$. Given non-zero $s\in
R$, define 
\refeq[zzb]$$\eqalign{ &   q_{s,\pi}(\lambda )=    \left(
\pi +  {\lambda -is\over \lambda +is} (I - \pi ) \right) \left(
\bar{\pi} + {\lambda +is\over \lambda -is}  (I - \bar{\pi} )
\right) =\cr       & {1\over \lambda^2 + s^2} \left(  \lambda^2 I
+ 2s\lambda \pmatrix{ 0 &  WZ^t J \cr -ZW^t & 0 \cr}    + s^2
\pmatrix{ I -2   WW^t & 0  \cr 0 & I
 -2 ZZ^t J\cr} \right),
\cr}$$ where $J=I_{n,1}=\diag(1, \cdots, 1, -1)$.  A direct computation implies that
$q_{s,\pi}$ satisfies the reality condition \refgw{}.  So $q_{s, \pi}\in G_-$.

\ms

\refclaim[zzc]Theorem.  Let $\xi$ be a  solution of the
$G_{m,n}^1$-system \refaaa{},
$\theta_\lambda$ its Lax connection as in \refaab{}, and $E$ a frame
for $\xi$. Let $W$ and $Z$ be unit vectors in $R^m$ and $R^{n,1}$,
\refeq[io]$$\pmatrix{\tilde{W} \cr i\tilde{Z}\cr}(x)=
E(x,-is) ^{-1}\pmatrix{W \cr iZ \cr},$$ 
and $\tilde{\pi}(x)$ the
orthogonal projection onto the span of
 $\pmatrix{\tilde{W} \cr i\tilde{Z}\cr}(x)$ with respect to $\li\ , \ \ri_1$.  
Let $0\not=s\in R$ be a constant, $\pi$ a projection onto
$\Cx\pmatrix{W\cr iZ\cr}$, and $q_{s,\pi}$ defined by
\refzzb{}.  Set $$\eqalign{\tilde{E}(x,\l)&= q_{s,\pi}(\l)E(x,\l)
q_{s,\tilde{\pi}(x)}^{-1}(\l),\cr \tilde \xi &= \xi - 2s(\hat
W\hat Z^t J)_\ast,\cr}$$ where $\hat{W}(x)= {1 \over \vert \vert
\tilde{W}(x) \vert \vert_{m}}\tilde{W}(x)$ and $\hat{Z}(x)=
 {1 \over \vert \vert \tilde{Z}(x)\vert \vert_{n,1}}\tilde{Z}(x)$, and
$(y_\ast)_{ij}= y_{ij}$ for $i\not=j$ and $(y_\ast)_{ii}=0$ for
all $i$.  Let 
$$\tilde {E^\na}(x,\l) = E(x,\l) q_{s, \tilde \pi(x)}^{-1}(\l).$$
Then
\item {(i)} $\tilde \xi$ is a new solution of \refaaa{},
\item {(ii)}  $\tilde E^\na$ is  a frame for $\tilde \xi$,
\item {(iii)} $E(x,0)=\pmatrix{A(x)&0\cr 0& B(x)\cr}$ and $ \tilde E^\na
(x,0)=\pmatrix{\tilde A^\na(x)& 0 \cr 0 & \tilde B^\na(x)\cr}$ for
some $A, B$, $\tilde A^\na, \tilde B^\na$, and
\refeq[ha]$$\cases{\tilde A^\na= A(I-2\hat W\hat W^t), &\cr \tilde
B^\na= B(I-2\hat Z\hat Z^t)I_{n,1}.&\cr}$$

\ms
When $m=n$, Theorem \refzzc{} (i) was proved in [Zh]. 
\ms

Write $\xi=\pmatrix{F\cr G\cr}$, $F=(f_{ij})$, $G=(g_{ij})$, $\tilde F=(\tilde
f_{ij})$, $\tilde G=(\tilde g_{ij})$, and
$\tilde \xi=\pmatrix{\tilde F\cr \tilde G\cr}$.  Then formula for $\tilde \xi$  in the
above theorem is
 \refeq[gy]$$\cases{ \tilde{f}_{ij} =
f_{ij} -2s \hat{w}_i \hat{z}_j \epsilon_j ,& \cr \tilde{g}_{ij} =
g_{ij} -2s \hat{w}_{n+1+i} \hat{z}_j \epsilon_j,& }$$ where
$\e_j=1$ for $1\leq j\leq n$ and $\e_{n+1}=-1$.

\ms
Differentiate \refio{} to get the ODE for $(\tilde W,\tilde Z)$:

\refclaim[zzf] Corollary.  Let $\xi = \pmatrix{F \cr G \cr}$ be a
solution of the $G_{m,n}^1$-system \refaaa{},
$E$ a frame of $\xi$, and $(\tilde W, \tilde Z)$ be as in Theorem
\refzzc{}.  Then $(\tilde W, \tilde Z)$ is a solution of the
following system 
\refeq[gx]$$  \cases{(\tilde{w}_i )_{x_i }
=-\sum_{j\not= i }^{n+1} f_{ji}\tilde{w}_j -\sum_{j=1}^{m-n-1}
g_{ji}\tilde{w}_{j+n+1} +s\tilde{z}_i \epsilon_i , \quad i\leq n+1&\cr 
(\tilde{w}_i
)_{x_j} = f_{ij}  \tilde{w}_j,\quad i\leq n+1, j\not=i,& \cr
(\tilde{w}_i )_{x_j} =   g_{ij}\tilde{w}_j, \quad i>n+1, &\cr 
(\tilde{z}_i )_{x_j} =  
f_{ji} \tilde{z}_j \epsilon_i \epsilon_j, \quad j\not= i, &\cr 
(\tilde{z}_i )_{x_i}=-\sum_{j\not= i}^{n+1} f_{ij} \tilde{z}_j + s\tilde{w}_i .&} $$

\refclaim[gz] Corollary.  Let  $ \xi=\pmatrix{F \cr G \cr}$ be a
solution of the $G_{m,n}^1$-system \refaaa{}. Then system \refgx{} is always solvable.
Moreover, if $(\tilde{W},\tilde{Z})$ is a solution of \refgx{}
such that $\tilde{W}(0)$ and $\tilde{Z}(0)$ are unit vectors in
$R^m$ and $R^{n,1}$ respectively, then $\pmatrix{\tilde{F}\cr
\tilde{G}\cr}$ is a solution of  \refaaa{}, where $\tilde{F}$ and
$\tilde{G}$ are defined by \refgy{}.

\ms 

Write $$A=(A_1, A_2), \quad \tilde A^\na=(\tilde A_1^\na, \tilde
A_2^\na)$$ with $A_1, \tilde A_1^\na\in \cm_{m\times (n+1)}$ and
$A_2, \tilde A_2^\na\in \cm_{m\times (m-n-1)}$.  Then $(A_1, F)$
and $(\tilde A_1^\na, \tilde F)$ are solutions of the $G_{m,n}^1$-system I \refcc{},
and $(F,G,B)$ and $(\tilde F, \tilde G, \tilde B^\na)$ are solutions
of the $G_{m,n}^1$-system II \refaad{}, and their frames are related by 
$$\eqalign{\tilde E^\na(x,\l) &= E(x,\l) q_{s,\tilde \pi(x)}(\l)^{-1},\cr \tilde
{E^\na}^I(x,\l)& = E^I(x,\l) \left(I - {2\over
\l^2+s^2}\pmatrix{\l A \hat W\cr s\hat Z\cr} \pmatrix{\l \hat W^t
A^t & s \hat Z^t J\cr} \right), \cr \tilde {E^\na}^{\II} (x,\lambda
) &= E^{\II} (x,\lambda ) \left( I - {2 \over \lambda^2 + s^2 }
\pmatrix{
  s^2 \hat{W}\hat{W}^t & -s\lambda \hat{W}\hat{Z}^t B^t J \cr
  -s \lambda B \hat{Z}\hat{W}^t & \lambda^2 B \hat{Z}\hat{Z}^t B^t J \cr}
  \right).}$$

\ms

For the partial  $G_{m,n}^1$-system we have

\refclaim[zzk]Theorem.  Let $(F , G  , b ):R^n \to
gl_\ast(n)\times \cm_{(m-n)\times n}\times \cM_{n\times 1}$  be a
solution of the partial
  $G_{m,n}^1$-system \refkb{}, $\tau_{\lambda}$ its
Lax connection as in \refka{}, and  $E$ a frame for $(F,G,b)$. Let $W$
and $Z$ be unit vectors in $R^m$ and $R^{n,1}$ respectively,
$$\pmatrix{\tilde{W} \cr i\tilde{Z}\cr}(x)= E(x,-is)^{-1}\pmatrix{
W \cr iZ \cr},$$ and  $\tilde{\pi}(x)$  the orthogonal projection
onto the span of
 $\pmatrix{ \tilde{W} \cr i\tilde{Z} \cr} (x)$ with respect to the bilinear
form $\li\  , \ \ri_1$. Let $0\not=s\in R$ be a real number, $\pi$ the
projection onto $\Cx\pmatrix{W\cr iZ\cr}$, and $q_{s,\pi}$ as
defined in \refzzb{}.  Let $$\tilde{E}(x,\l) =
q_{s,\pi}(\l)E(x,\l) q_{s,\tilde{\pi}(x)}^{-1}(\l), $$
\refeq[hc]$$\cases{\tilde{f}_{ij} = f_{ij} - 2s  \hat{w}_i
\hat{z}_j ,&\cr \tilde{g}_{ij} = g_{ij} - 2s  \hat{w}_{n+i}
\hat{z}_j,&\cr \tilde{b}_i = b_i + 2s \hat{w}_i
\hat{z}_{n+1},&\cr} $$ where  $\hat{W}={1\over \vert \vert
\tilde{W} \vert \vert_m } \tilde{W}$, $\hat{Z}={1\over \vert \vert
\tilde{Z} \vert \vert_{n,1} } \tilde{Z}$ and $\hat w_i$ and $\hat
z_i$ are the i-th coordinate of $\hat W, \hat Z$ respectively.
Then $(\tilde F, \tilde G, \tilde b)$ is a solution of \refkb{}
and $\tilde E$ is its frame, where $\tilde F=(\tilde f_{ij})$,
$\tilde G= (\tilde g_{ij})$, and $\tilde b=(\tilde b_i)$.

\refclaim[zzl] Corollary. Let $\pmatrix{F , G ,b}$ be a solution
of the partial
 $G_{m,n}^1$-system \refkb{},  $E$ a frame of $(F,G,b)$, and
$\pmatrix{\tilde{W} \cr i \tilde{Z}}$ as in Theorem \refzzk{}.
Then $(\tilde W, \tilde Z)$ is a solution of
\refeq[he]$$\cases{(\tilde{w}_i )_{x_i } =-\sum_{j\not= i }^n
f_{ji}\tilde{w}_j -\sum_{j=1}^{m-n} g_{ji}\tilde{w}_{j+n}
+s\tilde{z}_i   & $i\leq n$ ,\cr (\tilde{w}_i )_{x_j} = f_{ij}
\tilde{w}_j & $i\leq n , j\not= i$ ,\cr (\tilde{w}_i )_{x_j} =
g_{ij}\tilde{w}_j   & $i>n$ ,\cr (\tilde{z}_i )_{x_j} =   f_{ji}
\tilde{z}_j & $j\not= i$ , $i\leq n$\cr (\tilde{z}_i
)_{x_i}=-\sum_{j\not= i}^n  f_{ij} \tilde{z}_j - b_i
\tilde{z}_{n+1}+  s\tilde{w}_i, &\cr (\tilde{z}_{n+1})_{x_j} = -
b_j   \tilde{z}_j. &}$$ Conversely, if  $(F , G ,b)$ is a solution
of system \refkb{}, then system \refhe{} is solvable. Moreover, if
$(\tilde W, \tilde Z)$ is a solution of \refhe{} such that
$\tilde{W}(0)$ and $\tilde{Z}(0)$ are unit vectors in $R^m $ and
$R^{n,1}$ respectively,  then $(\tilde{F}, \tilde{G},\tilde{b} )$
is also a solution of \refkb{}, where
$\tilde{F},\tilde{G},\tilde{b}$ are defined by formula \refhc{}.

\ms 

 Let $$\tilde E^\na(x,\l)= E(x,\l) q_{s,\tilde \pi(x)}^{-1}(\l).$$ Then
$\tilde E^\na$ is also a frame of $(\tilde F, \tilde G, \tilde
b)$.  Write $$E(x,0)=\pmatrix{A&0\cr 0&B\cr}\in O(m)\times O(n,1),
\quad \tilde E^\na(x,0)=\pmatrix{\tilde A^\na(x)&0\cr 0& \tilde
B^\na(x)\cr}.$$ Then
 $$\cases{\tilde{A^\na}=A(I  -2  \hat{ W}\hat{W}^t ), \cr \tilde B^\na= B (I  -2
\hat{Z}\hat{Z}^t J ),\cr }$$ where $J=I_{n,1}$.

Let $$E^I(x,\l) =E(x,\l)\pmatrix{A^t&0\cr 0& I_{n+1}\cr}, \quad \tilde
{E^\na}^I(x,\l)=\tilde E(x,\l) \pmatrix{\tilde A^t&0\cr 0&
I_{n+1}\cr}.$$ A direct computation gives $$\tilde {E^\na}^I(x,\l) =
E^I(x,\l) \left(I - {2\over \l^2+s^2}\pmatrix{\l A \hat W\cr s\hat
Z\cr} \pmatrix{\l \hat W^t A^t & s \hat Z^t J\cr} \right).$$ Note
that $E^I$ and $\tilde {E^\na}^I$ are frames of $(A_1, F, b)$ and
$(\tilde A_1^\na, \tilde F, \tilde b)$ respectively.  We use the
same type of notations as in previous sections:
$$\eqalign{&(\tilde F, \tilde G, \tilde b, \tilde E^\na)=
q_{s,\pi}\cdot (F,G, b,E), \cr &(\tilde {A_1^\na}^, \tilde F,
\tilde b, \tilde {E^\na}^I)= q_{s,\pi}\cdot (A_1, F, b, E),\cr
&(\tilde A^\na, \tilde B^\na)=q_{s,\pi}\cdot (A, B).\cr}$$

\bs

\newsection  Ribaucour Transformations for $G_{m,n}^1$-Systems.\par

In this section, we describe the corresponding geometric transformations on
submanifolds associated to the action of $q_{s,\pi}$ on the space of solutions of
various $G_{m,n}^1$-systems described in section 11.

\refclaim[xxt]Theorem.  Let $X$ be a local isometric immersion of $H^n$ in $H^{n+m}$
with flat and non-degenerate normal bundle,
 $(x_1,\dots,x_n)$     line of curvature coordinates, and $(A_1, F, b)$ the
solution of the $G_{m,n}^1$-system I \refaaj{} corresponding to $X$ as defined in
Theorem
\refbbb{}.  Let $E^I(x,\l)$ be a frame of $(A_1, F, b)$, and $g(x)=E^I(x,1)$.  Let
$0\not=s\in R$, $\pi$ the projection onto $\Cx\pmatrix{W\cr iZ\cr}$, $q_{s,\pi}$
defined by \refzzb{}, and $\tilde{E^\na}^I, \tilde A_i^\na$ as in section 11. Write  
$$\eqalign{&(\tilde A_1^\na, \tilde F, \tilde b,\tilde{E^\na}^I)=
q_{s,\pi}\cdot (A_1, F, b, E^I),\cr & E^I(x,1)=(e_{n+1}(x), \cdots, e_{n+m}(x),
e_1(x),\cdots, e_n(x), X(x)),\cr 
&\tilde{E^\na}^I(x, 1) =(\tilde
e_{n+1}(x), \cdots, \tilde e_{n+m}(x), \tilde e_1(x), \cdots,
\tilde e_n(x), \tilde X(x)).\cr}$$
 Then:
\item {(i)} $\tilde X$ is a local isometric immersion of $H^n$ in $H^{n+m}$
with non-degenerate, flat normal bundle,
and $\{e_\a\}_{\a=n+1}^{n+m}$ is a parallel normal frame.
\item {(ii)} The solution of \refaaj{} corresponding to $\tilde X$ given in
Theorem \refbbb{} is $(\tilde A^\na_1, \tilde F, \tilde b)$.
\item {(iii)} The bundle morphism $P(e_\a(x))= \tilde e_\a(x)$ for $n+1\leq
\a\leq n+m$ is a Ribaucour transformation covering the map
$X(x)\mapsto \tilde X(x)$.

Recall that if $v_0$ is a null vector in $R^{2n,1}$ and $c\not=0$, then by
Proposition \refdl{} that 
$$N_{v_0, c}=\{x\in ^{2n,1}\n \li x,
x\ri_1= -1, \li x,v_0\ri_1= c\}$$ 
is a totally umbilic hypersurface of $H^{2n}$, and $N_{v_0,c}$ is isometric to
$R^{2n-1}$. Next we give  a condition on $s, \pi$ so that the Ribaucour transformation
corresponding to $q_{s,\pi}$ for local isometric immersions of $H^n$ in $H^{2n}$
preserves local isometric immersions of $H^n$ in $R^{2n-1}$. 

\refclaim[xxr]Corollary. Let $v_0$ be a null vector in $R^{2n,1}$,
$X$  a local isometric immersion of $H^n$ in $N_{v_0,c}\subset H^{2n}$ with flat and
non-degenerate normal bundle, and
$W, Z$, $q_{s,\pi}$  and $\tilde X$ as in Theorem \refxxt{}.   If
$$b^t (0){W} =-sz_{n+1},$$ then $\tilde X$ also lies in some flat
totally umbilic hypersurface of $H^{2n}$.

\proof By Theorem \refbbb{} (vi), there exists some constant
vector $w\in \cM_{n\times 1}$ such that  $b=A^t w$. It follows
from \refbb{} that $$\cases{d\tilde w =-\sum(-FC_i +C_i F^t)\tilde
w dx_i
 + s\sum \tilde z_i C_i dx_i,\cr
     d\tilde z_{n+1}=-\sum b_i \tilde z_idx_i.}$$
Since $A^t dA=\sum(-FC_i+C_iF^t)dx_i$, we have $$d(b^t\tilde
W)=w^t d(A\tilde W)=w^t (dA)\tilde W + w^t A d\tilde W=b^t(s\sum
\tilde z_iC_i dx_i)=-sd\tilde z_{n+1}.$$ Hence $b^t \hat{W}
=-s\hat{z}_{n+1}$ for all $x \in R^n$.A straightforward
calculation now implies that $$\tilde{A^\na}^t w
=A(I-2\hat{W}\hat{W}^t )A^t w = b +2s \hat{z}_{n+1}
\hat{W}=\tilde{b}, $$ which finishes the proof. \qed

\refclaim[zzv]Theorem.  Let $(F,G,B)$ be a solution of the $G_{m,n}^1$-system II
\refaad{}, and  $X=(X_1, \cdots, X_{n+1})$ the 
(n+1)-tuple in
$R^m$ of type
$O(n,1)$ corresponding to $(F,G,B)$ of \refaad{} as described in Theorem
\refcca{}. Let $q_{s,\pi}$ be the rational map defined by
\refzzb{} for some constant $s$ and unit vectors $W$ in $R^m$ and $Z$ in $R^{n,1}$,
and
$$(\tilde F,
\tilde G,
\tilde B^\na, \tilde A^\na)= q_{s,\pi}\cdot (F,G,B, A).$$ Let
\refeq[hu]$$\tilde X= X+{2\over s} A\hat W\hat Z^t B^t J =(\tilde
X_1, \cdots, \tilde X_{n+1}),$$ where $J=I_{n,1}$. Then:
\item {(i)} $\tilde X$ is a 
(n+1)-tuple in $R^m$ of type $O(n,1)$.
\item {(ii)} The solution of \refaad{} corresponding to $\tilde X$ is
$(\tilde F, \tilde G, \tilde B^\na) = q_{s,\pi}\cdot (F,G, B)$.
\item {(iii)} Write $A=(e_1, \cdots, e_m)$, $\tilde A^\na=(\tilde e_1,
\cdots, \tilde e_m)$.  Then $\{e_i\}_{i=1}^n$ and $\{\tilde
e_i\}_{i=1}^n$ are principal curvature frames, and
$\{e_\a\}_{\a=n+1}^m$ and $\{\tilde e_\a\}_{\a=n+1}^m$ are
parallel normal frames of $X_j$ and $\tilde X_j$ respectively.
\item {(iv)} The bundle morphism $P(e_\a(x))= \tilde e_\a(x)$ for $n+2\leq
\a\leq m$ is a Ribaucour transformation covering the map
$X_j(x)\mapsto \tilde X_j(x)$ for each $1\leq j\leq n$.  
\item {(v)} There exist smooth functions $\phi_{ij}$ such that \refeq[hz]$$X_j
+ \phi_{ij} e_i = \tilde X_j + \phi_{ij}\tilde e_i$$ for all
$1\leq j\leq n$ and $1\leq i\leq n$.

\bs

\newsection Darboux Transformations for $G_{m,1}^1$-Systems.\par

Let $M, \tilde M$ be two surfaces in $R^m$ with flat and
non-degenerate normal bundle, and $P:\nu(M)\to \nu(\tilde M)$ a
Ribaucour transformation that covers $\ell:M\to \tilde M$.  If in
addition $\ell$ is a conformal diffeomorphism, then  $P$ is called
a {\it Darboux transformation\/}. Such transformations for
surfaces in $R^3$ were studied by Darboux and Bianchi and in $R^m$ by Burstall
[Bu].  In this section, we show that the transformation constructed in Theorem
\refzzv{} for  2-tuples in $R^3$ of type $O(1,1)$ is a
Darboux transformation for isothermic surfaces.

\refclaim[xxs]Theorem. Let $(Y_1, Y_2)$ be an isothermic pair in
$R^m$ corresponding to the solution $(u,G)$ of \refdq{}, and  $\xi=\pmatrix{F\cr
G\cr}$ the corresponding solution of the $G_{m,1}^1$-system, where 
$F=\pmatrix{0& u_x\cr u_y&0\cr}$. 
Let
$0\not=s\in R$, $\pi$ the projection onto $\Cx\pmatrix{W\cr iZ}$,
 $q_{s,\pi}$ the rational map defined by \refzzb{},  and $\hat
W, \hat Z$ as in Theorem \refzzc{} for the solution
$\xi=\pmatrix{F\cr G\cr}$ of the $G_{m,1}^1$-system.  Let
$$B=\pmatrix{\cosh u& \sinh u\cr \sinh u& \cosh u\cr}, \quad
(\tilde {E^\na}^{\II}, \tilde A^\na, \tilde B^\na)= q_{s,\pi}\cdot
(E^{\II},A, B).$$ Write $A=(e_1,\cdots, e_m)$ and $\tilde
A^\na=(\tilde e_1, \cdots, \tilde e_m)$. Set
\refeq[xxw]$$\eqalign{\tilde{Y_1} &= Y_1+{2\over s}(\hat{z}_1
+\hat{z}_2) e^u \sum_{i=1}^m \hat w_i e_i,\cr \tilde{Y_2} &= Y_2
+{2\over s} (\hat{z}_1 -\hat{z}_2 ) e^{ -u}\sum_{i=1}^m \hat w_i
e_i.\cr}$$ Then
\item {(i)} $(\tilde Y_1,\tilde Y_2)$ is an isothermic pair,
\item {(ii)} $(\tilde u, \tilde G)$ is the solution of \refdq{}
corresponding to $(\tilde Y_1, \tilde Y_2)$, where $\tilde u=
2\a-u$, $\sinh \a=-\hat z_2$ and $\tilde G=(\tilde g_{ij})$ is
defined by \refgy{},
\item {(iii)} the fundamental forms of $(\tilde Y_1, \tilde Y_2)$ are
$$\eqalign{&\cases{\tilde I_1 = e^{-2\tilde u}(dx^2+dy^2), &\cr
\tilde{\II}_1= e^{-\tilde u}\sum_{j=1}^{m-2}(\tilde
g_{j1}dx^2-\tilde g_{j2}dy^2)\tilde e_{2+j},&\cr}\cr
&\cases{\tilde I_2=e^{2\tilde u}(dx^2+dy^2),&\cr \tilde {\II_2} =
e^{\tilde u}\sum_{j=1}^{m-2} (\tilde g_{j1}dx^2 + \tilde
g_{j2}dy^2)\tilde e_{j+2},&\cr}\cr}$$
\item {(iv)} the bundle morphism $P(e_\a(x))= \tilde e_\a(x)$ for $3\leq
\a\leq m$ covering $Y_i\mapsto\tilde Y_i$ is a Darboux
transformation for each $1\leq i\leq 2$.

\proof Let $X_1=(Y_1+Y_2)/2$, and $X_2=(Y_2-Y_1)/2$.   By
Propositions \refdp{}, \refel{} and \refed{}, $(F,G, B)$ is a
solution of the $G_{m,1}^1$-system II, and $X=(X_1,
X_2)$ is the corresponding  2-tuples in
$R^m$ of type $O(1,1)$. Let $\tilde X=(\tilde X_1, \tilde X_2)$ be
as in Theorem \refzzv{}, and $$\tilde Y_1=\tilde X_1 - \tilde X_2,
\quad \ti Y_2= \ti X_1+\ti X_2.$$ It follows from \refhu{} that
$\ti Y_1, \ti Y_2$ are given by \refxxw{}.

Since $\hat{z}_1^2 - \hat{z}_2^2 =1$,  there exists a function
$\alpha :R^2 \rightarrow R$ such that $$\hat{z}_1   =\cosh
\alpha,\quad  \hat{z}_2   =-\sinh \alpha.$$ Since $\ti u= 2\a -u$,
\refeq[kkg]$$ e^{-\tilde{u}}=e^u {\hat{z}_1 +\hat{z}_2 \over
\hat{z}_1 -\hat{z}_2 }.$$ Use \refha{}, \refkkg{} and a direct
computation to get \refeq[xxu]$$\eqalign{\tilde{B^\na} &=B (I-2
\hat{Z}\hat{Z}^t J)   = \pmatrix{-\cosh ( 2\alpha -u ) &  -\sinh (
2\alpha -u) \cr  \sinh (2\alpha -u ) &  \cosh ( 2\alpha -u)
\cr}\cr &=\pmatrix{-\cosh \ti u& -\sinh \ti u\cr \sinh \ti u&
\cosh \ti u\cr}.\cr}$$ So (i)-(iii) follows.

To prove (iv), we first note that the map $Y_i\mapsto \tilde Y_i$ is conformal
because $(x,y)$ are isothermic coordinates for $Y_i$ and $\ti Y_i$. It remains to
prove that $P$ is a Ribaucour transform.  To see this, we use 
\refhz{} to get 
$$\eqalign{&Y_1 + (\phi_{i1}-\phi_{i2})e_i =
\tilde Y_1 + (\phi_{i1}-\phi_{i2})\ti e_i,\cr & Y_2+
(\phi_{i1}+\phi_{i2})e_i = \ti Y_2 + (\phi_{i1} + \phi_{i2}) \ti
e_i\cr}$$ for each $1\leq i\leq m$.  So $P$ is a Darboux
transformation. \qed

\ms

Let $(Y_1, Y_2)$ be an isothermic pair in $R^3$ corresponding to
the solution $(u, G)$ of \refdq{}, and $(\ti Y_1, \ti Y_2)$ as in
Theorem \refxxs{}.  Write $G=(g_1, g_2)^t$. Then the mean
curvature for $Y_i$ and $\ti Y_i$ are given as follows:
\refeq[wwa]$$\eqalign{& H_1   = {-g_{1}-g_{2}\over e^u }\cr & H_2
={-g_1+g_2\over e^{-u}}\cr &\tilde{H_1} = {\tilde g_1-\tilde
g_2\over e^{ -\tilde{u}}}   = { \hat{z}_1 -\hat{z}_2\over
e^u(\hat{z}_1 +\hat{z}_2 ) }\{ g_1-g_2 -{2s \hat w_3 (\hat{z}_1
+\hat{z}_2) } \} \cr & \tilde{H_2}   = {
\tilde{g}_1+\tilde{g}_2\over e^{\tilde{u}}}
  = {e^u(\hat{z}_1 +\hat{z}_2) \over \hat{z}_1 -\hat{z}_2 }\{
  g_1+g_2
-{2s \hat w_3 ( \hat{z}_1 -\hat{z}_2) } \} . \cr }$$

\ms 
\refpar[ij] Example. 
 A plane in $R^3$ is the isothermic surface corresponding to the trivial solution
$(0,0,0)$ of \refdq{}.  Let $W=(w_1, w_2,
w_3)^t$, $Z=(z_1,z_2)^t$ be unit vectors in $R^3$ and $R^{1,1}$ respectively. 
Use Theorem \refzzc{} and a direct computation to get $$   \pmatrix{
\tilde{w}_1 \cr \tilde{w}_2 \cr \tilde{w}_3 \cr \tilde{z}_1 \cr
\tilde{z}_2 \cr}(x,y) = \pmatrix{ w_1 \cosh sx + z_1 \sinh sx \cr
w_2 \cos sy +z_2  \sin sy \cr  w_3 \cr z_1 \cosh sx  +w_1 \sinh sx
\cr - w_2  \sin sy + z_2\cos sy \cr}.$$ 
Then the isothermic pair  $(\ti Y_1, \ti Y_2)$ corresponding to $(\ti u, \ti g_1,
\ti g_2)= q_{s,\pi}\cdot (0,0,0)$ is given by 
$$ \tilde{Y_1} =
\pmatrix{ -  x \cr   y \cr 0 \cr}
 +  {2\over s\{ \cosh
\gamma \cosh (sx )-\cos (sy  )\} } \pmatrix{  \cosh \gamma \sinh
(sx ) \cr \sin  (sy
 ) \cr -\sinh \gamma\cr},$$
$$ \tilde{Y_2} = \pmatrix{ -  x \cr   -y \cr 0 \cr}
 +  {2\over s\{ \cosh
\gamma \cosh (sx )+\cos (sy  )\} } \pmatrix{  \cosh \gamma \sinh
(sx ) \cr \sin  (sy
 ) \cr -\sinh \gamma\cr},$$
where $\cosh \gamma =  {\sqrt{z_1^2 - w_1^2}\over \sqrt{z_2^2 +
w_2^2}}$  is a constant parameter.  See Figure 4. 

The first and second
fundamental forms of $Y_2$ are given by
$\tilde{I}=e^{2\tilde{u}}(dx^2 +dy^2 )$ and $\tilde{\II} =
e^{2\tilde{u}}(\tilde{k}_1 dx^2 +\tilde{k}_2   dy^2 )$ where
$$\eqalign{ e^{\tilde{u}}  = & {\cosh \gamma \cosh (sx ) -\cos (sy
  )\over \cosh \gamma \cosh (sx  ) + \cos (sy  )},\cr
 \tilde{k}_1  = &  { 2 s \sinh \gamma \cosh \gamma \cosh (sx
 )\over (\cosh \gamma \cosh (sx )-\cos(sy   ))^2},\cr
 \tilde{k}_2   = &  {2 s \sinh \gamma   \cos (sy
  )\over (\cosh \gamma \cosh (sx  )-\cos(sy ))^2}.
\cr}$$ 

\ms

 Suppose $Y_1(x,y)$ is an immersion of a surface in $R^3$ with constant mean
curvature $2c$ and no umbilical points, and $(x,y)$ is the canonical isothermic
coordinates.  It follows from Proposition
\refeo{} that $(u, g_1, g_2)$ is a solution of
\refdq{}, where $$g_1= e^{-u}+ce^u, \quad
g_2= -e^{-u} +ce^u.$$  In
particular, $Y_1$ is isothermic.   It follows from Theorem \refcca, Propositions
\refdp{} and
\refed{} that there exists an isothermic surface $Y_2$ such that $(Y_1,
Y_2)$ is an isothermic pair corresponding to $(u, g_1,g_2)$, and the
principal curvatures of $Y_2$ are 
$$-1- ce^{2u}, \quad  -1 + c e^{2u}.$$
So the mean curvature of $Y_2$ is
$ -2$.  If $Y_1$ is minimal, then $Y_2$ is totally umbilic with mean curvature
$-2$. Hence in this case, $Y_2$ is a standard unit sphere.  In fact, 
$Y_2=e_3$ is the standard sphere parametrized by the isothermic coordinates of
$Y_1$. 

Below we give a condition on $s, \pi$ so that the induced action of $q_{s,\pi}$ on
the space of isothermic surfaces preserves the subset of constant mean curvature
surfaces parametrized by canonical isothermic coordinates.  

\refclaim[wwb]Corollary. Let $(Y_1, Y_2)$ be a isothermic pair in $R^3$ such
that $Y_1$ has constant mean curvature $2c$, no umbilical points, and is
parametrized by the canonical isothermic coordinates.  Let
$s\not=0$ be a constant,  $W, Z$ constant unit vectors in
$R^m$ and
$R^{1,1}$ respectively, and $\tilde{Y}=(\ti Y_1, \ti Y_2) = q_{s,\pi}\cdot (Y_1,
Y_2)$ as in Theorem \refxxs{}. If $$s{w}_3 = ce^{u(0,0)} (z_1 +z_2 )
+e^{-u(0,0)}(z_1 -z_2),$$ then $\tilde{Y_1}$ has constant mean
curvature $2c$ and $\tilde Y_2$ has mean curvature $-2$.

\proof By \refhc{}, we have $$\ti g_1= g_1 - 2s \hat w_3 \hat z_1,
\quad \ti g_2 = g_2 + 2s \hat w_3 \hat z_2.$$ Write $\hat z_1=
\cosh \a$ and $\hat z_2= -\sinh \a$ for some function $\a$.
 It follows from the differential equation \refgx{} for $\tilde W, \ti Z$
that the derivative of \refeq[ia]$$s\tilde{w}_3 -ce^{u}
(\tilde{z}_1 + \tilde{z}_2 ) -e^{-u}(\tilde{z}_1 -\tilde{z}_2 )$$
is zero. Since $\tilde W(0)=W$ and $\tilde Z(0)= Z$,   \refia{} is
identically zero on $R^2$. Use Theorem \refxxs{}, formulas \refxxu{},
\refia{}, and a direct computation to prove that  the mean
curvature of $\ti Y_1$ and $\ti Y_2$  are $2c$ and $-2$
respectively. \qed

\ms 

Similar computation gives

\refclaim[wwc]Corollary. Let  $Y_1$ be a minimal surface without umbilic points in
$R^3$ parametrized by the canonical isothermic coordinates, and $\tilde{Y_1}$ as in
Corollary 
\refwwb{}. If
$s{w}_3 = e^{-u(0,0)}({z}_1 -{z}_2)$, then $\tilde{Y_1}$ is minimal.

\ms \refpar[jo] Example.

 Let $Y_1=\pmatrix{-x \cr \sin y \cr 1-\cos y}$ and $Y_2=\pmatrix{-x \cr -\sin y \cr -1+\cos
 y}$. This pair of cylinders $(Y_1,Y_2)$ is an isothermic pair with
$(u, g_1, g_2)=(0, 0, -1)$ as the corresponding solution to \refdq{}.
Using Theorem \refzzc{} and putting $s=\sinh c$, we obtain
$$   \pmatrix{
\tilde{w}_1 \cr \tilde{w}_2 \cr \tilde{w}_3 \cr \tilde{z}_1 \cr
\tilde{z}_2 \cr}= \pmatrix{ w_1 \cosh (x\sinh c) + z_1 \sinh (x\sinh c) \cr
w_2 \cos (y\cosh c) +w_3 \cosh c \sin (y\cosh c) \cr
-w_2\sech c \ \sin (y\cosh c) +w_3 \cos (y\cosh c) \cr z_1 \cosh (x\sinh c)  +w_1 \sinh (x\sinh c)
\cr \sinh c \ w_2\sech c \ \sin (y\cosh c) -\sinh c \ w_3 \cos (y\cosh c)  \cr}.$$ Let  $(\tilde Y_1, \tilde Y_2)$
be the isothermic pair with constant mean curvatures corresponding to $(\tilde u, \tilde g_1,
\tilde
g_2)= q_{s,\pi}\cdot (0,0 , -1)$. Then $$ \tilde{Y_1} =
\pmatrix{-x \cr \sin y \cr 1-\cos y}
 +r_1\pmatrix{a \sinh
(x\sinh c) \cr -\cosh c\ \cos y \cos (y\cosh c) -\sin y \sin (y\cosh c)\cr
-\cosh c\ \sin y \cos (y\cosh c) +\cos y \sin (y\cosh c)\cr},$$
$$ \tilde{Y_2} = \pmatrix{-x \cr -\sin y \cr -1+\cos y}
 +r_2 \pmatrix{a\sinh
(x\sinh c) \cr -\cosh c\ \cos y \cos (y\cosh c) -\sin y \sin (y\cosh c)\cr
-\cosh c\ \sin y \cos (y\cosh c) +\cos y \sin (y\cosh c)\cr},$$
where $a =\{(z_1^2 - w_1^2)/(w_3^2+w_2^2\sech^2 c)\}^{1\over 2}$,
and
$$\eqalign{r_1&={2\over \sinh c \{ a\cosh (x\sinh c) +\sinh c \ \sin (y\cosh c)\}},\cr
r_2&={2\over  \sinh c \{a\cosh (x\sinh c) - \sinh c \ \sin (y\cosh c)\}}.}$$
If $a>\sinh c$ and $\cosh c$ is a rational number $n/k$, then $\ti
Y_1, \ti Y_2$ are immersed cylinder with n bubbles (see Figures 5, 6, 7). 

Note that $\tilde u$ and the principal curvatures $\tilde k_{i,1}$, $\tilde k_{i,2}$ of
$\tilde Y_i$ are given by
$$\eqalign{ e^{\tilde{u}}  = & {a\cosh (x\sinh c) + \sinh c \ \sin (y\cosh c)
\over a\cosh (x\sinh c) -  \sinh c \ \sin (y\cosh c)},\cr
 \tilde k_{1,1}  = &  -{2a\sinh c \cosh (x\sinh c) \sin (y\cosh c)
\over (a\cosh (x\sinh c) - \sinh c \ \sin(y\cosh c))^2},\cr
 \tilde k_{1,2}  = & {a^2\cosh (x\sinh c) + \sinh^2 c \ \sin (y\cosh c)
\over (a\cosh (x\sinh c) - \sinh c \ \sin(y\cosh c))^2},\cr
 \tilde k_{2,1}  = &  -{2a\sinh c \cosh (x\sinh c) \sin (y\cosh c)
\over (a\cosh (x\sinh c) + \sinh c \ \sin(y\cosh c))^2},\cr
 \tilde k_{2,2}   = & -{a^2\cosh (x\sinh c) + \sinh^2 c \ \sin (y\cosh c)
\over (a\cosh (x\sinh c) + \sinh c \ \sin(y\cosh c))^2}.\cr}
$$ 
Note that $\tilde H_1=1$ and $\tilde H_2=-1$. This does not contradict Corollary
\refwwb{} because the isothermic coordinate system $(x,y)$ is not the canonical
one for the cylinder viewed as a CMC surface. The canonical isothermic coordinates 
are $x/2, y/2$.

\bs


\newsection B\"acklund transformations and loop group factorizations.\par

Terng and Uhlenbeck ([TU2]) show that the classical B\"acklund transformation of
surfaces in $R^3$ with constant curvature $-1$ corresponds to the action of a
rational map with only one simple pole on the space of solutions of the SGE.   
In this section, we give a generalization of this result to the $G_{n,n}$-system. 
We find a rational map that satisfies the
$G_{n,n}$-reality condition and has only one simple pole, and show that the
corresponding geometric action gives rise to a B\"acklund transformation of
n-dimensional submanifolds.

B\"acklund's theorem is generalized 
by Tenenblat and Terng in [TT] to n-submanifolds in $R^{2n-1}$,
and by Tenenblat in [Ten] to submanifolds in $S^{2n-1}$ and $H^{2n-1}$.  These
generalizations arise naturally if we reformulate the classical B\"acklund
transformations in terms of orthonormal frames as follows:  Note
that $\ell:M\to \tilde M$ is a B\"acklund transformation with
constant $\o$ for surfaces $M, \tilde M$ in $R^3$, then there
exist $O(3)$-frame $e_A$ and $\tilde e_A$ such that $\ell(p)= p+
\sin \o \, e_1(p)$ and $$(\tilde e_1, \tilde e_2, \tilde
e_3)(\ell(x))= (e_1, e_2, e_3)(x)\pmatrix{1&0 & 0\cr 0& \cos\o
&-\sin\o\cr 0& \sin\o & \cos\o\cr}$$ for all $x\in M$.

\ss \refpar[ic] Definition ([Ten]). \ni Let $M, \tilde M$ be
n-dimensional submanifolds in $S^{2n-1}$. A diffeomorphism
$\ell:M\to \tilde M$ is a {\it B\"acklund transformation\/} with
constant $\o$ if there exist local orthonormal $O(2n)$-frames
$\{e_A\}, \{\tilde e_A\}$ of $M, \tilde M$ respectively such that
\item {(i)} $e_1=X$ and $\tilde e_1=\tilde X$ are the immersions,
\item {(i)} $\{e_\a\}_{\a=n+1}^{2n-1}$ and
$\{\tilde e_\a\}_{\a=n+1}^{2n-1}$ are parallel normal frames for
$M, \tilde M$ respectively,
\item {(iii)} $$\eqalign{&(\tilde X, \tilde e_{n+1},\cdots, \ti e_{2n-1},
\ti e_1, \tilde e_n)(\ell(x))\cr &\quad =(X, e_{n+1}, \cdots,
e_{2n-1}, e_1, \cdots, e_n)(x) \pmatrix{\cos\o \, I_n & -\sin\o \,
I_n\cr \sin\o \, I_n & \cos\o\, I_n\cr}\cr}$$ for all $x\in M$.

\refclaim[if] Theorem ([Ten]). If $\ell:M^n\to \tilde M^n$ is a
B\"acklund transformation
 in $S^{2n-1}$, then both $M$ and $\tilde M$ are flat.
Moreover, if $M^n$ is a flat submanifold of $S^{2n-1}$, then given
any constant $\o$ and a unit vector $v_0\in TM_{p_0}$, there
exists a submanifold $\tilde M$ of $S^{2n-1}$ and a B\"acklund
transformation $\ell:M\to \tilde M$ such that $\ell(p_0)=\cos\o
\,\ell(p_0) + \sin\o\, v_0$.

Next, we explain the relation between the
geometric B\"acklund transformation and the dressing action on the space of
$G_{n,n}$-systems.
By Theorem \refaa{}, the Gauss-Codazzi equations for flat
n-dimensional submanifold in $S^{2n-1}$ are the $G_{n,n}$-system I, which is gauge
equivalent to the $G_{n,n}$-system for $F:R^n\to gl_\ast(n)$
such that
\refeq[hw]$$\o_\l=\l\pmatrix{0&-\d\cr \d & 0\cr} +
\pmatrix{-F\d + \d F^t &0\cr 0& -F^t\d + \d F\cr}$$ is flat for
all $\l\in C$, where $\d=\diag(dx_1, \cdots, dx_n)$. So the
$G_{n,n}$-system is the equation for $F$:
\refeq[hv]$$\cases{C_i F^t_{x_j}-F_{x_j}C_i - C_j F^t_{x_i} +
F_{x_i} C_j = [C_iF^t-FC_i, C_j F^t-F C_j],&\cr C_i
F_{x_j}-F^t_{x_j}C_i - C_j F_{x_i}+ F^t_{x_i}C_j =[C_i F-F^tC_i,
C_jF-F^t C_j],&\cr}$$ where $C_i=e_{ii}=\diag(0, \cdots, 1, 0, \cdots, 0)$ as before.

Let $s\in R$ be a non-zero constant,  $\b\in O(n)$ a constant, and
\refeq[hh]$$k_{s,\b}(\l)= {1\over \l-is}\pmatrix{s\b & \l \cr -\l
& s\b^t\cr}.$$ 
Note that $k_{s,\b}$ is holomorphic at $\l=\infty$,
 $k_{s,\pi}(\infty)\not=I$, and $k_{s,\b}$ only satisfies
the $G_{n,n}$-reality condition up to a scalar
function, i.e., $$\overline{k(\bar\l)}={\l-is\over \l +
is}\,k(\l), \quad I_{n,n}k(\l)I_{n,n}=k(-\l), \quad
k(\l)^tk(\l)={\l^2+s^2\over (\l-is)^2} \, I,$$ where $k=k_{s,\b}$.
So $k_{s,\b}$ does not belong to $G_-$. 
But the following is true: (i) an element of the form $f(\l)I$ in $G_-$ lies in the center,
(ii) the dressing action of such element on the space of solutions of the
$G_{n,n}$-system is trivial, and (iii) the factorization still works. In fact,  if
$E$ is a frame of a solution $F$ of the $G_{n,n}$-system
\refhv{}, then $k_{s,\pi}(\l)E(x,\l)$ can still be
factored as $\tilde E(x,\l)k_{s,\tilde \b(x)}$ for some functions
$\tilde \b(x)$ and $\tilde E(x,\l)$ so that $E$ is holomorphic for $\l\in C$.  Hence 
we get

\refclaim[hf] Theorem.  Let $F$ be a solution of the
$G_{n,n}$-system \refhv{}, $E$ a frame of $F$, $s\in
R$ a constant, $\b\in O(n)$ a constant matrix, and  $k_{s,\b}$
defined by \refhh{}.  Write $E(x,-is)=\pmatrix{\eta_1(x)&
\eta_2(x)\cr \eta_3(x)& \eta_4(x)\cr}$ with $\eta_i\in
\cm_{n\times n}$.   Set \refeq[hi]$$\eqalign{\tilde \b
&=(i\eta_4-\b \eta_2)^{-1}(i\b \eta_1 +\eta_3), \cr \tilde E(x,\l)
&= k_{s,\b}E(x,\l)k_{s, \tilde \b(x)}(\l)^{-1},\cr \tilde F &= F^t
+ s\tilde \b_\ast,\cr}$$ where $y_\ast$ is the matrix whose
$ij$-th entry is $y_{ij}$ for $i\not=j$ and is $0$ for $i=j$. Then
\item {(i)} $\tilde F$ is a solution of the $G_{n,n}$-system
\refhv{} and $\tilde E$ is a frame of $\tilde F$,
\item {(ii)} $\tilde \b$ is a solution of
\refeq[hj]$$\cases{d\tilde\beta= \sum_{i=1}^n (-\tilde\beta  (
FC_i   - C_i F^t ) +(F^t C_i -C_i F) \tilde{\beta}  -sC_i
+s\tilde\beta C_i \tilde\beta )dx_i ,\cr
   \tilde \beta\tilde \beta^t=I}.$$

\proof First we prove that $\tilde E(x,\l)$ is holomorphic in
$\l\in C$. By definition,
 $\tilde E(x,\l)$ is holomorphic for all $\l\in C$ except at
$\l=\pm is$, and has simple poles at $is$ and $-is$.  A direct
computation implies that  the residue at $-is$ is a constant times
$$\pmatrix{\b& -iI_n\cr i I_n& \b^t\cr}\pmatrix{\eta_1&\eta_2\cr
\eta_3& \eta_4\cr} \pmatrix{\tilde \b^t & iI_n\cr -i I_n& \tilde
\b\cr},$$ which is zero by the definition of $\ti \b$. Similar
computation shows that the residue of $\tilde E(x,\l)$ at $\l=is$
is also zero.  Hence $\tilde E(x,\l)$ is holomorphic in $\l\in C$.

Let $\tilde \o_\l= \tilde E^{-1} d\tilde E$, and $\o_\l =
E^{-1}dE$.  Then $\tilde \o_\l$ is holomorphic for $\l\in C$. It
remains to prove that $\tilde \o_\l$ is of the form \refau{} for
some $\tilde F$.   But \refeq[ho]$$\tilde \o_\l= \ti E^{-1} d\ti E
= k_{s,\tilde \b} \o_\l k_{s,\tilde \b}^{-1} - dk_{s, \tilde \b}
k_{s,\tilde \b}^{-1}.$$ Because $k_{s,\tilde \b}$ is holomorphic
at $\l=\infty$ and $\o_\l$ has only simple poles at $\l=\infty$,
so is $\tilde \o_\l$.  Compare coefficients of $\l^i$ in \refho{}
to conclude that $\tilde \o_\l$ is of the form \refau{} for some
$\tilde F$ and $\tilde F$ is given by \refhi{}.  This proves (i).

 Multiply \refho{} by $k_{s,\ti \b}$ on the right to get
$$\tilde \o_\l k_{s,\ti \b} = k_{s,\ti \b} \o_\l - d k_{s,\ti
\b}.$$ Multiply the above equation by $(\l-is)$ then compare
coefficients of $\l^i$ to get the differential equation for
$\tilde \b$ in (ii).
 \qed

\refclaim[ya] Corollary. Let $ F$ be a solution  of the
$G_{n,n}$-system  \refhv{}, and $E$ a frame of
$F$. Then  system \refhj{} is solvable for $\tilde \b$. Moreover,
if $\tilde \b$ is a solution of \refhj{} with initial condition
$\tilde{\beta}(0)=\beta$, then $\tilde F= F^t+ s\tilde \b_\ast$
is a solution of \refhv{}.

Since solutions of the $G_{n,n}$-system I correspond
to flat n-submanifolds in $S^{2n-1}$, Theorem \refhf{} and Corollary \refya{} give
a method of constructing new flat n-submanifolds in $S^{2n-1}$ from
a given one. Geometrically, this gives the geometric B\"acklund
transformation constructed by Tenenblat [Ten]:

\refclaim[yc] Theorem. Let $F, E, k_{s, \b}, \tilde \b, \tilde E,
\tilde F$ be as in Theorem \refhf{}. Write
$E(x,0)=\pmatrix{A(x)&0\cr 0& B(x)\cr}$.  Let
\refeq[hn]$$\eqalign{N(x)&=E^I(x,1)= E(x,1)\pmatrix{A(x)^t&0\cr 0&
I_n\cr},\cr \tilde N (x)&= \tilde {E^\na}^I(x,1)=
N(x)\pmatrix{A&0\cr 0&I_n\cr}\pmatrix{\cos \rho \,\tilde \b^t&
-\sin\rho\, I_n\cr\sin\rho\, I_n & \cos\rho\, \tilde
\b\cr}\pmatrix{\tilde {A^\na}^t&0\cr 0 &I_n\cr},\cr}$$ where
$\rho= \arctan (1/s)$, and $\ti A^\na= A \ti \b^t$.
 Let $v_i$ and $\tilde v_i$ denote the $i$-th column of $N$
and $\tilde N$ respectively.  Then the map $v_1(x)\mapsto \tilde
v_1(x)$ is a B\"acklund transformation of n-dimensional
submanifolds in $S^{2n-1}$ with constant $\rho$ defined by
Definition \refic{}.

\proof  By \refhi{}, $$\tilde
E(x,0)= \pmatrix{\b&0\cr 0& \b^t\cr}\pmatrix{ A(x) \tilde \b^t(x)
&0\cr 0 & B(x)\tilde\b(x)\cr}.$$ Since $\pmatrix{\b&0\cr 0&
\b^t\cr}$ is constant, $\pmatrix{ A\tilde \b^t &0\cr 0 &
B\tilde\b\cr}$ is also a trivialization of the Lax connection \refhw{}
for $\tilde F$ at $\l=0$.
 By Proposition \reffm{}, $(A, F)$ and $(\tilde A, \tilde F)$ are solutions
of the $U/K$-system I and $N, \ti N$ are the corresponding frames
at $\l=1$.  It follows from Theorem \refaa{} that $v_1, \tilde
v_1$ are flat n-dimensional submanifolds of $S^{2n-1}$, $(v_2,
\cdots, v_n)$ and $(\tilde
v_2, \cdots, \tilde v_n)$ are parallel normal frames for $v$ and $\tilde
v_1$ respectively. Multiply both sides of \refhn{} by $\pmatrix{I_n& 0\cr 0&
A^t\cr}$ on the right to get \refeq[hk]$$\eqalign{ &\tilde N
\pmatrix{I_n& 0\cr 0& A^t\cr} = N\pmatrix{\cos\rho I_n & -\sin\rho
I_n\cr \sin\rho\, \tilde \b A^t & \cos \rho\tilde \b A^t\cr}\cr
&=N \pmatrix{I_n& 0\cr 0 & \tilde \b A^t\cr} \pmatrix{\cos\rho \,
I_n & -\sin\rho\, I_n\cr \sin\rho \, I_n & \cos\rho\,
I_n\cr}.\cr}$$ Let $$\tilde N_b=\tilde N \pmatrix{I_n& 0\cr 0&
A^t\cr}, \quad N_b= N \pmatrix{I_n& 0\cr 0 & \tilde
{A^\na}^t\cr}.$$ Then \refhk{} can be rewritten as
\refeq[hm]$$\tilde N_b= N_b \pmatrix{\cos\rho \, I_n & -\sin\rho\,
I_n\cr \sin\rho \, I_n & \cos\rho\, I_n\cr}.$$  The first n column
vectors of $N$ and $N_b$ are the same, the first n column vectors
of $\ti N$ and $\tilde N_b$ are the same and they are parallel
normal frames.  The last n columns of
 $N_b$ and $\ti N_b$ are tangent frames for $v_1$ and $\ti v_1$
respectively (they are not principal curvature directions).
Geometrically, \refhm{} means that the map $v_1\mapsto \tilde v_1$
is a B\"acklund transformation with constant $\rho$. \qed

\ss

\refpar[ye] Example. We apply Theorem \refyc{} to the trivial
solution $F=0$ to get explicit immersions of flat n-submanifolds
in $S^{2n-1}$.   The Lax connection of $F=0$ of the $G_{n,n}$-system
\refhw{} is
$$\theta_\lambda = \lambda \sum_i
\pmatrix{ 0 &  -C_i \cr C_i & 0} dx_i.$$ So $$E(x,\l)=\pmatrix {C
(x,\l) &  -S (x,\l)\cr
                    S (x,\l) & C (x,\l)}$$ is a frame of $F=0$,
where $C (x,\l)=\diag (\cos (\lambda x_1), \dots, \cos (\lambda
x_n))$ and $S (x,\l)=\diag(\sin(\l x_1), \cdots, \sin (\l x_n))$.
Note that $$E(x,0)= \pmatrix{I_n&0\cr 0&I_n\cr}=\pmatrix{A&0\cr 0&
B\cr},$$ so $A= I_n$. It follows from Theorem \refhf{} that
$$\eqalign{&\tilde\beta = ( p_1-\beta p_2)^{-1}( \beta p_1- p_2),
\quad {\rm where\,}\cr &p_1=\diag(\cosh(sx_1), \cdots,
\cosh(sx_n)), \quad p_2=\diag(\sinh (sx_1), \cdots,
\sinh(sx_n)).\cr}$$  Since $A=I_n$, $\tilde A^\na=A\ti \b^t= \ti
\b^t$ and formula \refhn{} gives $$\tilde N(x)= \pmatrix{C(x,1)&
-S(x,1)\cr S(x,1)& C(x,1)\cr} \pmatrix{\cos \rho \, I_n&
-\sin\rho\, I_n\cr \sin\rho \, \ti\b& \cos \rho\, \ti \b\cr}.$$
The first column of $\tilde N(x)$ gives an explicit immersion of
flat n-submanifolds in $S^{2n-1}$ with flat and non-degenerate
normal bundle.

\ms

Next we recall a generalization of B\"acklund's theorem to n-dimensional
submanifolds in $R^{2n-1}$ with constant sectional curvature -1 proved by
Tenenblat and Terng in [TT].  First we recall the definition.

\ss \refpar[ib] Definition ([TT]).  Let $M, \tilde M$ be
n-dimensional submanifolds in $R^{2n-1}$ with flat normal bundle.
A diffeomorphism $\ell:M\to \tilde M$ is called a {\it B\"acklund
transformation\/} with constant $\o$ if there exist local
orthonormal frames $\{e_A\}$ and $ \{\tilde e_A\}$ of $M, \tilde
M$ respectively such that
\item {(i)} $\{e_\a\}_{\a=n+1}^{2n-1}$ and
$\{\tilde e_\a\}_{\a=n+1}^{2n-1}$ are parallel normal frames,
\item {(ii)} $\ell(x)=x+ \sin \o \ e_1(x)$ for all $x\in M$,
\item {(iii)}
$$(\tilde e_1, \cdots, \tilde e_{2n-1})(\ell(x)) =(e_1, \cdots,
e_{2n-1})(x) \pmatrix{1&0&0\cr 0& \cos\o \, I_{n-1}& -\sin\o\,
I_{n-1}\cr 0& \sin\o \, I_{n-1}& \cos\o \, I_{n-1}\cr}$$ for all
$x\in M$.

\refclaim[ie] Theorem ([TT]).  If $\ell:M^n\to \tilde M^n$ is a
B\"acklund transformation
 in $R^{2n-1}$ with constant $\o$, then both $M$ and
$\tilde M$ have constant sectional curvature $-1$.  Moreover, if
$M^n$ is a submanifold in $R^{2n-1}$ with sectional curvature
$-1$, then given any constant $\o$ and a unit vector $v_0\in
TM_{p_0}$, there exist a submanifold $\tilde M$ of $R^{2n-1}$ and
a B\"acklund transformation $\ell:M\to \tilde M$ such that
$\ell(p_0)= p_0 + \sin \o \, v_0$.

To explain the analytic version of the above theorem, we first recall that ([TT]) if
$M$ is an n-dimensional submanifold of
$R^{2n-1}$ with sectional curvature $-1$ then the normal bundle is
flat and there exist line of curvature coordinates such that
the two fundamental forms are \refeq[hp]$$\eqalign{I&=
\sum_{i=1}^n a^2_{1i} dx_i^2, \cr \II& = \sum_{i,j=1}^{i=n, j=n-1}
a_{1i}a_{ji}dx_i^2 e_{n+j-1}\cr}$$ for some $O(n)$-valued map
$A=(a_{ij})$ and parallel orthonormal normal frame $e_{n+1}, \cdots,
e_{2n-1}$. The Gauss-Codazzi equations for $M$ are the following
system for $(A, F)$: 
\refeq[hx]$$\eqalign{&\cases{A^{-1}dA = -F\d + \d F^t,
&\cr d\w+\w\wedge\w =- \d A^t e_{11} A \d, & where \cr} \cr
& \ \ \d=\diag(dx_1, \cdots, dx_n), \ \w=\d F - F^t \d, \ {\rm and}\ \ 
e_{11}=\diag(1, 0,
\cdots, 0).\cr }$$ Equation \refhx{} is a called the
{\it generalized sine-Gordon equation\/} (GSGE).
The analytic version of Theorem \refie{} is 

\refclaim[jc] Theorem ([TT]).  Let $(A,F)$ be a solution of the GSGE \refhx{}.  Then
the following equation for $\tilde A$ is solvable:
\refeq[jd]$$d\tilde A= -\tilde A (-F^t \d + \d F) + \tilde A \d A^t D \tilde A - DA\d.$$
Moreover, 
\item {(i)} for all $j\not=k$ and $i$, we have  $(\tilde a_{ij})_{x_k}/\tilde a_{ik}=
(\tilde a_{rj})_{x_k}/\tilde a_{rk}$, which will be denoted by $\tilde f_{jk}$,
\item {(ii)} let $\tilde F$ be the $gl(n)$-valued map whose $jk$-th entry is
$\tilde f_{jk}$ for $j\not=k$ and $\tilde f_{jj}=0$, then $(\tilde A, \tilde F)$ is
again a solution of the GSGE \refhx{}.   

Next we explain the relation between GSGE and the dressing action of
$G_{n,n}^1$-systems.  It follows from Theorem \refbbb{} that local isometric
immersions of $H^n$ in $R^{2n-1}$  correspond to solutions
$(A,F,b^t)$ of the
$G_{n,n}^1$-system I such that $b^t=Ae_{11}$.  It can be easily seen that the GSGE
\refhx{} is the $G_{n,n}^1$-system \refaaa{} with $v=(F,b^t)$.  Hence GSGE has a Lax
connection
\refeq[iy]$$\o_\l=\pmatrix{\w&  -\d \l &0\cr \d \l &\tau &\d b^t\cr 0&b\d&0\cr}, \ \
{\rm where}\ \
\w= -F\d +\d F^t, \ \tau= - F^t\d + \d F.$$
However, unlike the
$G_{n,n}$-system, we are  not able to find a rational map with only one simple pole
satisfying the
$G_{n,n}^1$-reality condition up to scalar function. 
But the following element satisfies the first two conditions of the
$G_{n,n}^1$-reality condition \refgw{}:
\refeq[jf]$$g=\l\pmatrix{0&-I&0\cr I&0&0\cr 0&0&0\cr} + \pmatrix{\b &0&0\cr
0&\g&\xi^t\cr 0&\eta& \e\cr},$$
 where $\b,\g,\xi,\eta,\e$ are matrix valued functions on $R^n$. 
We explain below how to use the gauge transformation of $g$ (all entries are
assumed to be functions of $x$) of the Lax connection $\o_\l$ of the
$G_{n,n}^1$-system  so that $g\ast \o_\l=\tilde \o_\l$ for some other solution
$(\tilde F, \tilde b^t)$.    Suppose there exists a solution
$\tilde v= (\tilde F,
\tilde b^t)$ of the $G_{n,n}^1$-system \refiy{} such that 
\refeq[ix]$$g\o_\l - dg= \tilde \o_\l g,$$
where $\tilde \o$ is the Lax connection \refiy{} corresponding to $\tilde v$. 
Equate coefficient of $\l^j$ of \refix{} for each $j$ to get
\refeq[iz]$$\cases{\tilde w= \tau +\b \d - \d \g,&\cr
\tilde \tau =\w +\g \d -\d \b,&\cr
\eta = \tilde b, \ \ \xi = b, &\cr 
d\b =\b \w -\tilde w\b,&\cr
d\g= \g\tau + \xi^t b\d - \tilde \tau \g + \d \tilde b^t \eta,&\cr 
d\xi^t=  \g\d b^t -\tilde \tau\xi^t -\d \tilde b^t \e,&\cr
d\eta = \eta \tau +\e b \d - \tilde b \d \g, &\cr
d\e=0.&\cr}$$
The first and second equations of this system imply that $\b=\g^t$, and the fourth
and fifth equations imply that 
$$\d(\b^t\b- b^tb)= (\b\b^t -\tilde b^t \tilde b)\d.$$
Equate each entry of the above equation to get 
\refeq[jb]$$\b^t\b- b^tb = \b\b^t - \tilde b^t \tilde b= \D$$
for some diagonal matrix.  We make an ansatz that $\D=r I$ for some constant
$r$.  Note that  the fourth equation of \refiz{} means that $\tilde w$ is the gauge
transformation of $\w$ by $\b$.  Hence there exists a constant matrix $C$ so that
the trivialization $A$ of $\w$ and $\tilde A$ of $\tilde w$ are related by
\refeq[ja]$$CA=\tilde A \b.$$
Substitute \refja{} into \refjb{} to get
$$A^t(C^tC-e_{11})A= \tilde A^t (CC^t-e_{11})\tilde A = r I.$$
Hence $C^tC=CC^t= r I + e_{11}$.   Let $r=\cot^2 \o$. Then 
$\b= \tilde A^t D A$, where $D=\diag(\csc\o, \cot\o, \cdots, \cot\o)$.  Substitute
this to the fourth equation of \refiz{}, we get the analytic B\"acklund transformation
\refjd{}.  

\ms

Notice that the element $g$ defined as in \refjf{} can be written as
$$g= \pmatrix{\tilde A^t&0&0\cr 0& A^t&0\cr 0&0&1\cr}\pmatrix{D&-\l I& 0\cr \l
I & D & e_1^t\cr 0 & e_1& \csc\o\cr}\pmatrix{A&0&0\cr 0&\tilde A&0\cr
0&0&1\cr},$$
where $e_1=(1,0, \cdots, 0)$ and $D=\diag(\csc\o, \cot\o, \cdots, \cot\o)$. 
\bs


\newsection Permutability Formula for Ribaucour transformations.\par

Bianchi proved a Permutability Theorem for B\"acklund transformations for surfaces in
$R^3$: If
$\ell:M_0\to M_i$ is a B\"acklund transformation with constant $\o_i$ for surfaces in
$R^3$ for $i=1, 2$ and $\sin^2\o_1\not= \sin^2\o_2$, then there exist a unique
surface $M_3$ and B\"acklund transformations $\tilde \ell_1:M_2\to M_3$ and
$\tilde \ell_2:M_1\to M_3$ with constant $\o_1, \o_2$ respectively such that 
$\tilde\ell_1\circ \ell_2= \tilde \ell_2\circ \ell_1$.  Moreover, if $q_i$ is the
solution of the SGE corresponding to $M_i$ respectively, then 
$$\tan {q_3-q_0\over 4} = {s_1+s_2\over s_1-s_2} \tan {q_1-	q_2\over 4},$$
where $s_i= \tan(\o_i/2)$.  Terng and Uhlenbeck proved in [TU2] that the
Permutability theorem is a consequence of a relation between generators $h_{s,\pi}$
defined by
\reffc{} and fact that the dressing is a group action.  In this section, we find the relation
among
$g_{s,\pi}$'s and use the same proof as in [TU2] to get the Permutability Theorem
for Ribaucour transformations for submanifolds associated to $G_{m,n}$- and
$G_{m,n}^1$-systems.  

Recall that given unit vectors $W$ in $R^m$ and $Z$ in $R^n$,
$g_{s,\pi}$ was defined by $$\eqalign{g_{s,\pi}(\lambda ) &= \left( \pi +
{\lambda -is\over \lambda +is} (I - \pi ) \right) \left( \bar{\pi}
+ {\lambda +is\over \lambda -is}  (I - \bar{\pi} ) \right)\cr
&= I + {2is\pi\over \l-is} -{2is\bar\pi\over \l+is}.\cr}$$
where $\pi$ is the Hermitian projection of $C^{n+m}$ onto $
\Cx\pmatrix{ W \cr iZ \cr} $.

\refclaim[ma] Proposition.  Let $W_k$ and $Z_k$ be unit vectors in $R^m, R^n$
respectively,
$v_k^t=(W_k^t, iZ_k^t)$, and $\pi_k$  the Hermitian projections onto
$v_k$ for
$k=1,2$ respectively. Let $s_1, s_2\in R $ be constants such that $s_1^2 \not=
s_2^2$ and $s_1 s_2\not= 0$. Let $u_k$ denote the unit direction of 
$$g_{s_j,\pi_j}(-is_k)(v_k)$$ for $j\not=k$, and $\tau_k$ the  Hermitian
projection onto $u_k$.  Then:
\item {(i)} $u_k$ is of the form ${1\over \sqrt{2}}(U_k^t, iV_k^t)$
for some unit vectors $U_k\in R^n$ and $V_k\in R^m$.
\item {(ii)}
 \refeq[mb]$$ g_{s_2,\tau_2}(\lambda )g_{s_1,\pi_1}(\lambda
)=g_{s_1,\tau_1}(\lambda )g_{s_2,\pi_2}(\lambda ).$$  
\item {(iii)} $\tau_1, \tau_2$
are unique projections satisfying \refmb{}.

\proof 

(i) follows from the fact that $g_{s_j,\pi_j}$ satisfies the $G_{m,n}$-reality
condition.

(ii) Let $g_i=g_{s_i,\pi}$ and $\ti g_i= g_{s_i, \tau_i}$ for $i=1, 2$.  The residue
of $\ti g_1 g_2 g_1^{-1}$ at $\l= is_1$ is 
$$R_{is_1}=2is_1 (\tau_1 g_2(is_1)(I-\pi_1) + (I-\bar\tau_1)g_2(is_1)\bar\pi_1).$$
Since $u_1$ is parallel to $g_2(-is_1)(v_1)$, 
$$\li g_2(is_1)(v_1^\perp), u_1\ri = \li v_1^\perp,
g_2(is_1)^*g_2(-is_1)(v_1)\ri,$$ where $v_1^\perp$ is the orthogonal
complement of $v_1$.  But $g(\bar z)^*g(z)=I$.  So the above inner product is
zero, i.e., $\tau_1g_2(is_1)(I-\pi_1)=0$.  Since  $\bar u_1$ is parallel to 
$$ \overline{g_2(-is_1)(v_1)}= g_2(is_1)(\bar v_1),$$
a similar argument gives
$(I-\bar\tau_1)g_2(is_1)\bar\pi_1=0$.  Hence $R_{is_1}=0$, which implies that 
$\ti g_1 g_2 g_1^{-1}$ is holomorphic at $\l= is$. Since $\ti g_1 g_2 g_1^{-1}$
satisfies the reality condition, it is also holomorphic at $\l= -is_1$.  So 
$h_1= \ti g_1 g_2g_1^{-1}\ti g_2^{-1}$ is holomorphic at $\pm is_1$.  
Use a similar
argument to prove that $h_2= \ti g_2g_1g_2^{-1}\ti g_1^{-1}$ is holomorphic at
$\l=\pm is_2$.  But $h_1=h_2^{-1}$, which is holomorphic for all $\l\in C$ except
at $\l=\pm is_1, \pm is_2$.  Since $h_i(\l)^{-1}=h_i(\bar \l)^*$, $h_2^{-1}$ is
holomorpic at $\pm is_2$.  This proves that $h_1$ is holomorphic for all $\l\in C$. 
But $h_1(\infty)=I$.  So $h_1=I$.  This proves (ii).

(iii) Let $\s_i$ denote the projection onto $y_i$, and $\phi_i=g_{s_i,
\s_i}$. Suppose $\phi_1 g_2=
\phi_2g_1$. We want to prove $\s_i=\tau_i$.  But $\phi_1=\phi_2g_1g_2^{-1}$
implies that 
$\phi_2g_1g_2^{-1}$ is holomorphic at $\l=is_2$.  Hence the residue at $is_2$
must be zero, i.e.,
\refeq[jp]$$\s_2 g_1(is_2)(I-\pi_2) + (I-\bar\s_2)g_2(is_2)\bar\pi_2=0.$$  
Since 
$$\li v_2,\bar v_2\ri=0, \quad \li y_2, \bar y_2\ri=0,$$ \refjp{} implies that
$\s_2g_1(is_2)(I-\pi_2)=0$.  So 
$$0=\li g_1(is_2)(v_2^\perp), y_2\ri = \li v_2^\perp, g_1(is_2)^*(y_2)\ri 
=\li v_2^\perp, g_1(-is_2)^{-1}(y_2)\ri.$$
This proves that $g_1(-is_2)^{-1}(y_2)\in Cv_2$.  Hence $y_2$ is parallel to
$g_1(-is_2)(v_2)$. \qed

\ms

It is known (cf. [TU2]) that an analogue of Bianchi's Permutability Theorem can
be obtained easily from the above Proposition.  We sketch the reason here.  Let $\xi$
be a solution of the $G_{m,n}$-system, and $E$ the frame of $\xi$ with $E(0,\l)=I$. 
Let $g_j=g_{is_j, \pi_j}$, $h_j= g_{is_j, \tau_j}$ as in Proposition \refma{}. So
$h_1g_2=h_2g_1$. Since dressing $\sharp$ is an action, 
$$(h_1g_2)\sharp \xi = (h_2g_1)\sharp \xi.$$
Let $\xi_i=g_i\sharp \xi$. Then 
$$h_1\sharp \xi_2= h_2\sharp \xi_1,$$ which will be denoted by $\xi_3$.  Factor
$$g_1E = E^1 \tilde g_1, \quad g_2E= E^2\tilde g_2, \quad h_2E^1 =  E^4
\tilde h_2, \quad h_1 E^2 = E^3\tilde h_1$$ such that $E^i(x,\l)$ are holomorphic for
$\l\in C$ and
$\tilde g_i(x,\l), \tilde h_i(x,\l)$ are meromorphic for $\l\in S^2$ and are equal to $I$
at $\l=\infty$.  Then 
$$h_2g_1E= E^4\tilde h_2 \tilde g_1, \quad h_1g_2E= E^3 \tilde h_1\tilde g_2.$$
It follows from the assumption $h_1g_2=h_2g_1$ and the uniqueness of factorization
that $E^3=E^4$, and
$\tilde h_1\tilde g_2=\tilde h_2
\tilde g_1$.  So by Proposition \refma{} and Theorem \refec{}, we have
$$g_{s_1, \tilde \tau_1(x)}  g_{s_2, \tilde \pi_2(x)} = g_{s_2, \tilde \tau_2(x)}
g_{s_1, \tilde \pi_1(x)},$$ and the projections $\tilde \tau_i$ and $\tilde \pi_i(x)$ are
related the same way as $\tau_i$ and $\pi_i$.  

To summarize, let $\xi$ be a solution of the $G_{m,n}$-system, and $E$ the
frame of $\xi$ such that $E(0,\l)=I$. Let
$\pi_k$ be the projection onto $v_k$, $\xi_i= g_{s_i,\pi}\sharp \xi$, and $E_i$ frame of
$\xi_i$.  
$$\ti v_k(x)= E(x, -is_k)^{-1}(v_k), \quad \hat v_k = \ti v_k/\N\ti v_k\N,$$
and $\ti \pi_k(x)$ the projection onto $\hat v_k$.  
Let $\Phi:gl(n+m,C)\to \cm_{m,n}$ be the map defined by 
$$\Phi_{m,n}\left(\pmatrix{P&Q\cr R&S\cr}\right)= Q,$$ where a $gl(n+m)$ matrix is
blocked into $(m,n)$ blocks.  Then 
$$\xi_k= \xi - 4is\Phi_{m,n}(\hat v_k\hat v_k^*), \quad k=1, 2.$$
Suppose $s_1^2\not= s_2^2$. Let 
$$\ti u_k = g_{is_j, \ti\pi_j(x)}(-is_k)(\hat v_k), \quad k\not=j.$$
Then we get the Permutability Formula:
$$\eqalign{\xi_3 &=\xi- 4i\Phi_{m,n}(s_1\hat u_1\hat u_1^* + s_2 \hat v_2\hat
v_2^*) \cr
& =\xi - 4i \Phi_{m,n}(s_2\hat u_2\hat u_2^* + s_1\hat v_1\hat v_1^*)\cr
&=g_{is_2,\tau_2}\sharp \xi_1 = g_{is_1,\tau_1}\sharp \xi_2,\cr}$$
where $\tau_j$ is the projection onto $\hat u_j(0)$.  
We have shown in section 10 that the action $g_{s, \pi}$ on the space of solutions
of the $G_{m,n}$-system gives rise to Ribaucour transformations for submanifolds. 
Hence  we get

\refclaim[mc] Corollary.  Let $P_i:\nu(M)\to \nu(M_i)$ be the Ribaucour transformation
for flat n-submanifolds in $S^{m+n-1}$ corresponding to the action of $g_{s_i,\pi_i}$.  If
$s_1s_2\not=0$ and
$s_1^2\not=s_2^2$, then there exist unique flat n-submanifold $M_3$ in
$S^{n+m-1}$ and Ribaucour transformations $\ti P_1:\nu(M_2)\to \nu(M_3)$ and $\ti
P_2:\nu(M_1)\to \nu(M_3)$ such that $\ti P_1\circ P_2 = \ti P_2\circ P_1$. 

The Permutability Formula and Theorem for flat n-submanifolds  in
$R^{m+n}$, n-tuples in
$R^m$ of type $O(n)$, 
$G_{m,n}^1$-systems, and the submanifolds associated to the
$G_{m,n}^1$-system can be obtained exactly the same way.  So we will not write
them down.

\bs


\newsection The $U/K$-hierarchy and Finite type solutions.\par

In this section, we give a short review how the $U/K$-system fits into the
$U/K$-hierarchy of soliton equations.  We also explain the relation between the
ODE method of constructing finite type solutions of the $U/K$-system and the
dressing action.

First we give a quick review of the ZS-AKNS $n\times n$ hierarchy of commuting
flows.  Let $a_1$ be a diagonal matrix with distinct eigenvalues, and $a_1, \cdots,
a_{n-1}$ linearly independent diagonal matrices in
$sl(n,C)$.   Since $a_1$ has distinct eigenvalues, the map $\ad(a_1)(x)=[a_1,x]$
is a linear isomorphism of $sl_\ast(n,C)$, where
$sl_\ast(n,C)$ is the space of all $x\in sl(n,C)$ with zero on all diagonal
entries.   It is known (cf. [Sa], [TU1]) that given $1\leq i\leq n-1$ and a positive
integer
$j$ there exists a polynomial differential operator $Q_{a_i,j}$ of order $j-1$ for
maps
$u:R\to sl_\ast(n,C)$ such that
\item {(i)} $Q_{a_i,0}= a_i$ and $Q_{a_i,1}(u)= \ad(a_i)\ad(a_1)^{-1}(u)$,
\item {(ii)} $(Q_{a_i,j})_x + [u, Q_{a_i,j}] = [Q_{a_i, j+1}, a]$ for all $j$.

\ni  {\it The $j$-th flow\/} on $C^\infty(R, sl_\ast(n,C))$ defined by $a_j$ is the
following evolution equation
\refeq[is]$$u_t= (Q_{a_i,j}(u))_x + [u, Q_{a_i, j}(u)].$$ 
All these flows commute. The hierarchy of these commuting flows is called the
$SL(n,C)$-{\it hierarchy.\/} 
It follows from (ii) that  $u$ is a solution of \refis{} if and only if
$$\W_\l= (a_1\l + u) dx + \left(\sum_{k=0}^j Q_{a_i, k}\l^{j-k}\right) dt$$
is flat for all $\l\in C$.  

For example, when $j=1$, equation \refis{} is 
\refeq[it]$$u_t=\ad(a_i)\ad(a_1)^{-1}(u_x) + [u, \ad(a_i)\ad(a_1)^{-1}(u)].$$

 Many well-known soliton hierarchies come from restricting the flows in the
$SL(n,C)$-hierarchy to various invariant submanifolds.  
If $a_i\in su(n)$, then $\cm_1=C^\infty(R, su_\ast(n))$ is
invariant  under the j-th flow for all $j\geq 1$, where $su_\ast(n)= sl_\ast(n,C)\cap
su(n)$.  The restriction of the
$SL(n,C)$-hierarchy of flows to $\cm_1$ is called the $SU(n)$-hierarchy.  
For example, the second flow in the $SU(2)$-hierarchy is the non-linear Schr\"odinger
equation:
$$q_t= {i\over 2} \left(q_{xx} + 2\n q\n^2 q\right),\eqno({\rm NLS})$$
where $u=\pmatrix{0&q\cr -\bar q& 0\cr}$ and $a_1=\diag(i,-i)$.

Let $u:R^2\to su_\ast(n)$ be a solution of the $j$-th flow of the $SU(n)$-hierarchy,
and  $E$ the solution
for $E^{-1}dE=\W_\l$ with $E(0,\l)=I$. Then $\W_\l$ and $E$ satisfy the
$SU(n)$-reality condition.  

If $a_i\in su(n)$ is pure imaginary for all $i$, then $\cm_2=C^\infty(R, so(n))$ is
invariant under all the odd flows.  The restriction of  the
$SU(n)$-hierarchy to $\cm_2$ is the $SU(n)/SO(n)$-hierarchy.  
For example, the third flow in the $SU(2)/SO(2)$-hierarchy is the modified KdV equation:
$$q_t=-{1\over 4}\left(q_{xxx} + 6 q^2 q_x\right),$$
where $u=\pmatrix{0&q\cr -q&0\cr}$ with $q$ real and $a_1=\diag(i,-i)$.  

 If $v:R^2\to so(n)$ is a solution of the j-th flow in the $SU(n)/SO(n)$-hierarchy and
$E$ is the solution for $E^{-1}dE=\W_\l$ with $E(0,\l)=I$, then $\W_\l$ and $E$
satisfy the
$SU(n)/SO(n)$-reality condition.

We can replace $sl(n,C)$ by any complex semi-simple Lie algebra $\cg$, $su(n)$ by a
real form $\cu$ of $\cg$, and $SU(n)/SO(n)$ by a symmetric space $U/K$ to construct
$G$-, $U$- and $U/K$-hierarchies similarly.  

\ms

Let $G_+$ denote the group of holomorphic maps from $C$ to $U_C$ that
satisfies the $U/K$-reality condition \refgg{}, and $G_-$ the group of germs of
holomorphic map $g$ from 
$\in S^2$ to $U_C$  at $\l=\infty$ that satisfies the $U/K$-reality condition and
$g(\infty)=I$ as before.    Given $g\in G_-$, factor
$$g^{-1}(\l)\exp(a_1\l x+ a_i\l^j t)= E(x,t,\l)m(x,t,\l)^{-1}$$
with $E(x,t,\cdot)\in G_+$ and $m(x,t,\cdot)\in G_-$. 
It is known (cf. [TU1]) that  
\item {(i)}$E^{-1}E_x=a_1\l + u(x,t)$ and $u$
 is a solution of the $j$-th flow \refis{} defined by $a_i$,
\item {(ii)} $Q_{a_i,k}(u)$ is the coefficient of $\l^{-k}$ in  $m^{-1}a_im$. 

\ni Let $e_{a_i,j}(x)$ denote the one-parameter subgroup of $G_+$ generated by
$\eta(\l)=a_i\l^j$, i.e., $e_{a_i,j}(x)(\l)= \exp(a_i\l^j x)$. Then the
$j$-th flow defined by
$a_i$ can be viewed as the flows corresponding to the dressing action of the two
parameter abelian subgroup $\{e_{a_1,1}(x) e_{a_i,j}(t)\n (x,t)\in R^2\}$ of
$G_+$ on $G_-$.  The $U/K$-system \refat{} can be viewed as
given by the n commuting first flows in the
$U/K$- hierarchy (c.f. [TU1]).  In other words,  it corresponds to the dressing
action of the n-dimensional abelian subgroup
$$\{e_{a_1,1}(x_1)\cdots e_{a_n,1}(x_n)\n (x_1, \cdots, x_n)\in R^n\}$$
 of $G_+$  on
$G_-$.  We explain this more precisely below.  Let $g\in G_-$, and 
$$E_0(x,\l)=  \exp\left(\sum_{i=1}^n a_i\l x_i\right).$$
 The Birkhoff factorization
Theorem implies that  there exists unique $E(x,\l)$ and $m(x,\l)$ such that
\refeq[jq]$$g(\l)^{-1} E_0(x,\l) = E(x,\l) m(x,\l)^{-1}$$ 
with $E(x,\cdot)\in G_+$ and $m(x,\cdot)\in G_-$ for each $x$. 
Equation \refjq{} implies that $g E = E_0 m$.  So we have
\refeq[jr]$$E^{-1}E_{x_i} = m^{-1} a_i m \l + m^{-1} m_{x_i}$$
for each $1\leq i\leq n$.  Since the left hand side of \refjr{} is holomorphic for $\l\in C$,
so is the right hand side.   Expand $m(x, \l)$ at
$\l=\infty$: 
$$m(x,\l)= I + m_1(x) \l^{-1} + m_2(x) \l^{-2} + \cdots .$$
Then we have 
$$m^{-1}am = a + [a, m_1]\l^{-1} + \cdots.$$
For $h(\l)=\sum_{j\leq n_0} h_j\l^j$, define
$$h_+(\l)=\sum_{j\geq 0} h_j\l^j, \quad h_-(\l)=\sum_{j\leq -1} h_j \l^j.$$
Since $E^{-1}E_{x_i}$ is holomorphic for $\l\in C$ and $m^{-1}m_{x_i}=\sum_{j\leq
-1} b_j\l^{j}$, 
equation \refjr{}  implies that 
\refeq[js]$$E^{-1}E_{x_i}= (m^{-1} a_im\l)_+ = a_i\l + [a_i, m_1]$$
\refeq[jt]$$m^{-1}m_{x_i} = -(m^{-1}a_im\l)_-.$$
By Proposition \refas{}, $p(m_1)$ is a solution of the $U/K$-system \refat{}, where
$p$ is the projection onto $\cp\cap \ca^\perp$.  In fact, $p(m_1)= g\sharp 0$, the
action of $g$ at the solution $v=0$.  

\ms

Let
$\cf$ denote the space $g\in G_-$ such that 
$g(\l)^{-1}ag(\l)$ is a polynomial in $\l^{-1}$.  Solution $g\sharp 0$ is called a {\it
finite type solution\/}.  This is motivated by the finite type solutions constructed for
CMC tori in $R^3$ by Pinkall and Sterling [PiS], in $N^3(c)$ by Bobenko [Bo1], and
for harmonic maps from a torus to a symmetric space by Burstall, Ferus, Pedit and
Pinkall [BFPP].     We explain below the relation between  the dressing action and the
ODE construction of finite type solutions.

First we derive a system of equations for $m^{-1}a_1m$.  
It follows from a direct computation and \refjt{} that
$$(m^{-1}a_1m)_{x_i} = [m^{-1}a_1 m, m^{-1}m_{x_i}]=[m^{-1}a_1m,
-(m^{-1}a_im\l)_-].$$
But $[a_1, a_i]=0$ implies that $[m^{-1}a_1m, m^{-1}a_i m]=0$.  So we have
\refeq[ju]$$(m^{-1}a_1m)_{x_i}= [m^{-1}a_1m, (m^{-1}a_im\l)_+].$$

Write
$$m^{-1}(x,\l)a_1m(x,\l) =\sum_{j=0}^\infty \xi_j(x) \l^{-j}.$$ 
Then $\xi_1=[a,m_1]$ and 
$$(m^{-1}a_im\l)_+ = a_i\l+ [a_i, m_1] = a_i\l + \ad(a_i)\ad(a_1)^{-1}(m_1).$$
So equation \refju{} becomes
\refeq[jw]$$\left(\sum_{j=0}^\infty\xi_j\l^{-j}\right)_{x_i} = \left[\sum_{j=0}^\infty
\xi_j\l^{-j},\ \  a_i\l +
\ad(a_i)\ad(a_1)^{-1}(\xi_1)\right].$$ 

Since $m$ and $a_1\l$ satisfy the $U/K$-reality condition,
 $m^{-1} a_1\l m$ satisfies the $U/K$-reality condition. It follows from Proposition
\refjh{}  that the coefficient of $\l^{-j}$ of $m^{-1}a_1m$ satisfies the condition
that $\xi_j\in
\ck$ for $j$ odd and $\xi_j\in \cp$ for $j$ even. 

Compare coefficient of $\l^{-j}$ in \refjw{} for each $j$ to get 
\refeq[jva]$$(\xi_j)_{x_i} = [\xi_j, \ \ \ad(a_i)\ad(a_1)^{-1}(\xi_1)] +
[\xi_{j+1}, a_i],  \qquad 0\leq  j,$$ 
where $\xi_j\in \ck$ if $ j$ is odd,  $\xi_j\in \cp$ if $j$ is even, 
 and $\xi_0=a_1$.  
Note that for a fixed positive integer $k$,  system \refjva{} leaves the subset defined by
$\xi_j=0$ for all
$j\geq k+1$ invariant.  On this invariant subset, system \refjva{} becomes the following
system for $(\xi_1, \cdots, \xi_k)$:
\refeq[jv]$$\cases{(\xi_1)_{x_i} = [\xi_1, \  [a_i, \ad(a_1)^{-1}(\xi_1)]] + [\xi_2,
a_i],&\cr
(\xi_2)_{x_i} =[\xi_2, \  [a_i, \ad(a_1)^{-1}(\xi_1)]] + [\xi_3, a_i], &\cr
\cdots &\cr
(\xi_k)_{x_i} = [\xi_k, \ [a_i, \ad(a_1)^{-1}(\xi_1)]]. &\cr}$$
Thus we get

\refclaim[jx] Theorem. Let $k$ be a positive integer,  $a_1+\sum_{j=1}^k \eta_j
\l^{-j}= g(\l)^{-1}a_1 g(\l)$ for some $g\in \cf$, $V_1= \ck\cap \ca^\perp$,
$V_j=\ck$ if
$j$ is odd, and $V_j=\cp$ if $j$ is even.  If $(\xi_1, \cdots, \xi_k):R^n\to
V_1\times V_2\times \cdots \times V_k$ is a solution of 
\refjv{} with $\xi_0=a_1$, $\xi_{k+1}=0$, and $\xi_j(0,\cdots, 0)= \eta_j$ for
$1\leq j\leq k$, then
$\ad(a_1)^{-1}(\xi_1)$ is a solution of the $U/K$-system \refat{} and
$\ad(a_1)^{-1}(\xi_1) = g\sharp 0$.

\vfil\eject

\centerline{\bf Figure 1}
\ms
\centerline{Example \refil{} (i). A 2-tuple in $R^3$ of type
$O(2)$  obtained by applying }
\centerline{a Ribaucour trasformation to a
pair of lines}
\centerline{\BoxedEPSF{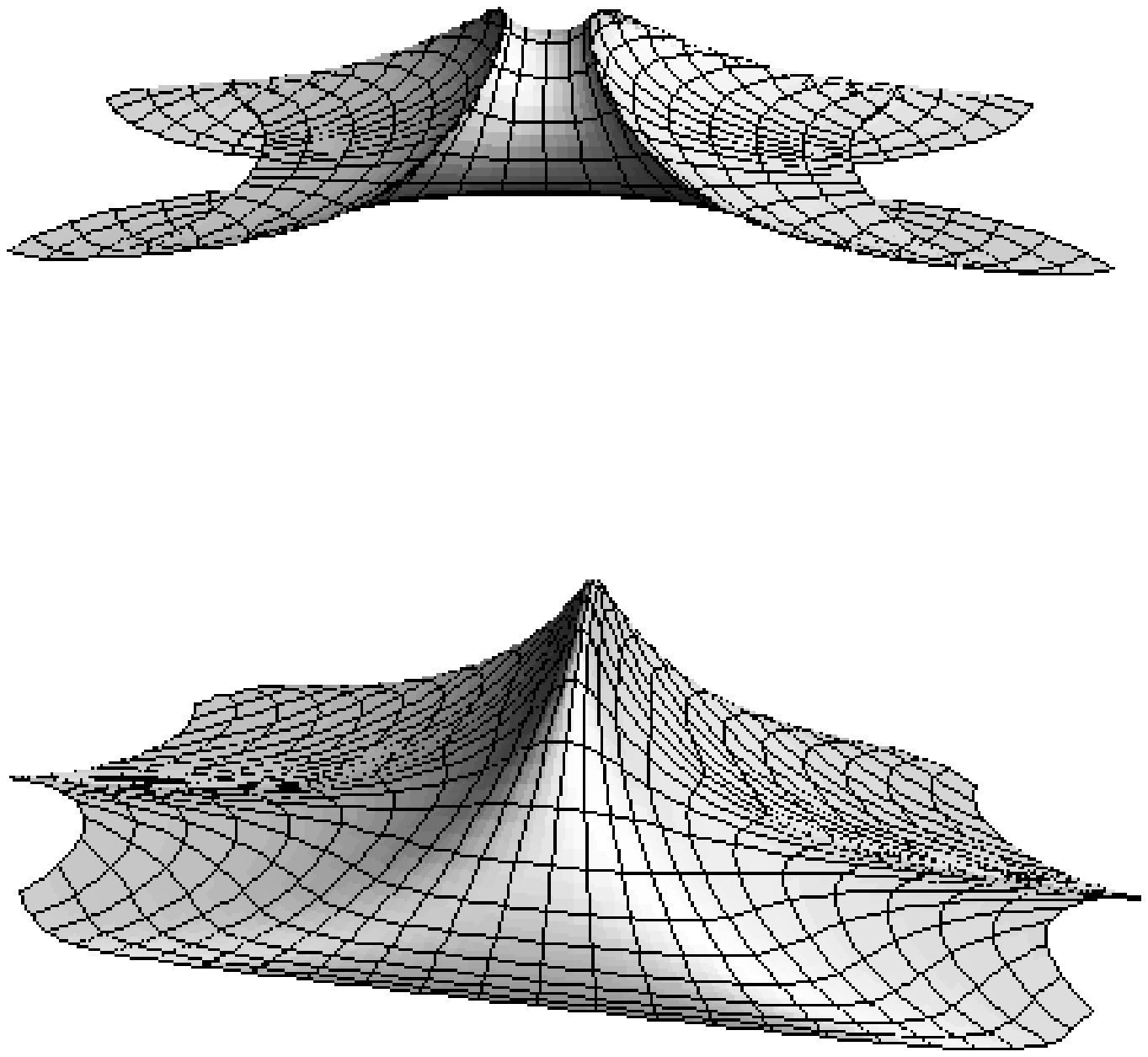 scaled 650}}
\vfil\eject

\centerline{\bf Figure 2}\ms
\centerline{Example \refil{} (ii). A 2-tuple in $R^3$ of type
$O(2)$ obtained by applying }
\centerline{a Ribaucour trasformation to a
pair of lines}
\centerline{\BoxedEPSF{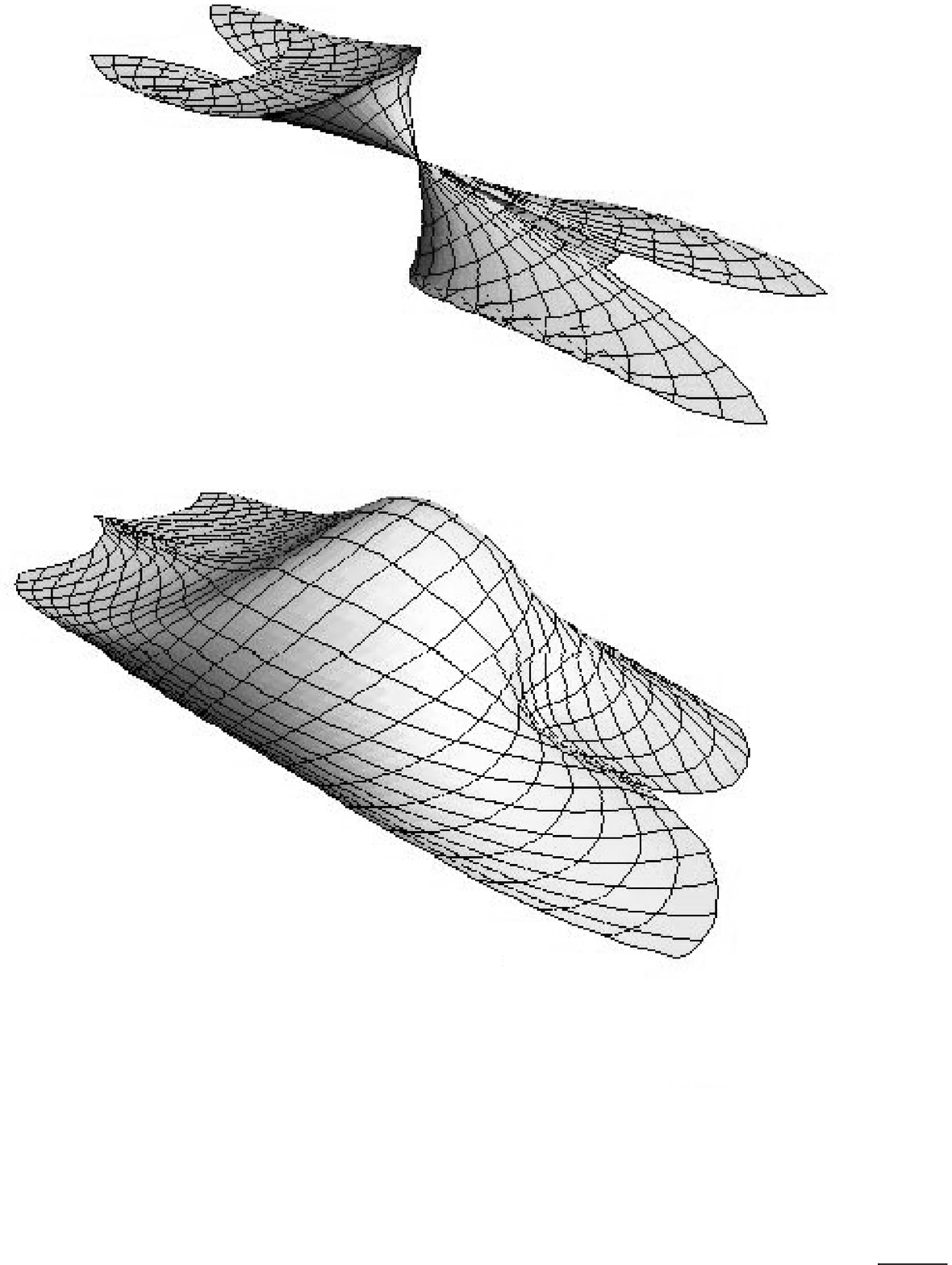 scaled 650}}
\vfil\eject

\centerline{\bf Figure 3}\ms
\centerline{Example \refgt{}. $K= -1$ surfaces, top: case (1),
bottom: case (2)}
\bs
\centerline{\BoxedEPSF{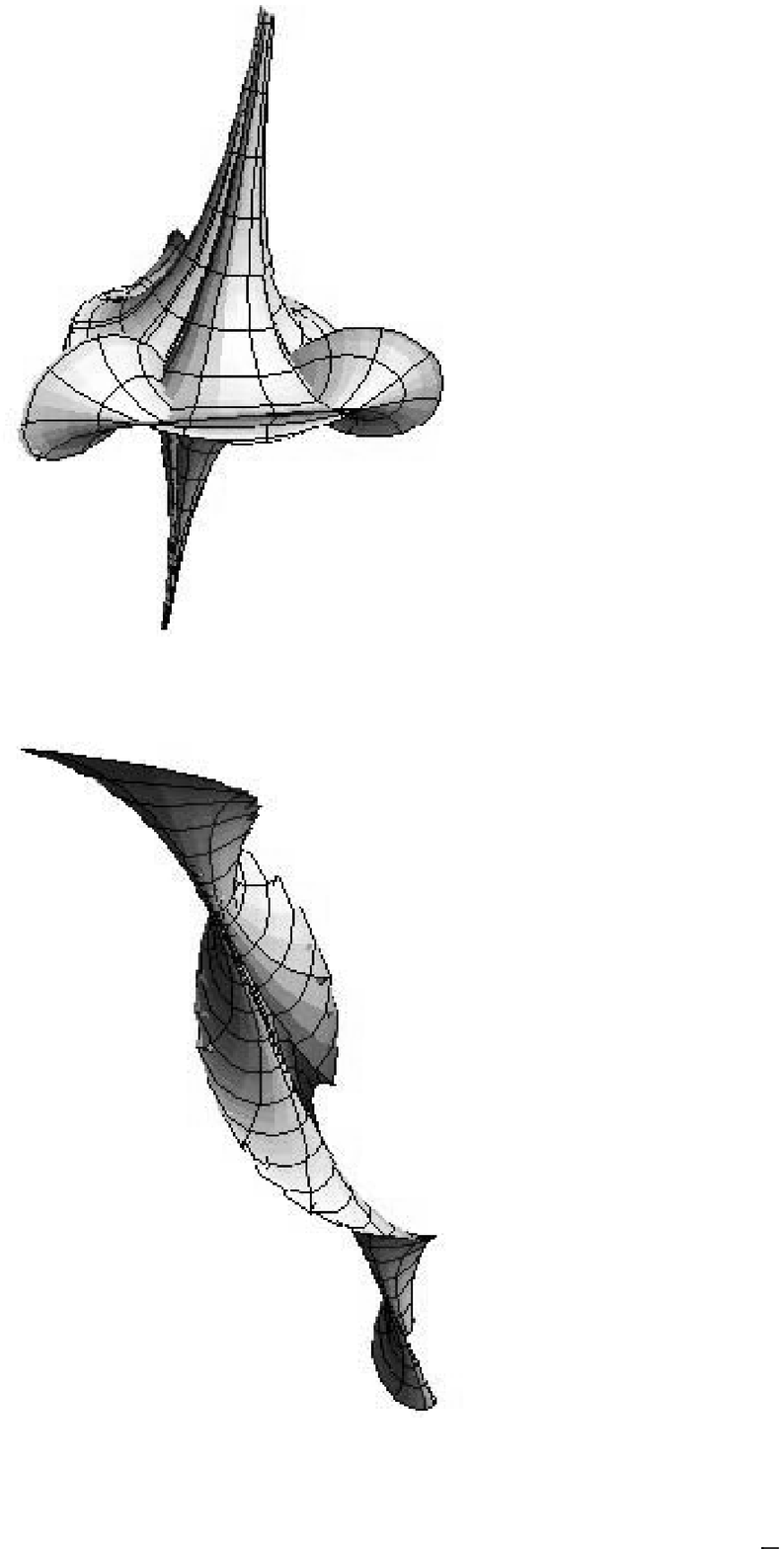 scaled 600}}
\vfil\eject

\centerline{\bf Figure 4}\ms
\centerline{Example \refij{}.  An isothermic pair obtained by
applying }
\centerline {a Darboux transformation to the isothermic pair
of planes}
\centerline{\BoxedEPSF{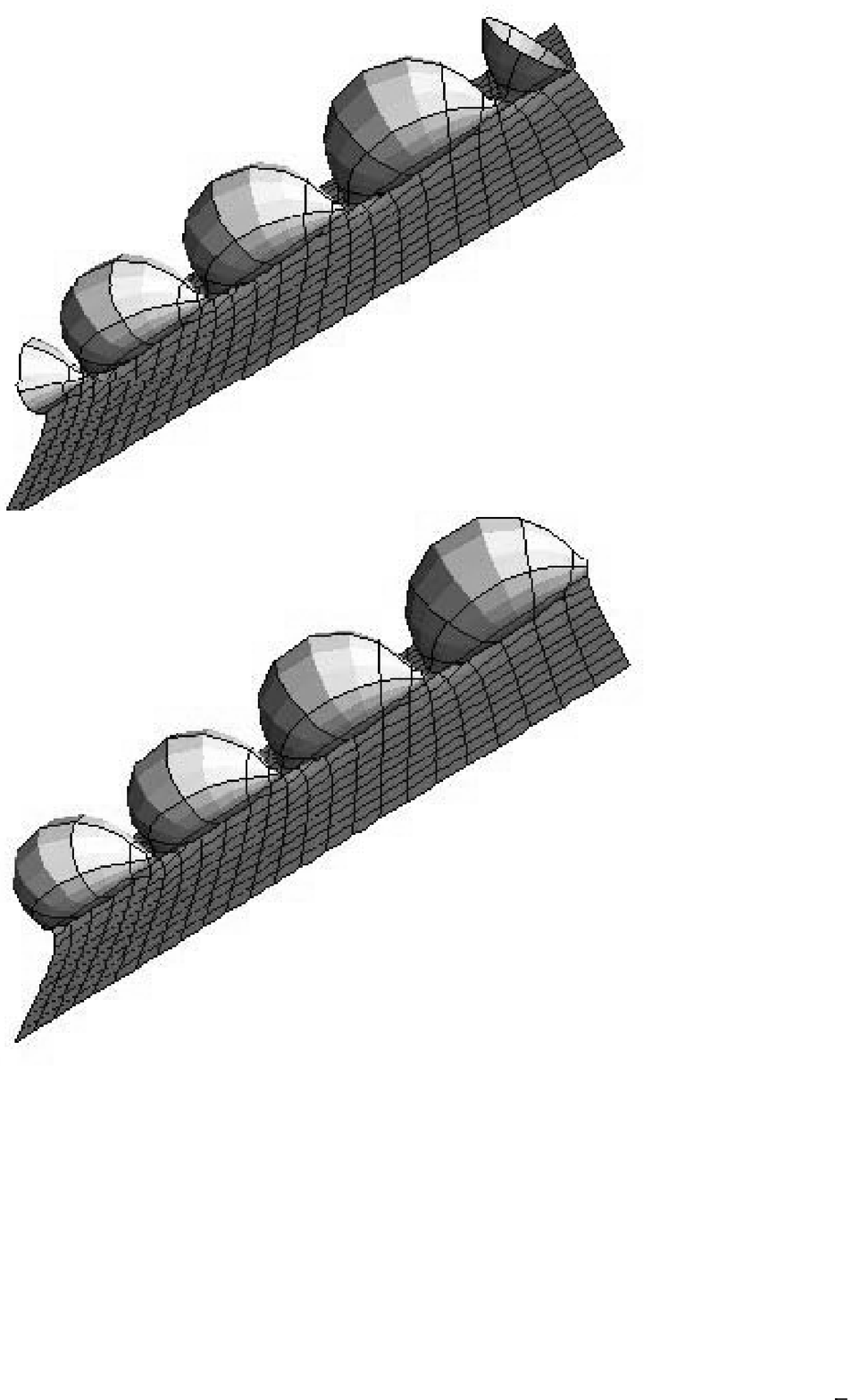 scaled 650}}
\vfil\eject

\centerline{\bf Figure 5}\ms
\centerline{Example \refjo{}. A CMC surface obtained by
applying a Darboux transformation to a cylinder}
\centerline{ Top: surface $\tilde Y_1$ for $a=2, \cosh c=
2$. Bottom: surface
$\tilde Y_1$ for $a=2, \cosh c=4$.}
\centerline{\BoxedEPSF{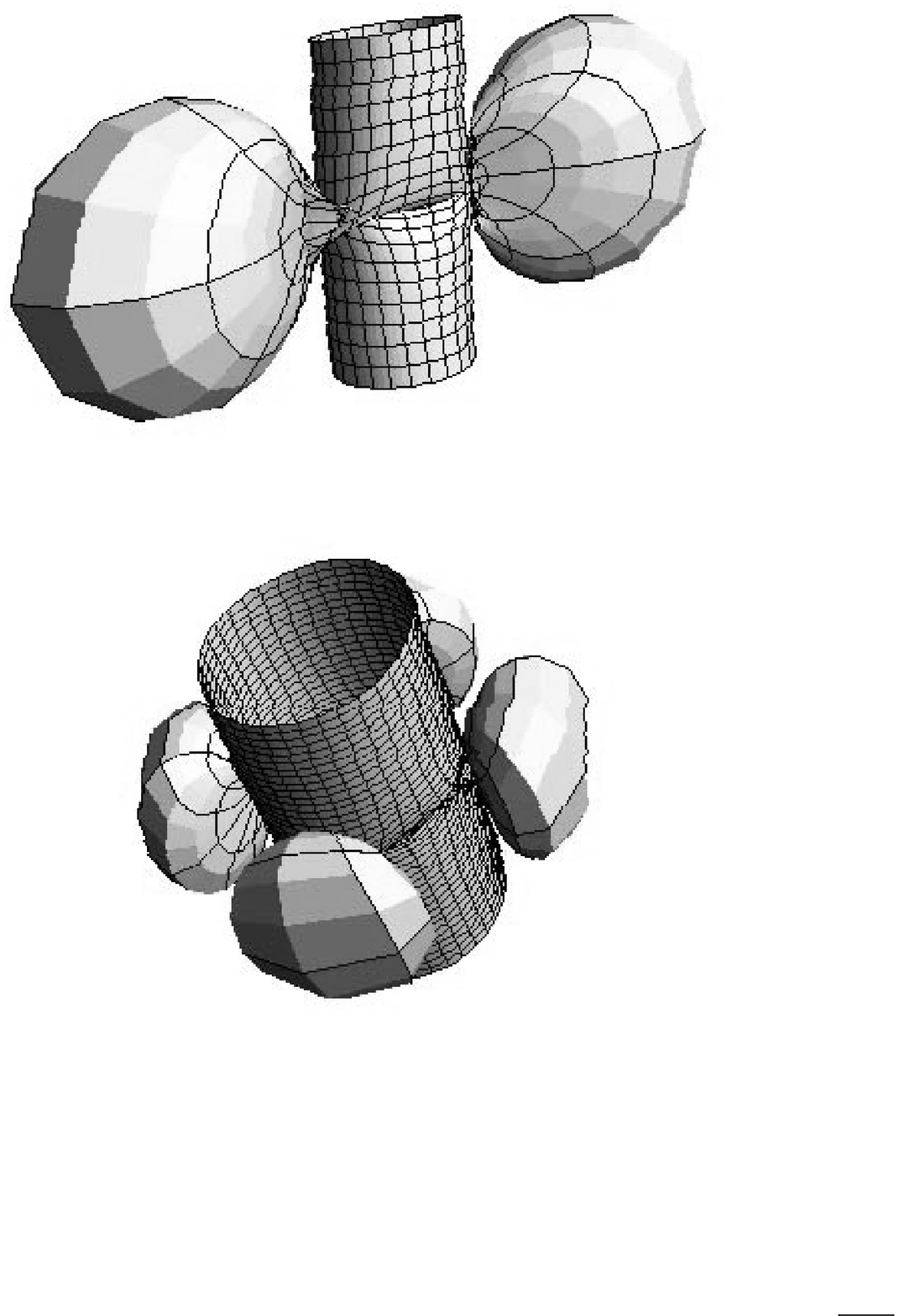 scaled 650}}
\vfil\eject

\centerline{\bf Figure 6}\ms
\centerline{Example \refjo{}. A CMC surface obtained by
applying a Darboux transformation to a cylinder}
\centerline{ Top and bottom: two views of the surface $\tilde
Y_1$ for
$a=2.75,
\cosh c= {9\over 4}$}
\centerline{\BoxedEPSF{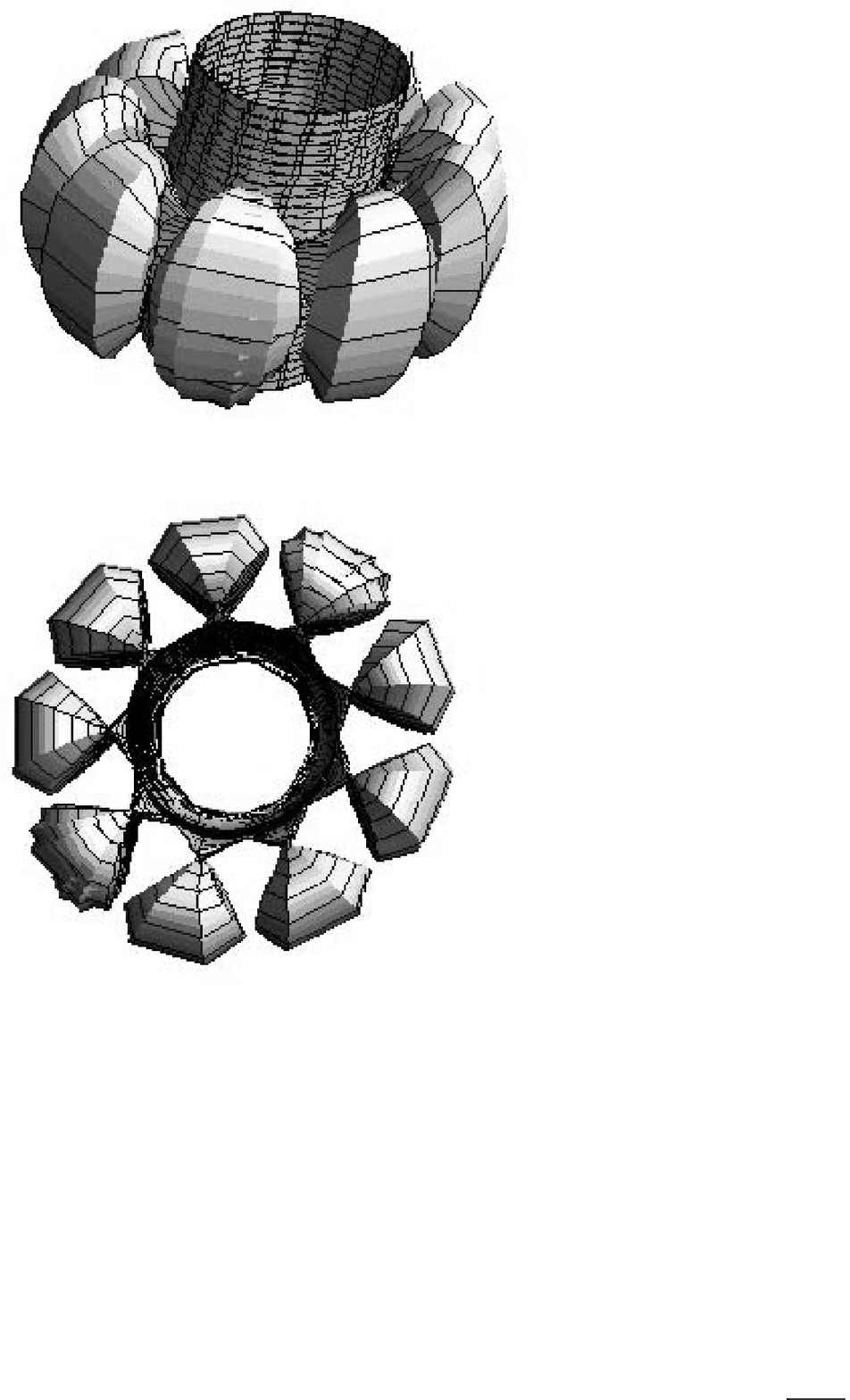 scaled 650}}
\vfil\eject

\centerline{\bf Figure 7}\ss
\centerline{Example \refjo{}. A CMC surface obtained by
applying a Darboux transformation to a cylinder}
\centerline{ Top left and right are two views
of the surface $\tilde Y_1$ for $a=2, \cosh c= {3\over 2}$}
\centerline{Bottom left and right are two views of the surface
$\tilde Y_1$ for $a=2, \cosh c={7\over 5}$.}
\centerline{\BoxedEPSF{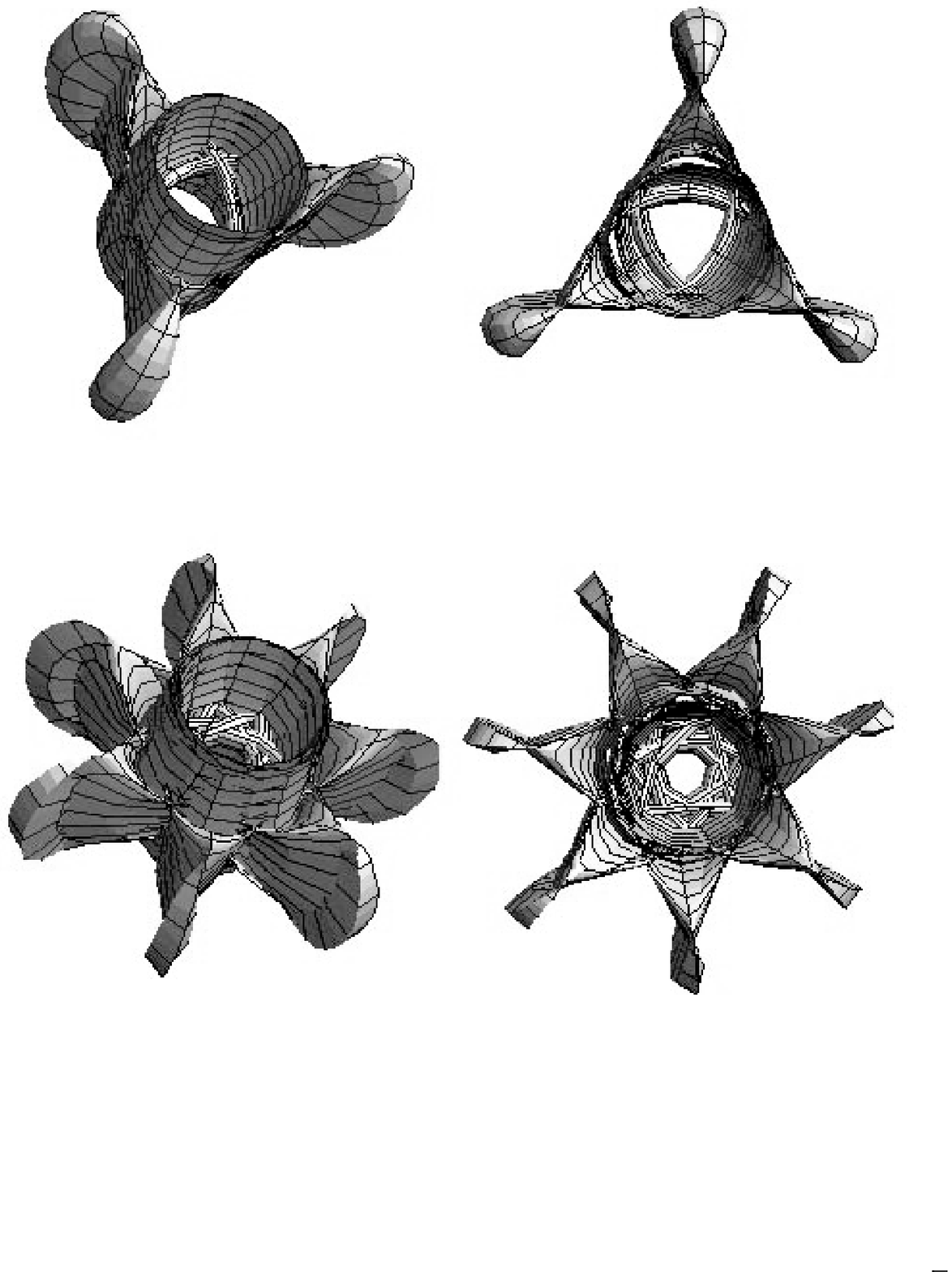 scaled 650}}
\vfil\eject

{\Bibliography

\b //Ba//B\"acklund, A.V.//Concerning surfaces with constant
negative curvature//New Era Printing Co., Lancaster,
PA////1905//////

\a  //BC//Beals, R., Coifman, R.R.//Inverse scattering and evolution
equations//Comm. Pure Appl. Math.//38 //1985//39-90////

\b //Bi1//Bianchi, L.//Leconzioni di Geometria Differenziale//
Bologna////1927//////

\a //Bi2//Bianchi, L.//Ricerce sulle superficie isoterme e sulla
deformazione delle quadriche//Annali Math.
III//11//1905//93-157////

\a//Bo1//Bobenko, A. I.// All constant mean curvautre tori in $R\sp 3, S\sp
3, H\sp 3$ in terms of theta-functions//Math. Ann.//290//1991//209-245////

\a//Bo2//Bobenko, A. I.//Surfaces in terms of 2 by 2 matrices, old and new
integrable cases//Harmonic maps and Integrable
systems, edited by A.P. Fordy and J. Wood////1994//83-127//Vieweg//

\p//Bo3//Bobenko, A.I.//Exploring surfaces through methods for
the theory of integrable systems, lectures on the Bonnet
porblem//////////preprint, math-dg/9909003//

\a//BFPP//Burstall, F., Ferus, D., Pedit, F., Pinkall, U.//Harmonic tori in
symmetric spaces and commuting Hamiltonian systems on loop algebras//Ann. of
Math.//138//1993//173-212////

\a //BHPP//Burstall, F., Hertrich-Jeromin, U., Pedit, F., Pinkall,
U.//Curved flats and isothermic surfaces//Math. Z.
//225//1997//199-209////

\p//Bu//Burstall, F.//Isothermic surfaces: conformal geometry, Clifford
algebras and integrable systems//////////preprint, math-dg/0003096//

\a //Ca// E. Cartan//Sur les vari\'et\'es de courbure constante
d'un espace euclidien ou non-euclidien//Bull. Soc. Math.
France//47//1919//132-208////

\b //Ch//Cherednik, I.V.//Basic methods of soliton theory, Adv. ser. 
Math. Phys. v. 25//World Scientific////1996//////

\a //CGS//Cie\'sli\'nski, J., Goldstein, P., Sym, A.//Isothermic
surfaces in $E^3$ as soliton surfaces//Phys. Lett.
A//205//1995//37-43////

\a //Ci//Cie\'sli\'nski, J.//The Darboux-Bianchi transformation for isothermic
surfaces//Differential Geom. Appl.//7//1997//1-28////

\a//C//Chen, B.//Surfaces with a parallel isoperimetric section//Bull. Amer.
Math. Soc.//79//1973//599-600////

\p //DT//Dajczer, M., Tojeiro, R.//Commuting Codazzi tensors and
the Ribaucour transformation for submanifolds//////////preprint//

\b //Da//G. Darboux//Lecon sur la th\'eorie g\'en\'erale des
surfaces//Chelsea////1972////3rd edition//

\b //FT//Faddeev, L.D., Takhtajan, L.A.//Hamiltonian
methods in the theory of solitons//Springer-Verlag////1987//////

\a //FP1//Ferus, D., Pedit, F.//Curved flats in symmetric
spaces//Manuscripta Math.//91//1996//445-454////

\a //FP2//Ferus, D., Pedit, F.//Isometric immersions of space forms and
soliton theory//Math. Ann.//305//1996//329-342////

\b//H// Helgason, S.//Differential Geometry, Lie Groups, and
Symmetric Spaces//Academic Press////1978//////

\a //HP//Hertrich-Jeromin, U., Pedit, F.//Remarks on the Darboux transform
of isothermic surfaces//Doc. Math.//2//1997//313-333////

\p //HMN//Hertrich-Jeromin, U., Musso, E., Nicolodi, L.//M\"obius
geometry of surfaces of constant mean curvature 1 in hyperbolic
space//////////preprint math/9810157//

\a //M// J. D. Moore//Isometric immersions of space forms in space
forms//Pacific Jour. Math.//40//1972//157-166////

\a//PiS//Pinkall. U., Stering, I.//On the classification of constant mean
curvature tori//Ann. of Math.//130//1989//407-451////

\b //PS//Pressley, A. and Segal, G. B.//Loop Groups//Oxford
Science Publ., Clarendon Press, Oxford////1986//////

\b//PT//Palais, R.S., Terng, C.L.//Critical Point Theory and
Submanifold Geometry//Lecture Notes in
Math.//1353//Springer-verlag//1988////

\a //Sa//Sattinger, D.H.//Hamiltonian hierarchies on semi-simple Lie
algebras//Stud. Appl. Math.//72//1984//65-86////

\a //SZ1//Sattinger, D.H.,Zurkowski, V.D.//Gauge theory of
B\"acklund transformations. I//Dynamics of infinite dimensional
systems, Nato Sci. Inst. Ser. F. Comput. Systems
Sci.//37//1987//273-300//Springer-Verlag//

\a //SZ2//Sattinger, D.H.,Zurkowski, V.D.//Gauge theory of
B\"acklund transformations. II//Physica//26D//1987//225-250////

\a //Sy//Sym, A//Soliton surfaces//lett. Nuovo Cim.//33//1982//394-400////

\a //Ten//Tenenblat, K.//B\"acklund's theorem for submanifolds of
space forms and a generalized wave equation//Boll. Soc. Brasil.
Mat.//16//1985//67-92////

\a //TT//Tenenblat, K., Terng, C.L.//B\"acklund's theorem for
n-dimensional submanifolds of $R^{2n-1}$//Ann.
Math.//111//1980//477-490////

\a //Te1//Terng, C.L.//A higher dimensional generalization of the
sine-Gordon equation and its soliton theory//Ann.
Math.//111//1980//491-510////

\a //Te2//Terng, C.L.//Soliton equations and differential
geometry//J. Differential Geometry//45//1997//407-445////

\a //TU1//Terng, C.L., Uhlenbeck, K.//Poisson actions and
scattering theory for integrable systems//Surveys in Differential
Geometry: Integrable systems (A supplement to J. Differential
Geometry)//4//1998//315-402//preprint dg-ga 9707004//

\a //TU2//Terng, C.L., Uhlenbeck, K.//B\"acklund transformations
and loop group actions//Comm. Pure. Appl. Math.//53//2000//1-75////

\a //TU3//Terng, C.L., Uhlenbeck, K.//Geometry of
solitons//Notice, A.M.S.//47//2000//17-25////

\a //U1//Uhlenbeck, K.//Harmonic maps into Lie group (classical
solutions of the Chiral model)//J. Differential
Geometry//30//1989//1-50////

\a //U2//Uhlenbeck, K.//On the connection between harmonic maps
and the self-dual Yang-Mills and the sine-Gordon
equations//Geometry \& Physics//2//1993//////

\a //ZS//Zakharov, V.E., Shabat, A.B.//Integration of non-linear
equations of mathematical physics by the inverse scattering
method, II//Funct. Anal. Appl.//13//1979//166-174////

\a //Zh//Zhou, Z.X.//Darboux transformations for the twisted
$so(p,q)$-system and local isometric immersions of space
forms//Inverse Problem//14//1998//1353-1370////

}

\vfil\eject
\noindent Authors' addresses:

\ss
\ni {(1)} Martina Br\"uck, Institut f\"ur Mathematik, Universit\"at Augsburg,
Universit\"at str. 8, D-86135 Augsburg, Germany, email:mbrueck@MI.Uni-Koeln.DE

\ni {(2)} Xi Du, Department of Mathematics, Northeastern Univ., Boston, MA
02115

\ni {(3)} Joonsang Park, Department of Mathematics, Dongguk
Univ., Seoul, 100-715 Korea, email:jpark@cakra.dongguk.ac.kr

\ni {(4)} Chuu-Lian Terng,  Department of Mathematics, Northeastern
Univ., Boston, MA 02115, email: terng@neu.edu

\end